\documentclass{article}

\usepackage[english]{babel}

\usepackage[letterpaper,top=2cm,bottom=2cm,left=3cm,right=3cm,marginparwidth=1.75cm]{geometry}

\usepackage{amsmath}
\usepackage{amssymb}
\usepackage{graphicx}
\usepackage[colorlinks=true, allcolors=blue]{hyperref}
\usepackage{amsthm}
\usepackage{tikz}
\usepackage{caption}
\usepackage{subcaption}
\usepackage{float}
\usepackage{svg}
\usepackage{mathtools}
\usetikzlibrary{math}
\usetikzlibrary{shapes}
\usepackage{stmaryrd}
\usepackage{hyperref}
\usepackage{comment}

\usepackage{tikz}
\usetikzlibrary{arrows.meta,calc,positioning,decorations.pathreplacing}

\tikzset{
  >={Latex[length=3mm]},
  mybox/.style={draw=black, thick, rounded corners=2pt},
  support/.style={draw=black!65, dashed, line width=0.9pt},
  neggrad/.style={-Latex, very thick, blue!70!black},
  rawgrad/.style={-Latex, very thick, red!75!black, dashed},
  projgrad/.style={-Latex, very thick, red!75!black},
  cone/.style={fill=red!10, draw=black!45},
  smalllabel/.style={font=\small},
  titlelabel/.style={font=\small\bfseries}
}

\newcommand{\cve}[1]{\textcolor{blue}{[\texttt{#1}]}}
\newcommand{\ve}[1]{\textcolor{blue}{#1}}

\DeclareMathOperator{\barron}{\mathcal{B}}
\DeclareMathOperator{\Prob}{\mathcal{P}}
\DeclareMathOperator{\F}{\mathcal{F}} 
\DeclareMathOperator{\R}{\mathbb{R}} 
\DeclareMathOperator{\N}{\mathbb{N}} 
\DeclareMathOperator{\E}{\mathcal{E}}
\DeclareMathOperator{\C}{\mathcal{C}}

\newtheorem{definition}{Definition}
\newtheorem{theorem}{Theorem}
\newtheorem{proposition}{Proposition}
\newtheorem{lemma}{Lemma}
\newtheorem{remark}{Remark}
\newtheorem{corollary}{Corollary}

\newtheorem{hypothesis}{Hypothesis}

\everymath{\displaystyle}

\title{Two-layer neural networks for Schrödinger eigenvalue problems}
\author{Dus Mathias$^{*\dag}$, Ehrlacher Virginie\footnote{CERMICS, École nationale des Ponts et Chaussées, IPP, CNRS, Marne-la-Vallée, France} \footnote{INRIA Paris, France, PARMA team-project, \newline \indent \indent emails: mathias.dus@enpc.fr, virginie.ehrlacher@enpc.fr} 
}

\begin{document}

\maketitle

\begin{abstract}
We study an infinite-width two-layer neural-network approach for Schrödinger eigenvalue problems with Neumann boundary conditions on the unit cube in arbitrary dimension. The ground-state problem is recast as the minimization of the Schrödinger energy under an L2-normalization constraint. Through a Barron-type representation, the neural network is described by a probability measure on the parameter space, and the resulting constrained optimization problem is analyzed as a Wasserstein gradient dynamics.

The normalization constraint substantially changes the structure of the flow: admissible directions involve a nonlocal Lagrange multiplier, and the evolution is governed by a constrained nonlinear transport equation in parameter space. We establish existence and uniqueness of global measure-valued  solutions by a local fixed-point construction for the characteristic
flow, combined with energy and moment estimates that prevent finite-time escape of the support. We then prove that any trajectory converging in 2-Wasserstein distance has, under a suitable density assumption on the feature class, a limit whose represented function is an eigenfunction of the Schrödinger operator. The result is therefore a conditional identification of convergent equilibria, rather than a global convergence theorem toward the ground state.

Numerical experiments on smooth benchmark potentials illustrate the empirical
behavior of a finite-width stochastic training scheme inspired by the analyzed constrained Wasserstein dynamics.
\end{abstract}


\section{Introduction}
The main interest of this work is to introduce a rigorous framework for solving eigenvalue problems of a certain type of elliptic operators \textit{i.e.} Schrödinger operators, using neural networks.
In recent years, the use of deep neural networks to solve partial differential equations has sparked massive interest, giving rise to highly popular methodologies such as Physics-Informed Neural Networks (PINNs) \cite{raissi2019physics} and the Deep Galerkin Method \cite{sirignano2018dgm}.
A first step towards this objective came from ideas introduced in~\cite{BachChizat2018, DusEhrlacher2023} where the authors propose a numerical method to solve a different but close problem:

\begin{equation}\label{PoissonEquationIntro}
\left\{
\begin{aligned}
- \Delta u^\star &= f \text{ on } \Omega, \\
 \partial_n u^\star &= 0 \text{ on } \partial \Omega
\end{aligned}
\right.
\end{equation}
where $f$ is a source term. Under some conditions on $f \in L^2(\Omega)$, it is possible to approximate numerically the minimizer $u^\star$ as a convex combination of what we call feature functions. In the same article, it is shown that minimizing a certain variational energy $\E_\tau$, it is possible to solve the Poisson-Neumann problem efficiently when the source term $f$ belongs to the Barron space $\barron^0(\Omega)$. To do so, the authors propose to analyze the convergence of the gradient curve associated with $\E_\tau$. More precisely, they prove that if the gradient curve converges, then it is necessarily towards an optimum. Some numerical experiments show the potential of the method on simple canonical examples.

For the Schrödinger eigenvalue problem, things are a bit different. It consists of finding one/the ground eigencouple of the elliptic operator $- \Delta + W$ where $W$ is a potential, supplemented with boundary condition that we do not specify here. In \cite{LuLu2022}, the authors adapted the work \cite{Jianfeng2021} to give conditions on the potential $W$ such that the ground eigenfunction belongs to the Barron space of order $s>0$ denoted $\barron^s$. More specifically, they prove that if the potential $W$ belongs to $\barron^s$  and if the ground eigenvalue is negative, then the eigenfunction $u^\star$ belongs to $\barron^s$. A two-layer neural network approximation is then particularly relevant when $s \geq 2$, as the Barron space naturally corresponds to the set of functions represented by infinite-width two-layer neural networks, which can be elegantly analyzed through the lens of probability measures and Wasserstein gradient flows.

A major goal of studying the Schrödinger problem is to get the ground state of the multi-body Schrödinger operator given by:

\begin{equation}
        H :=   -\frac 1 2 \sum_{n=1}^N \Delta_{x_n} + \sum_{n=1}^N \sum_{m=1}^M \frac{-Z_m}{|x_n - X_m|} + \sum_{1 \leq j < n \leq N} \frac{1}{|x_j - x_n|}
\end{equation}
where:

\begin{itemize}
\item $(x_n)_{1\leq n \leq N}$ is the position of $N$ electrons,
\item $(X_m)_{1\leq m \leq M}$ is the position of $M$ nuclei of charge $(Z_m)_{1\leq m \leq M}$
\item The form domain of $H$ is $H^1(\R^{dN})$.
\end{itemize}

The difficulties posed by the Coulombic singularities, the unbounded domain and the antisymmetry condition (due to spins) on the domain of definition of $H$ are far beyond the scope of this paper. The strong asymmetry and singularities present in the full many-body problem hinder a direct two-layer neural network approximation. Indeed, while the Jastrow factor ansatz shows that the ground eigenfunction is highly regular in the $W^{2, \infty}$ sense \cite{Fournais2005}, its lack of Barron regularity \cite{yserentant2025} rules out the use of this architecture in that specific context. Nevertheless, the mixed derivative regularity described in \cite{Yserentant2004} provides evidence for a deep neural approximation as described in \cite{MontanelliDu2018, GriebelHamaekers2007, Smo63, Kor63}. Furthermore, the experimental results presented in references \cite{Pfau_2020, Hermann_2020}, building upon the seminal introduction of Neural Quantum States for many-body systems \cite{carleo2017solving}, suggest that, at least for simple systems, the eigenfunction can be represented by complex neural networks, for which a clearly defined mathematical framework is lacking. Consequently, this work is a theoretical starting point that can be applied directly to simpler problems, such as Henon-Heiles's systems in cosmology.

We will follow and adapt the strategy developed in~\cite{BachChizat2018, DusEhrlacher2023}, but the present setting
introduces several substantial additional difficulties. In contrast with the unconstrained
framework considered by Chizat and Bach, the Schr\"odinger ground-state problem is posed
under the normalization constraint
\[
\|u\|_{L^2(\Omega)}=1,
\]
which changes the geometry of the problem at the level of measures. The relevant
dynamics is therefore not the plain Wasserstein gradient flow of the lifted energy, but a
constrained gradient dynamics evolving on a nonlinear constraint manifold in probability
space.

Our first contribution is to formulate this constrained infinite-width neural-network
problem at the level of probability measures and to identify the associated constrained
local slope. The main novelty here is that admissible perturbations cannot be obtained by
following only the negative energy gradient. In order to remain on the constraint set, one
must combine a displacement along the energy gradient with a correction along the
constraint gradient. 
This leads to a constrained dynamics. In the interior of the truncated parameter
domain, the correction is given by a nonlocal Lagrange multiplier. The truncation is
used as a compactness device in the existence proof; on each finite time interval, the
truncation is chosen large enough so that the corresponding trajectory remains strictly
inside the truncated domain.

Our second contribution is a global well-posedness theory for the constrained measure-valued dynamics. Local solutions are constructed by a fixed-point argument at the level of the characteristic flow. Energy, second-moment, and characteristic estimates then prevent finite-time blow-up and yield a globally defined solution with bounded support on every finite time interval. The truncated minimizing  movement remains useful for identifying the interior metric slope and the corresponding velocity field, but is not used to propagate support.

Our third contribution concerns the asymptotic behavior of the flow. Assuming that the
trajectory converges in Wasserstein distance, we show that its limit is a critical point of
the constrained problem and, under the density assumption on the feature space, that the
represented function is an eigenfunction of the Schr\"odinger operator. The proof of this
identification result combines the structure of the reduced potential
$V_{\tau,\mu}^C$, support propagation along the transported slice, and a dynamical
exclusion mechanism for noncritical limits. In particular, the analysis of the limiting
behavior differs markedly from the one in the unconstrained setting of Chizat--Bach:
here the limit is not characterized as the minimizer of an unconstrained training
functional, but through the Euler--Lagrange condition associated with a normalized
eigenvalue problem.

The paper is organized as follows. In Section~\ref{sec:notation}, we introduce the Schr\"odinger
eigenvalue problem, the neural-network parametrization, and the approximation results
used throughout the paper. Section~\ref{sec:mainres} is devoted to statement of the main theoretical results proved in the present work: the existence of a constrained Wasserstein gradient curve and its convergence towards a solution of the Schr\"odinger eigenvalue problem. In Section~\ref{proofTheorem3}, we prove the existence of a constrained Wasserstein gradient curve: we define the lifted constrained energy, identify its local slope, and prove the
existence of global measure-valued solutions.  In Section~\ref{sec:proofConvergence}, we study the asymptotic
behavior of the flow and prove the identification of convergent limits as critical points,
and then as eigenfunctions. Section~\ref{sec:sim} presents numerical experiments illustrating the proposed approach.

\medskip
\underline{\textbf{Notation:}} The notation $C$ is used for universal constants. If there exists a dependence on some parameter $\iota$, the constant is denoted $C(\iota)$. 

\section{Preliminaries and notation}\label{sec:notation}

The aim of this section is to present the eigenvalue problem considered in this work and fix some notation.

\subsection{Problem setting}

\medskip

Let $\Omega := [0,1]^d$ be the $d$-dimensional hypercube and $\partial \Omega$ denote the boundary of $\Omega$. We first introduce the quadratic functional $\mathcal E: H^1(\Omega) \to \mathbb{R}$ defined as follows:
$$
\forall u \in H^1(\Omega), \quad \mathcal{E}(u):=  \frac 12 \int_\Omega |\nabla u|^2 + W u^2dx
$$
where $W\in L^\infty(\Omega; \mathbb{R})$. The functional $\mathcal E$ is then the quadratic form associated to the Schrödinger operator $-\Delta + W$ with Neumann boundary condition, which is a self-adjoint operator with compact resolvent on $L^2(\Omega)$ with form domain $H^1(\Omega)$ and operator domain $\{v \in H^2(\Omega): \; \partial_n v=0 \mbox{ on } \partial \Omega\}$.  We are interested in the computation of a solution $u^\star\in H^1(\Omega)$ of the following minimization problem
\begin{equation}\label{eq:Rayleigh}
u^\star \in \mathop{\rm argmin}_{u\in H^1(\Omega), \; \|u\|_{L^2(\Omega)} = 1} \mathcal E(u),
\end{equation}
solution of the Schrödinger eigenvalue problem with Neumann boundary condition: find $(u^\star, \lambda)\in H^1(\Omega) \times \mathbb{R}$ solution to:
\begin{equation}\label{PoissonEquation}
\left\{
\begin{aligned}
- \Delta u^\star + W u^\star &= \lambda u^{\star} \text{ on } \Omega, \\
 \partial_n u^\star &= 0 \text{ on } \partial \Omega.
\end{aligned}
\right.
\end{equation}
Moreover, since $u^\star$ is a minimizer of~\eqref{eq:Rayleigh}, it holds that $\lambda$ is the lowest eigenvalue of the operator $-\Delta + W$. Formulating the eigenvalue problem as a variational energy minimization task solved via neural network ansatz is deeply connected to the Deep Ritz Method \cite{e2018deep}, which tackles variational problems by minimizing energy functionals over neural network parameter spaces. For convenience, we use the following notation for the constraint:

$$
\forall u \in L^2(\Omega), \quad \mathcal{C}(u):=  \sqrt{\int_\Omega |u|^2 dx} - 1.
$$

%In all of the paper the manifold $\Theta$ is the form:
%
%$$
%\Theta:= B_{\R^3}(0, \sqrt{3}) \times S^{d-1}
%$$
%where $r>0$ will be chosen later. An element of $\Theta$ is denoted $\theta:= (c,a,b,w) \in \R^3 \times S^{d-1}$.

\medskip

The aim of the present paper is to propose and analyze a neural-network based numerical method for the resolution of (\ref{eq:Rayleigh}) in the spirit of~\cite{DusEhrlacher2023, BachChizat2018}, based on the use of inifinite-width two-layer neural networks.

\medskip

\subsection{Activation function}\label{sec:activation}

We introduce here the particular choice of activation function we consider in this work, which is the same as the one used in~\cite{DusEhrlacher2023}. Let $\sigma: \mathbb{R} \rightarrow \mathbb{R} $ be the classical Rectified Linear Unit (ReLU) function so that for all $y\in \mathbb{R}$, $\sigma(y) := \max(y,0)$. Let $\rho: \mathbb{R} \to \mathbb{R}$ be defined by 
\begin{equation}
\left\{
\begin{array}{cl}
Z \exp\left(- \frac{\tan(\frac{\pi}{2} y)^2}{2}\right) & \text{if } |y| \leq 1 \\
0 & \text{otherwise},
\end{array}
\right.
\end{equation}
where the constant $Z\in \mathbb{R}$ is defined such that the integral of $\rho$ is equal to one. For all $\tau >0$, we then define $\rho_\tau := \tau\rho(\tau \cdot)$ and $\sigma_\tau: \mathbb{R} \to \mathbb{R}$ the regularized ReLU function defined by
\begin{equation}\label{definitionSoftplus}
\forall y \in \mathbb{R}, \ \sigma_\tau(y) := (\rho_\tau \star \sigma)(y). 
\end{equation}

In this work, we will rather use a hat version of the regularized ReLU activation function. More precisely, we define:
\begin{equation}\label{definitionRegularTrelu}
\forall y \in \mathbb{R}, \ \sigma_{H,\tau}(y) := \sigma_\tau(y+1) - \sigma_\tau(2y) + \sigma_\tau(y-1),
\end{equation}
which is called hereafter the regularized HReLU (Hat ReLU) activation. When $\tau = +\infty$, the following notation is proposed:

\begin{equation}\label{definitionTrelu}
\forall y \in \mathbb{R}, \ \sigma_{H}(y) := \sigma(y+1) - \sigma(2y) + \sigma(y-1).
\end{equation}
Note that a direct corollary of \cite[Lemma~2]{DusEhrlacher2023} is that there exists a constant $C>0$ such that for all $\tau >0$, 
\begin{equation}\label{eq:sigmaerr}
\|\sigma_H - \sigma_{H,\tau}\|_{H^1(\mathbb{R})} \leq \frac{C}{\sqrt{\tau}}.
\end{equation}
We will also use the fact that there exists a constant $C>0$ such that for all $\tau>0$, 
\begin{equation}\label{eq:bound}
\|\sigma_{H,\tau}\|_{L^\infty(\mathbb{R})}\leq C, \; \|\sigma'_{H,\tau}\|_{L^\infty(\mathbb{R})}\leq C,\; \|\sigma''_{H,\tau}\|_{L^\infty(\mathbb{R})}\leq C\tau \; \mbox{ and } \|\sigma'''_{H,\tau}\|_{L^\infty(\mathbb{R})}\leq C\tau^2. 
\end{equation}

\begin{remark} Actually, we would like to stress here on the fact that all the theoretical results proved here remain true under the more general assumption that the activation function is $C^\infty$ with compact support.
\end{remark}

From now on, in all the article, we will assume that $\tau \geq 1$.

\subsection{Infinite width two-layer neural networks}

Let us introduce the set of parameter values of the neural network
$$
\Theta := \R \times S^{d-1} \times \R,
$$ 
where $S^{d-1}$ denotes the unit sphere embedded in $\R^d$. Moreover, let $\mathcal P_2(\Theta)$ be the set of probability measures on $\Theta$ with finite second-order moments. The space $\mathcal P_2(\Theta)$ is equipped with the 2-Wasserstein distance:
$$
\forall \mu, \nu \in \mathcal P_2(\Theta), \quad W_2^2(\mu, \nu):= \mathop{\inf}_{\gamma \in \Gamma(\mu,\nu)} \int_{\Theta^2} d(\theta, \widetilde{\theta})^2\,d\gamma(\theta, \widetilde{\theta}),
$$
where $\Gamma(\mu,\nu)$ is the set of probability measures on $\Theta^2$ with marginals given respectively by $\mu$ and $\nu$ and where $d$ is the geodesic distance in $\Theta$. The geodesic distance between $\theta :=(a,w,b) , \tilde{\theta}:=(\tilde{a}, \tilde{w}, \tilde{b}) \in \Theta$ is defined as:
$$
d(\theta, \tilde{\theta}) = \sqrt{ (a - \tilde{a})^2 + d_{S^{d-1}}(w,\tilde{w})^2 + (b - \tilde{b})^2},
$$
with $d_{S^{d-1}}$ the geodesic distance on the unit sphere $S^{d-1}$. We also denote by $\Gamma_o(\mu,\nu)$ the subset of probability measures $\gamma \in \Gamma(\mu,\nu)$ such that $W_2(\mu,\nu) = \sqrt{\int_{\Theta^2} d(\theta, \tilde{\theta})^2\,d\gamma(\theta, \tilde{\theta})}$. If $\gamma \not\in \Gamma_o(\mu,\nu)$, we use the notation $W_{2,\gamma}(\mu,\nu) := \sqrt{\int_{\Theta^2} d(\theta, \tilde{\theta})^2\,d\gamma(\theta, \tilde{\theta})}$. 

For any $\tau >0$, we introduce the function $\Phi_{\tau} : \Theta \times \Omega \to \mathbb{R}$ defined such that
\begin{equation}
\label{definitionPhi}
\forall \theta:=(a,w,b)\in \Theta, \; \forall x\in \Omega, \quad \Phi_\tau(\theta; x) := a \sigma_{H,\tau}(w \cdot x + b),
\end{equation}
and the operator $P_\tau : \Prob_2(\Theta) \rightarrow H^1(\Omega)$ defined by:
$$
\forall \mu \in \Prob_2(\Theta), \ P_\tau \mu := \int_\Theta \Phi_\tau(\theta; \cdot) d\mu(\theta).
$$

\medskip

The aim of this paper is to analyze a numerical method which consists in approximating a solution $u^*\in H^1(\Omega)$ to (\ref{eq:Rayleigh}) by a function of the form $P_\tau \mu^\star$ where $\mu^\star$ is a solution of a constrained minimization problem of the form
\begin{equation}\label{eq:constrained}
\begin{array}{rl}
\inf_{ \substack{ \mu \in \Prob_2(\Theta) \\ \C_\tau(\mu) = 0} } \E_\tau(\mu).
\end{array}
\end{equation}
where:
\begin{itemize}
\item the energy writes:
\begin{equation}\label{definitionEtau}
\E_\tau :
\left\{
\begin{array}{rl}
\Prob_2(\Theta) \mapsto& \R \\
\mu \rightarrow& \E(P_\tau \mu),
\end{array}
\right.
\end{equation}

\item the constraint cost is:

$$
\C_\tau :
\left\{
\begin{array}{rl}
\Prob_2(\Theta) \mapsto& \R \\
\mu \rightarrow& \| P_\tau \mu \|_{L^2(\Omega)} - 1.
\end{array}
\right.
$$ 
The associated set of constraint is denoted by $X_\tau := \{ \mu \in \Prob_2(\Theta) \ | \ \C_\tau(\mu) = 0 \}$.
\end{itemize} 
More precisely, we prove the existence of a constrained gradient curve $(\mu_t)_{t\geq 0}$ associated to (a regularized version of) the constrained minimization problem~\eqref{eq:constrained}. Then, under the assumption that such a curve converges as $t$ goes to infinity to some measure $\mu_\infty\in \mathcal P_2(\Theta)$, we prove that $P_\tau(\mu_\infty)$ is an eigenfunction of the Schrödinger operator $-\Delta + W$. We would like to stress on the fact that this eigenfunction may however be associated to an eigenvalue $\lambda_\infty$ which may not be the smallest eigenvalue of the operator $-\Delta + W$ in general.

\subsection{Barron spaces and regularity of eigenfunctions of Schrödinger operators}

In this section, we introduce the Barron space to justify our choice of neural network approximation of ground eigenfunctions. To do so, let us introduce the orthonormal basis in $L^2(\Omega)$ composed of the eigenfunctions $\{\phi_k\}_{k\in \mathbb{N}^d}$ of the Laplacian operator with Neumann boundary condition, where
\begin{equation}\label{baseL2}
\forall k = (k_1,\ldots, k_d)\in \N^d, \; \forall x:=(x_1,\cdots, x_d)\in \Omega, \quad  \phi_k(x_1, \ldots,x_d):= \prod_{i=1}^d c_{k_i}\cos(\pi k_i x_i),
\end{equation}
where $c_0 = 1$ and $c_\ell = \sqrt{2}$ for all $\ell \in \mathbb{N}\setminus\{0\}$. 

Notice that $\{\phi_k\}_{k\in \N^d}$ is also an orthogonal basis of $H^1(\Omega)$. Using this basis, we have the Fourier representation formula for any function $u \in L^2(\Omega)$:
$$
u = \sum_{k \in \N^d} \hat{u}(k) \phi_k,
$$
where for all $k\in \mathbb{N}^d$, $ \hat{u}(k):= \langle \phi_k, u \rangle_{L^2(\Omega)}$. This allows to define the (spectral) Barron space~\cite{Jianfeng2021} as follows:

\begin{definition}

For all $s>0$, the Barron space $\barron^s(\Omega)$ is defined as:

\begin{equation}
\label{def:DefinitionBarron}
\barron^s(\Omega) := \Big\{ u \in L^1(\Omega) : \sum_{k \in \N^d}(1 + \pi^s |k|^s_1)|\hat{u}(k)| < +\infty \Big\}
\end{equation}
and the space $\barron^2(\Omega)$ is denoted $\barron(\Omega)$. Moreover, the space $\barron^s(\Omega)$ is embedded with the norm:
\begin{equation}
\label{normBarron}
\| u\|_{\barron^s(\Omega)} := \sum_{k \in \N^d}(1 + \pi^s |k|^s_1)|\hat{u}(k)|.
\end{equation}

\end{definition}

One of the reason we consider problem \eqref{eq:constrained} is that under some regularity hypothesis on the potential, the ground eigenfunction is unique and belongs to some Barron space. More precisely, using Krein-Rutman theory, the authors of \cite{LuLu2022} proved the following result.

\begin{theorem}[\cite{LuLu2022}]\label{th:resultLuLu}
Let $s>0$. If $W \in L^\infty(\Omega) \ve{\cap} \barron^s(\Omega)$ is such that $(-\Delta + W, H^1(\Omega))$  has a spectral gap, then the eigenspace associated to the lowest eigenvalue of $-\Delta +W$ is one-dimensional; choosing the positive $L^2$-normalized representative yields a unique ground eigenfunction $u^*$ which belongs to $\barron^s(\Omega)$.
\end{theorem}

By the theorem of approximation of Barron function, it is then possible to approximate the ground state $u^\star$ by a two-layer neural network. To present such result, we introduce the concept of feature space $\F_{\chi,m}(B)$ defined as:
\begin{equation}
\label{eq:defFrelu}
\F_{\chi,m}(B) := \left\{ \sum_{i=1}^m a_i \chi(w_i \cdot x + b_i) : a_i, b_i\in \R, \, w_i \in \R^d, \;  |w_i| = 1, |b_i| \leq 1, \sum_{i = 1}^m |a_i| \leq 4B \right\}
\end{equation}
where $\chi : \R \rightarrow \R$ be measurable, $m\in \N^*$ and $B>0$.

\begin{theorem}[\cite{Jianfeng2021}]
\label{thapproximationTheorem}

For any $u \in \barron(\Omega)$, $m \in \N^*$ :

\begin{itemize}
\item[(i)] there exists $u_m \in \F_{\sigma_H,m}(\| u \|_{\barron(\Omega)})$ such that:
$$
\| u - u_m \|_{H^1(\Omega)} \leq  \frac{C \| u \|_{\barron(\Omega)}}{\sqrt{m}},
$$

\item[(ii)] there exists $\tilde{u}_m \in \F_{\sigma_{H,m},m}(\| u \|_{\barron(\Omega)})$ such that:
\begin{equation}
\| u - \tilde{u}_m \|_{H^1(\Omega)} \leq  \frac{C \| u \|_{\barron(\Omega)}}{\sqrt{m}}
\end{equation}

\end{itemize}
where for both items, $C$ is a universal constant which depends neither on $d$ nor on $u$.
\end{theorem}

\section{Main theoretical results}\label{sec:mainres}

We gather in this section our main theoretical results, namely the existence of a solution to a Wasserstein constrained gradient curve (Section~\ref{sec:existence}), and its convergence to a solution of the Schrödinger eigenvalue problem (Section~\ref{sec:convergence}). In all the rest of the document, we assume that $d\geq 2$ is fixed. 

\subsection{Existence of a constrained gradient curve}\label{sec:existence}

To prove the existence of a gradient curve, we first prove the existence of a gradient curve on compact domains which is stated in the following proposition. More precisely, for $r>0$, let us introduce $K_r:= [-r,r] \times S^{d-1} \times [-r, r]$. We introduce the penalized energy $\E^C_{\tau,r} : \Prob_2(\Theta) \mapsto \R \cup \{ +\infty\}$ which is defined by:
$$
\forall \mu \in \mathcal P_2(\Theta), \quad \E^C_{\tau,r}(\mu) := 
\left\{
\begin{array}{rl}
\E_{\tau}(\mu) & \text{ if } \C_\tau(\mu) = 0 \text{ and Supp}(\mu) \subset K_r, \\
+\infty & \text{otherwise,}
\end{array}
\right.
$$ 
where ${\rm Supp}(\mu)$ denotes the support of the probability measure $\mu$. Let us also introduce the notation:
$$
X_{\tau} := \{ \mu \in \Prob_2(\Theta) \ | \ \C_\tau(\mu) = 0 \} \quad \mbox{ and } \quad X_{\tau,r} := \{ \mu \in \Prob_2(\Theta) \ | \ \C_\tau(\mu) = 0 \text{ and Supp}(\mu) \subset K_r \}.
$$

In the following, we use the standard metric local slope for any ${\cal F}: \mathcal P_2(\Theta) \to \R \cup \{+\infty\}$ defined as
\begin{equation}\label{eq:localSlope}
  |\partial\mathcal F|(\mu)
  :=\limsup_{W_2(\nu, \mu)\to 0}
  \frac{(\mathcal F(\mu)-\mathcal F(\nu))^+}{W_2(\mu,\nu)}
\end{equation}
and the relaxed local slope \(|\partial\mathcal F|_{\rm rel}\), defined as the lower
semicontinuous envelope of \(|\partial\mathcal F|\) with respect to joint convergence of
\(\mu\) and \(\mathcal F(\mu)\). More precisely, we have
\begin{equation}\label{eq:relaxedSlope}
|\partial {\cal F}|_{\mathrm{rel}}(\mu)
:=
\inf
\left\{
\liminf_{n\to\infty} |\partial {\cal F}|(\mu_n)
:
\mu_n\to \mu \mbox{ in }W_2,\ 
{\cal F}(\mu_n)\to {\cal F}(\mu)
\right\}.
\end{equation}

It then holds:
\begin{proposition}[Existence of a truncated generalized minimizing movement]
\label{prop:existenceCurveSupportCond}
Let $r>0$ and $\mu_0\in X_{\tau,r}$. Then, there exists a curve $\mu^r\in AC^2_{\mathrm{loc}}\bigl([0,+\infty);\mathcal P_2(K_r)\bigr)$ such that $\mu^r_0=\mu_0$,
which is a generalized minimizing movement for ${\cal E}^C_{\tau,r}$ and satisfies
$\mu^r_t\in X_{\tau,r}$ for every $t\geq0$. Moreover, for every
$0\leq s\leq t<+\infty$,
\begin{equation}
\label{eq:truncated-edi}
 {\cal E}_\tau(\mu^r_t)
 +\frac12\int_s^t |\dot\mu^r_q|^2\,dq
 +\frac12\int_s^t
   |\partial {\cal E}^C_{\tau,r}|_{\mathrm{rel}}^2(\mu^r_q)\,dq
 \leq {\cal E}_\tau(\mu^r_s).
\end{equation}
For a.e. $t>0$, there exists a Borel vector field
$v^r_t\in L^2(\Theta,\mu^r_t;T\Theta)$ such that
\[
  \partial_t\mu^r_t+\operatorname{div}(v^r_t\mu^r_t)=0
\]
in the sense of distributions and
\[
  |\dot\mu^r_t|^2
  =\int_\Theta |v^r_t(\theta)|^2\,d\mu^r_t(\theta) \qquad \mbox{for a.e. }t>0.
\]

\end{proposition}

\begin{remark}
Proposition~\ref{prop:existenceCurveSupportCond} asserts neither uniqueness of the
truncated generalized minimizing movement nor that
$|\partial E^C_{\tau,r}|_{\mathrm{rel}}$ is a strong upper gradient; in
particular, it does not claim that $\mu^r$ is a curve of maximal slope.  The
integrated dissipation inequality \eqref{eq:truncated-edi}, which is the
estimate needed below, is obtained directly by passing to the limit in the
minimizing-movement scheme.  Uniqueness of the untruncated bounded-support
flow will follow instead from the stability estimate proved later.
\end{remark}

The previous proposition does not identify the velocity field \(v_t^r\). Such an
identification will only be needed in the regime where the support of \(\mu_t^r\) remains
strictly inside \(K_r\). This is precisely the content of Proposition~\ref{prop:propositionCharacteristics} below. Actually, in this case, the velocity field $v_t^r$ can be interpreted using characteristics as shown by Proposition~\ref{prop:propositionCharacteristics}. More precisely, in the following, we denote by $T\Theta$ the tangent bundle of $\Theta$, i.e. 
$$
T\Theta := \bigcup_{\theta\in \Theta} \{\theta\} \times T_\theta \Theta, 
$$
where $T_\theta \Theta$ is the tangent space to $\Theta$ at the element $\theta\in \Theta$. It is easy to check that for all $\theta := (a,w,b) \in \Theta$, it holds that $T_\theta \Theta = \mathbb{R}  \times {\rm Span}\{w\}^\perp \times \mathbb{R}$, where ${\rm Span}\{w\}^\perp$ is the subspace of $\mathbb{R}^d$ containing all $d$-dimensional vectors orthogonal to $w$.

For any $\theta := (a,w,b) \in \Theta$ and any function $f:\Theta \to \mathbb{R}$ differentiable at $\theta$, we denote by $\nabla f(\theta)\in \mathbb{R}^{d+2}$ the (Euclidean) gradient of $f$. We also denote by $\nabla_\Theta f(\theta)$ the orthogonal projection of $\nabla f(\theta)$ onto $T_\theta \Theta$ and by $\nabla_{\Theta^\perp} f(\theta)$ the orthogonal projection of $\nabla f(\theta)$ onto $(T_\theta\Theta)^\perp = \{0\} \times {\rm Span}\{w\} \times \{0\} = {\rm Span}\{(0,w,0)\}$. Therefore, we have the unique orthogonal decomposition of the gradient:
\begin{equation}\label{eq:gradientDecomposition}
\nabla f(\theta) = \nabla_\Theta f(\theta) + \nabla_{\Theta^\perp} f (\theta).
\end{equation}

For all $\theta\in \Theta$ and $\mu \in {\cal P}_2(\Theta)$, define $C_{\tau, \mu}(\theta) := \frac{\langle P_\tau \mu , \Phi_\tau(\theta) \rangle_{L^2(\Omega)}}{\| P_\tau \mu \|_{L^2(\Omega)}}$, $V_{\tau, \mu}(\theta) := \langle \nabla P_\tau \mu, \nabla \Phi_\tau(\theta; \cdot) \rangle_{L^2(\Omega)} + \langle W P_\tau \mu, \Phi_\tau(\theta; \cdot) \rangle_{L^2(\Omega)}$ and 
$$
V^C_{\tau,\mu} := V_{\tau,\mu} - \sigma_\mu C_{\tau,\mu},
$$
where 
$$
\sigma_\mu := \frac{\langle \nabla_\Theta V_{\tau,\mu}, \nabla_\Theta C_{\tau,\mu} \rangle_{L^2(\Theta; \mu)} }{ \| \nabla_\Theta C_{\tau,\mu}\|^2_{L^2(\Theta; \mu)}}.
$$

Whenever $u_\mu\neq0$ and
$\|\nabla_\Theta C_{\tau,\mu}\|_{L^2(\mu)}>0$, we use the same formulas
for $\sigma_\mu$, $V^C_{\tau,\mu}$, and $F_\mu$ even when
$\mu\notin X_\tau$.  This extension will only be used on sets where
$\|u_\mu\|_{L^2(\Omega)}$ is bounded from below by a positive constant.

We identify $T_\theta\Theta$ with its canonical subspace of
$\mathbb R^{d+2}$.  Differences of tangent vector fields based at distinct
points are always understood in this ambient Euclidean space.

\begin{proposition}[Characteristics for the truncated curve]\label{prop:propositionCharacteristics}
Let $r>0$ and let \((\mu_t^r)_{t\ge0}\) be given by
Proposition~\ref{prop:existenceCurveSupportCond}. Fix \(T>0\) and assume that
there exists \(\delta>0\) such that
$\operatorname{Supp}(\mu_t^r)\subset K_{r-\delta}$ for all $t\in[0,T]$.
Define
\begin{equation*}
  b_t^r(\theta):=-\nabla_\Theta V^C_{\tau,\mu_t^r}(\theta).
\end{equation*}
Then the non-autonomous Cauchy problem
\begin{equation}\label{eq:char_truncated_corrected}
  \begin{cases}
  \partial_t\chi^r(t,s,\theta_0)=b_t^r(\chi^r(t,s,\theta_0)),\\
  \chi^r(s,s,\theta_0)=\theta_0,
  \end{cases}
\end{equation}
has a unique solution on \([0,T]\). Moreover, for every \(R>0\) there is a constant
\(L_{T,R,r}\) such that
\begin{equation*}
  d\bigl(\chi^r(t,s,\theta_1),\chi^r(t,s,\theta_2)\bigr)
  \le e^{L_{T,R,r}|t-s|}d(\theta_1,\theta_2)
\end{equation*}
for all \(\theta_1,\theta_2\in K_R\) and all \(s,t\in[0,T]\). Finally,
\begin{equation}\label{eq:pushforward_truncated_corrected}
  \mu_t^r=\chi^r(t,0,\cdot)_\#\mu_0
  \qquad\text{for all }t\in[0,T].
\end{equation}
\end{proposition}

\begin{comment}
\begin{proposition}\label{prop:propositionCharacteristics}
Let $(\mu_t^r)_{t\ge 0}$ be as in Proposition~\ref{prop:existenceCurveSupportCond}. Fix $T>0$ and assume that there exists
$\delta>0$ such that
\[
\operatorname{Supp}(\mu_t^r)\subset K_{r-\delta}
\qquad\forall t\in[0,T].
\]
Define
\[
b_t^r(\theta):=-\nabla_\Theta V_{\tau,\mu_t^r}^C(\theta),
\qquad \theta\in\Theta.
\]
Then the following hold.

\smallskip

\noindent
(i) For every $(s,\theta_0)\in[0,T]\times\Theta$, the non-autonomous Cauchy problem
\[
\left\{
\begin{aligned}
\partial_t \chi^r(t,s,\theta_0)&=b_t^r(\chi^r(t,s,\theta_0)),
\\
\chi^r(s,s,\theta_0)&=\theta_0
\end{aligned}
\right.
\]
admits a unique solution on $[0,T]$.

\smallskip

\noindent
(ii) For every $R>0$, there exists a constant $L_{T,R}>0$ such that
\[
d\!\left(\chi^r(t,s,\theta_1),\chi^r(t,s,\theta_2)\right)
\le
e^{L_{T,R}|t-s|}\,d(\theta_1,\theta_2)
\]
for all $\theta_1,\theta_2\in K_R$ and all $s,t\in[0,T]$.

\smallskip

\noindent
(iii) For every $t\in[0,T]$,
\[
\mu_t^r=\chi^r(t,0,\cdot)_\#\mu_0.
\]
\end{proposition}
\end{comment}

Finally, Theorem \ref{th:existenceNoSupport} gives the existence of a gradient curve removing the support penalization in the energy.

\begin{theorem}[Existence of the untruncated constrained flow]\label{th:existenceNoSupport}
Let \(\mu_0\in X_\tau\) be compactly supported. Then there exists a locally absolutely
continuous curve \((\mu_t)_{t\ge0}\subset X_\tau\) such that, for every \(T>0\),
\begin{equation}\label{eq:untruncated_PDE_corrected}
  \partial_t\mu_t+
  \operatorname{div}\bigl((-\nabla_\Theta V^C_{\tau,\mu_t})\mu_t\bigr)=0
  \qquad\text{in }\mathcal D'((0,T)\times\Theta),
\end{equation}
\begin{equation}\label{eq:untruncated_dissipation_bound_corrected}
  \int_0^T
  \|\nabla_\Theta V^C_{\tau,\mu_t}\|_{L^2(\Theta;\mu_t)}^2\,dt<+\infty,
\end{equation}
and there exists \(R_T>0\) such that
\begin{equation}\label{eq:support_bound_corrected}
  \operatorname{Supp}(\mu_t)\subset K_{R_T}
  \qquad\text{for every }t\in[0,T].
\end{equation}
Moreover, if \(\chi(t,s,\theta)\) denotes the flow associated with the vector field
\(-\nabla_\Theta V^C_{\tau,\mu_t}\), then
\begin{equation}\label{eq:untruncated_pushforward_corrected}
  \mu_t=\chi(t,0,\cdot)_\#\mu_0
  \qquad\text{for all }t\ge0.
\end{equation}
\end{theorem}

The proofs of these existence results are given in Section~\ref{proofTheorem3}.

\subsection{Convergence of the constrained gradient curve}\label{sec:convergence}

Our convergence result towards an eigenfunction is based on the following hypothesis on the initial measure $\mu_0$:

\begin{hypothesis}\label{hyp:geometric_support}
Assume that $\mu_0\in X_\tau$ with compact support and 
\[
\{0\}\times B\subset \operatorname{Supp}(\mu_0),
\]
where
\[
B:=S_{\mathbb R^d}(1)\times[-\sqrt d-2,\sqrt d+2].
\]
\end{hypothesis}

As we will use Hypothesis \ref{hyp:geometric_support}, we need to be sure it is not empty.

\begin{proposition}\label{prop:existence_measure_hyp1}
There exists a measure $\mu_0\in\mathcal P_2(\Theta)$ satisfying
Hypothesis~\ref{hyp:geometric_support}.
\end{proposition}

\begin{proof}
Choose $s_0\in\mathbb R$ such that
\[
\sigma_{H,\tau}(s_0)\neq 0.
\]
Let
\[
x_0:=\Bigl(\frac12,0,\dots,0\Bigr)\in \Omega=[0,1]^d,
\qquad
w_\ast:=e_1\in S_{\mathbb R^d}(1),
\qquad
b_\ast:=s_0-\frac12.
\]
Then
\[
w_\ast\cdot x_0+b_\ast=s_0,
\]
hence
\[
\sigma_{H,\tau}(w_\ast\cdot x_0+b_\ast)\neq 0.
\]
Therefore the function
\[
x\longmapsto \sigma_{H,\tau}(w_\ast\cdot x+b_\ast)
\]
is not identically zero on $\Omega$.

Let $\nu\in\mathcal P(B)$ be any probability measure with full support:
\[
\operatorname{Supp}(\nu)=B.
\]
Define
\[
\widetilde\mu_0
:=
\frac12\,\delta_{(1,w_\ast,b_\ast)}
+
\frac12\,(\delta_0\otimes \nu),
\]
where $\delta_0\otimes \nu$ denotes the measure on $\Theta=\mathbb R\times B$
supported on the slice $a=0$ and with marginal $\nu$ in the variable $\omega=(w,b)$.

Then
\[
P_\tau\widetilde\mu_0
=
\frac12\,\sigma_{H,\tau}(w_\ast\cdot +b_\ast)\neq 0
\qquad\text{in }L^2(\Omega).
\]
Let
\[
\lambda:=\|P_\tau\widetilde\mu_0\|_{L^2(\Omega)}^{-1}>0
\]
and define the scaling map
\[
S_\lambda(a,w,b):=(\lambda a,w,b),
\qquad
\mu_0:=(S_\lambda)_\#\widetilde\mu_0.
\]
Since $P_\tau$ is linear in the variable $a$, one has
\[
\|P_\tau\mu_0\|_{L^2(\Omega)}=1,
\]
hence
\[
\mu_0\in X_\tau.
\]

Moreover, since $S_\lambda$ fixes the slice $\{0\}\times B$,
\[
\operatorname{Supp}(\mu_0)
=
\{(\lambda,w_\ast,b_\ast)\}\cup(\{0\}\times B).
\]
Therefore
\[
\{0\}\times B\subset \operatorname{Supp}(\mu_0).
\]
This proves Hypothesis~\ref{hyp:geometric_support}.
\end{proof}

\begin{comment}
\begin{hypothesis}\label{hypothesisSupport}
The support of the measure $\pi_{w,b} \# \mu_0$ verifies:
$$
S_{\R^d}(1) \times [-\sqrt{d} -2,\sqrt{d} + 2] \subset {\rm Supp}(\pi_{w,b} \# \mu_0)
$$
and $\mu_0 \in X_\tau$
\end{hypothesis}
As we will use Hypothesis \ref{hypothesisSupport}, we need to be sure it is not empty.

\begin{proposition}
There exists a measure $\mu_0 \in \Prob_2(\Theta)$ satisfying Hypothesis \ref{hypothesisSupport}.
\cve{TODO: check}
\end{proposition}

\begin{proof}
%Define $\mu_{00}(a,w,b) := \delta_0(a) \times U_{S^{d-1}}(w) \times U_{[-\sqrt{d} -2,\sqrt{d} + 2]}(b) $ so that $P_\tau \mu_{00} = 0$. 

Define $\mu_{00} := \delta_1 \times \delta_{(w,b)} $ for some arbitrary $(w,b) \in S^{d-1} \times [-\sqrt{d} -2,\sqrt{d} + 2]$. Finally, set $\mu_{0} := \delta_{1/\|P_\tau \mu_{00}\|_{L^2(\Omega)}} \times \delta_{(w,b)} $ such that $\|P_\tau \mu_{0}\|_{L^2(\Omega)} = 1$ and $\mu_0$ belongs to $X_\tau$.

%The initial condition $\mu_0$ is given by :
%
%
%\begin{equation}
%\mu_0 := \frac{1}{2} (\mu_{01} + \mu_{02})
%\end{equation}
%so that $\mu_0 \in X_\tau$. 

%This $\mu_0$ is a good candidate.
\end{proof}

\end{comment}

Under such hypothesis, one gets a result of convergence in the spirit of a previous work from Bach and Chizat~\cite{BachChizat2018}:
\begin{theorem}\label{theoremConvergence}
Assume Hypothesis~\ref{hyp:geometric_support}. Let $(\mu_t)_{t\ge0}$ be the solution
given by Theorem~\ref{th:existenceNoSupport}, and assume that $\mu_t\longrightarrow \mu^\star$ in $W_2(\Theta)$as $t\to+\infty$.
Then $V_{\tau,\mu^\star}^C\equiv 0$ on $\Theta$.
\end{theorem}

%\begin{proof}
%The proof is identical to the proof of \cite[Theorem 7]{DusEhrlacher2023} and given in Appendix for completeness.
%\end{proof}

In fact the equilibrium $\mu^\star$ given by Theorem \ref{theoremConvergence} represents an eigenfunction:
\begin{theorem}
If $\mu \in X_\tau$ is such that $V^C_{\tau, \mu} = 0$ everywhere, then $P_\tau \mu$ is an eigenfunction of the Schrödinger-Neumann problem with $\sigma_\mu$ as eigenvalue.
\end{theorem}

\begin{proof}
The assertion $V^C_{\tau, \mu} = 0$ is equivalent to:

$$
\forall \theta \in \Theta, \ \langle \nabla_x P_\tau \mu, \nabla_x \Phi_\tau(\theta; \cdot) \rangle_{L^2(\Omega)} + \langle W  P_\tau \mu, \Phi_\tau(\theta; \cdot) \rangle_{L^2(\Omega)}  - \sigma_\mu \langle P_\tau \mu, \Phi_\tau(\theta; \cdot) \rangle_{L^2(\Omega)} = 0.
$$

As the space of features is dense in $H^1(\Omega)$, the conclusion is straightforward.
\end{proof}

Section~\ref{sec:proofConvergence} is devoted to the proof of Theorem~\ref{theoremConvergence}.

\section{Proof of the existence of a constrained gradient curve}
\label{proofTheorem3}

We gather in this section the proofs of Proposition~\ref{prop:existenceCurveSupportCond},  Proposition~\ref{prop:propositionCharacteristics} and Theorem~\ref{th:existenceNoSupport}. The constants involved in the results and in the proofs depend in general on the value of the dimension $d$ but, for the sake of simplicity, we do not highlight this dependence explicitly in the two next sections. 

\subsection{Notation}

We introduce here some notation which will be used throughout the proof. Let us first consider the bilinear form $\mathfrak a: H^1(\Omega) \times H^1(\Omega) \to \mathbb{R}$ as follows:
\begin{equation}\label{eq:bilinear_form}
\forall u,v\in H^1(\Omega), \quad   \mathfrak a(u,v)
  := \int_\Omega \nabla u\cdot \nabla v\,dx
     + \int_\Omega Wuv\,dx,
\end{equation}
so that \(\mathcal E(u)=\frac12\mathfrak a(u,u)\) and
\(d\mathcal E(u)[h]=\mathfrak a(u,h)\). We then define for all $\mu \in \mathcal P_2(\Theta)$, $  u_\mu:=P_\tau\mu$ so that
\begin{equation*}
  {\cal E}_\tau(\mu)=\mathcal E(u_\mu),
  \qquad
  {\cal C}_\tau(\mu)=\|u_\mu\|_{L^2(\Omega)}-1 .
\end{equation*}
For \(\mu\in P_2(\Theta)\) such that \(u_\mu\ne 0\), we then have
\begin{align}
  C_{\tau,\mu}(\theta)
  &:= \frac{\langle u_\mu,\Phi_\tau(\theta;\cdot)\rangle_{L^2(\Omega)}}
          {\|u_\mu\|_{L^2(\Omega)}},
  \label{eq:C_potential_corrected}\\
  V_{\tau,\mu}(\theta)
  &:= \mathfrak a\bigl(u_\mu,\Phi_\tau(\theta;\cdot)\bigr)
   = \langle\nabla u_\mu,\nabla\Phi_\tau(\theta;\cdot)\rangle_{L^2(\Omega)}
     +\langle W u_\mu,\Phi_\tau(\theta;\cdot)\rangle_{L^2(\Omega)}.
  \label{eq:V_potential_corrected}
\end{align}
On the constraint set \(X_\tau\), where \(\|u_\mu\|_{L^2}=1\), the Lagrange multiplier is
\begin{equation}\label{eq:sigma_corrected}
  \sigma_\mu
  :=
  \frac{\langle \nabla_\Theta V_{\tau,\mu},\nabla_\Theta C_{\tau,\mu}
          \rangle_{L^2(\Theta;\mu)}}
       {\|\nabla_\Theta C_{\tau,\mu}\|_{L^2(\Theta;\mu)}^2},
\end{equation}
and we have
\begin{equation}\label{eq:VC_corrected}
  V^C_{\tau,\mu}:=V_{\tau,\mu}-\sigma_\mu C_{\tau,\mu},
  \qquad
  F_\mu:=\nabla_\Theta V^C_{\tau,\mu}.
\end{equation}
The tangent gradient is always understood with respect to the product Riemannian
structure of \(\Theta=\mathbb R\times S^{d-1}\times\mathbb R\).

For \(r>0\), recall that $K_r:=[-r,r]\times S^{d-1}\times[-r,r]$,  $X_{\tau,r}:=\{\mu\in {\cal P}_2(\Theta): {\cal C}_\tau(\mu)=0,
  \ \operatorname{Supp}(\mu)\subset K_r\}$, and
\begin{equation*}
  {\cal E}^C_{\tau,r}(\mu):=
  \begin{cases}
   {\cal E}_\tau(\mu), & \text{if } \mu\in X_{\tau,r},\\
    +\infty, & \text{otherwise.}
  \end{cases}
\end{equation*}
We shall use the following explicit measurable choice of a minimizing
logarithm. Let $(e_1,\ldots,e_d)$ be the canonical basis of $\mathbb{R}^d$.
For $w\in\mathbb{S}^{d-1}$, set
\[
  j(w):=\min\bigl\{j\in\{1,\ldots,d\}: |w_j|<1\bigr\},
  \qquad
  e(w):=\frac{e_{j(w)}-w_{j(w)}w}{\sqrt{1-w_{j(w)}^2}}\in T_w\mathbb{S}^{d-1}.
\]
Since $d\geq 2$, the index $j(w)$ is well defined, and the map
$w\mapsto e(w)$ is Borel measurable. For $w,\widetilde w\in\mathbb{S}^{d-1}$,
let $\alpha:=\arccos(w\cdot\widetilde w)\in[0,\pi]$ and define
\[
\operatorname{Log}_{w}(\widetilde w):=
\begin{cases}
  0, & \alpha=0,\\[0.4em]
  \displaystyle
  \frac{\alpha}{\sin\alpha}
  \bigl(\widetilde w-\cos(\alpha)w\bigr),
  & 0<\alpha<\pi,\\[0.8em]
  \pi e(w), & \alpha=\pi.
\end{cases}
\]
For $\theta=(a,w,b)$ and
$\widetilde\theta=(\widetilde a,\widetilde w,\widetilde b)$, we then set
\[
  \xi_{\theta,\widetilde\theta}
  :=\bigl(\widetilde a-a,\operatorname{Log}_{w}(\widetilde w),
          \widetilde b-b\bigr)\in T_\theta\Theta.
\]
The map $(\theta,\widetilde\theta)\mapsto
\xi_{\theta,\widetilde\theta}$ is Borel measurable and satisfies
\[
  \exp_\theta\bigl(\xi_{\theta,\widetilde\theta}\bigr)
  =\widetilde\theta,
  \qquad
  \bigl|\xi_{\theta,\widetilde\theta}\bigr|
  =d(\theta,\widetilde\theta).
\]
Only measurability is needed below; continuity at the cut locus is not
required.

\subsection{First-order expansions}

We begin by establishing some auxiliary lemmas which gives first-order expansions of the quantities involved in our problem. 

\begin{lemma}[First variation of the constraint]\label{lem:C_first_variation_corrected}
Let $r>0$. The map ${\cal C}_\tau:\mathcal P_2(K_r)\longrightarrow \mathbb R$ is Lipschitz continuous with respect to $W_2$. Moreover, for every
$c>0$, there exists a constant $C=C(\tau,r,c)>0$ such that, for every $\mu,\nu\in\mathcal P_2(K_r)$ satisfying
$\|u_\mu\|_{L^2(\Omega)}\geq c$, and every $\gamma\in\Gamma(\mu,\nu)$ with
$\int_{\Theta^2} d(\theta, \widetilde\theta)^2 \,d\gamma(\theta, \widetilde{\theta})\leq 1$, one has
\begin{equation}
\left|
{\cal C}_\tau(\nu)-{\cal C}_\tau(\mu)
-
\int_{\Theta^2}
\nabla_\Theta C_{\tau,\mu}(\theta)
\cdot\xi_{\theta,\widetilde\theta}
\,\mathrm d\gamma(\theta,\widetilde\theta)
\right|
\leq
C\int_{\Theta^2} d(\theta, \widetilde\theta)^2 \,d\gamma(\theta, \widetilde{\theta}).
\label{eq:constraint-second-order}
\end{equation}
In particular, if $\gamma\in\Gamma_o(\mu,\nu)$, then
\begin{equation}\label{eq:constraint-second-order2}
{\cal C}_\tau(\nu)
=
{\cal C}_\tau(\mu)
+
\int_{\Theta^2}
\nabla_\Theta C_{\tau,\mu}(\theta)
\cdot\xi_{\theta,\widetilde\theta}
\,\mathrm d\gamma(\theta,\widetilde\theta)
+
O\!\left(W_2(\mu,\nu)^2\right).
\end{equation}
The constant in the remainder is uniform when $\mu$ ranges in a
subset of $\mathcal P_2(K_r)$ on which
$\|u_\mu\|_{L^2(\Omega)}$ is bounded from below by a positive constant.
\end{lemma}

\begin{proof}
Since $\sigma_{H,\tau}$ is smooth with bounded derivatives, the map $K_r\ni\theta\longmapsto \Phi_\tau(\theta;\cdot)\in H^1(\Omega)$
is of class $C^2$, with uniformly bounded first and second order
derivatives. For $\theta,\widetilde\theta\in K_r$, let
\[
\kappa_{\theta,\widetilde\theta}(t)
:=\exp_\theta\bigl(t\xi_{\theta,\widetilde\theta}\bigr),
 \qquad t\in[0,1].
\]
This is a minimizing geodesic from $\theta$ to $\widetilde\theta$ and is
entirely contained in $K_r$. Taylor's formula along this geodesic gives
\[
\Phi_\tau(\widetilde\theta;\cdot)
 =\Phi_\tau(\theta;\cdot)
  +D_\Theta\Phi_\tau(\theta;\cdot)
      [\xi_{\theta,\widetilde\theta}]
  +R_\tau(\theta,\widetilde\theta;\cdot),
\]
where
\[
R_\tau(\theta,\widetilde\theta;\cdot)
 =\int_0^1(1-t)\,
   \operatorname{Hess}_\Theta\Phi_\tau
   \bigl(\kappa_{\theta,\widetilde\theta}(t);\cdot\bigr)
   [\dot\kappa_{\theta,\widetilde\theta}(t),
    \dot\kappa_{\theta,\widetilde\theta}(t)]\,dt.
\]
Consequently,
\[
\|R_\tau(\theta,\widetilde\theta;\cdot)\|_{H^1(\Omega)}
 \leq C(\tau,r)\,d(\theta,\widetilde\theta)^2.
\]

%Since $\sigma_{H,\tau}$ is smooth and has bounded derivatives,
%the map
%$K_r\ni\theta\longmapsto\Phi_\tau(\theta;\cdot)\in H^1(\Omega)$ is of class $C^2$, with uniformly bounded first and second order derivatives. Hence, for every $\theta,\widetilde\theta\in K_r$,
%\[
%\Phi_\tau(\widetilde\theta;\cdot)
%=
%\Phi_\tau(\theta;\cdot)
%+
%D_\Theta\Phi_\tau(\theta;\cdot)
% [\xi_{\theta,\widetilde\theta}]
%+
%R_\tau(\theta,\widetilde\theta;\cdot),
%\]
%where
%\[
%\|R_\tau(\theta,\widetilde\theta;\cdot)\|_{H^1(\Omega)}
%\leq
%C(\tau,r)d(\theta,\widetilde\theta)^2.
%\]
For $\gamma\in\Gamma(\mu,\nu)$, set
\[
h_\gamma
:=
\int_{\Theta^2}
D_\Theta\Phi_\tau(\theta;\cdot)
 [\xi_{\theta,\widetilde\theta}]
\,\mathrm d\gamma(\theta,\widetilde\theta)
\quad \mbox{ and }\quad 
\rho_\gamma
:=
\int_{\Theta^2}
R_\tau(\theta,\widetilde\theta;\cdot)
\,\mathrm d\gamma(\theta,\widetilde\theta).
\]
Then, $u_\nu-u_\mu=h_\gamma+\rho_\gamma$,
with
\[
\|h_\gamma\|_{H^1(\Omega)}
\leq C(\tau,r)\left(\int_{\Theta^2} d(\theta, \widetilde\theta)^2 \,d\gamma(\theta, \widetilde{\theta})\right)^{1/2}
\quad \mbox{ and } \quad
\|\rho_\gamma\|_{H^1(\Omega)}
\leq C(\tau,r)\int_{\Theta^2} d(\theta, \widetilde\theta)^2 \,d\gamma(\theta, \widetilde{\theta}).
\]

We shall use the following elementary estimate. For every $c>0$, there
exists a constant $C_c>0$ such that, for all $u,v\in L^2(\Omega)$ with
$\|u\|_{L^2(\Omega)}\geq c$,
\[
\left|
 \|u+v\|_{L^2(\Omega)}
 -\|u\|_{L^2(\Omega)}
 -\frac{\langle u,v\rangle_{L^2(\Omega)}}{\|u\|_{L^2(\Omega)}}
\right|
\leq C_c\|v\|_{L^2(\Omega)}^2.
\]
Indeed, if $\|v\|_{L^2(\Omega)}\leq c/2$, this follows from Taylor's
formula applied to the norm on the segment $[u,u+v]$, which remains at
distance at least $c/2$ from the origin. If
$\|v\|_{L^2(\Omega)}>c/2$, the expression in the absolute value is
bounded by $2\|v\|_{L^2(\Omega)}$, and hence by
$4c^{-1}\|v\|_{L^2(\Omega)}^2$. 
Applying this estimate with $u=u_\mu$ and
$v=h_\gamma+\rho_\gamma$, and absorbing the contribution of
$\rho_\gamma$ in the remainder, we obtain
\begin{align*}
{\cal C}_\tau(\nu)-{\cal C}_\tau(\mu)
&=\frac{\langle u_\mu,h_\gamma+\rho_\gamma\rangle_{L^2(\Omega)}}
        {\|u_\mu\|_{L^2(\Omega)}}
  +O_c\!\left(\|h_\gamma+\rho_\gamma\|_{L^2(\Omega)}^2\right)\\
&=\frac{\langle u_\mu,h_\gamma\rangle_{L^2(\Omega)}}
        {\|u_\mu\|_{L^2(\Omega)}}
  +O_c\!\left(
       \|\rho_\gamma\|_{L^2(\Omega)}
       +\|h_\gamma+\rho_\gamma\|_{L^2(\Omega)}^2
     \right).
\end{align*}

%Applying this estimate with $u=u_\mu$ and
%$v=h_\gamma+\rho_\gamma$, we obtain
%\[
%\begin{aligned}
%N(u_\mu+h_\gamma+\rho_\gamma)-N(u_\mu)
%&=
%\frac{\langle u_\mu,h_\gamma\rangle_{L^2(\Omega)}}
%     {\|u_\mu\|_{L^2(\Omega)}}  \\
%&\quad
%+O_c\left(
% \|\rho_\gamma\|_{L^2(\Omega)}
% +\|h_\gamma+\rho_\gamma\|_{L^2(\Omega)}^2
%\right).
%\end{aligned}
%\]
Since
\[
\|h_\gamma\|_{L^2(\Omega)}^2
 \leq C(\tau,r)\int_{\Theta^2}
 d(\theta,\widetilde\theta)^2\,d\gamma(\theta,\widetilde\theta)
\]
and
\[
\|\rho_\gamma\|_{L^2(\Omega)}
 \leq C(\tau,r)\int_{\Theta^2}
 d(\theta,\widetilde\theta)^2\,d\gamma(\theta,\widetilde\theta),
\]
the assumption
$\int_{\Theta^2}d(\theta,\widetilde\theta)^2\,
d\gamma(\theta,\widetilde\theta)\leq 1$ 
shows that the remainder is bounded by
$
C(\tau,r,c)
\int_{\Theta^2}d(\theta,\widetilde\theta)^2\,
d\gamma(\theta,\widetilde\theta)$.
Moreover,
\[
\frac{\langle u_\mu,h_\gamma\rangle_{L^2(\Omega)}}
     {\|u_\mu\|_{L^2(\Omega)}}
=
\int_{\Theta^2}
\nabla_\Theta C_{\tau,\mu}(\theta)
 \cdot\xi_{\theta,\widetilde\theta}\,
d\gamma(\theta,\widetilde\theta).
\]
This proves \eqref{eq:constraint-second-order} and its optimal-plan
consequence \eqref{eq:constraint-second-order2}.

Finally, for an optimal plan $\gamma\in\Gamma_o(\mu,\nu)$,
\[
\begin{aligned}
|{\cal C}_\tau(\nu)-{\cal C}_\tau(\mu)|
&\leq \|u_\nu-u_\mu\|_{L^2(\Omega)} \\
&\leq
\int_{\Theta^2}
\|\Phi_\tau(\widetilde\theta;\cdot)
      -\Phi_\tau(\theta;\cdot)\|_{L^2(\Omega)}
\,\mathrm d\gamma(\theta,\widetilde\theta) \\
&\leq C(\tau,r)W_2(\mu,\nu),
\end{aligned}
\]
which proves the claimed Lipschitz continuity, including at measures for which
$u_\mu=0$.
\end{proof}

\begin{lemma}[First variation of the energy]\label{lem:E_first_variation_corrected}
Let $r>0$. The map
$
{\cal E}_\tau:\mathcal P_2(K_r)\longrightarrow\mathbb R
$ is Lipschitz continuous with respect to $W_2$. Moreover, there exists
$C=C(\tau,r,\|W\|_{L^\infty(\Omega)})>0$ such that, for every
$\mu,\nu\in\mathcal P_2(K_r)$ and every
$\gamma\in\Gamma(\mu,\nu)$ with $\int_{\Theta^2} d(\theta, \widetilde\theta)^2 \,d\gamma(\theta, \widetilde{\theta})\leq 1$,
\begin{equation}
\left|
{\cal E}_\tau(\nu)-{\cal E}_\tau(\mu)
-
\int_{\Theta^2}
\nabla_\Theta V_{\tau,\mu}(\theta)
\cdot\xi_{\theta,\widetilde\theta}
\,\mathrm d\gamma(\theta,\widetilde\theta)
\right|
\leq
C\,\int_{\Theta^2} d(\theta, \widetilde\theta)^2 \,d\gamma(\theta, \widetilde{\theta}).
\label{eq:energy-second-order}
\end{equation}
In particular, for $\gamma\in\Gamma_o(\mu,\nu)$,
\begin{equation}
{\cal E}_\tau(\nu)
=
{\cal E}_\tau(\mu)
+
\int_{\Theta^2}
\nabla_\Theta V_{\tau,\mu}(\theta)
\cdot\xi_{\theta,\widetilde\theta}
\,\mathrm d\gamma(\theta,\widetilde\theta)
+
O\!\left(W_2(\mu,\nu)^2\right).
\label{eq:energy-second-order2}
\end{equation}
\end{lemma}

\begin{proof}
Using the notation introduced in the proof of the previous lemma, $u_\nu-u_\mu=h_\gamma+\rho_\gamma$, where
\[
\|h_\gamma\|_{H^1(\Omega)}
\leq C(\tau,r)\left(\int_{\Theta^2} d(\theta, \widetilde\theta)^2 \,d\gamma(\theta, \widetilde{\theta})\right)^{1/2}
\quad \mbox{ and } \quad
\|\rho_\gamma\|_{H^1(\Omega)}
\leq C(\tau,r)\int_{\Theta^2} d(\theta, \widetilde\theta)^2 \,d\gamma(\theta, \widetilde{\theta}).
\]
Since ${\cal E}_\tau(\mu)=\frac12 \mathfrak a(u_\mu,u_\mu)$,
\[
{\cal E}_\tau(\nu)-{\cal E}_\tau(\mu)
=
\mathfrak a(u_\mu,h_\gamma)
+
 \mathfrak a(u_\mu,\rho_\gamma)
+
\frac12
\mathfrak a(u_\nu-u_\mu,u_\nu-u_\mu).
\]
The last two terms satisfy
\[
|\mathfrak a(u_\mu,\rho_\gamma)|
+
\left|
\mathfrak a(u_\nu-u_\mu,u_\nu-u_\mu)
\right|
\leq
C\int_{\Theta^2} d(\theta, \widetilde\theta)^2 \,d\gamma(\theta, \widetilde{\theta}),
\]
because $u_\mu$ is uniformly bounded in $H^1(\Omega)$ for
$\mu\in\mathcal P_2(K_r)$. Finally,
\[
\begin{aligned}
\mathfrak a(u_\mu,h_\gamma)
&=
\int_{\Theta^2}
\mathfrak a\!\left(
u_\mu,
D_\Theta\Phi_\tau(\theta;\cdot)
 [\xi_{\theta,\widetilde\theta}]
\right)
\,\mathrm d\gamma(\theta,\widetilde\theta) \\
&=
\int_{\Theta^2}
\nabla_\Theta V_{\tau,\mu}(\theta)
\cdot\xi_{\theta,\widetilde\theta}
\,\mathrm d\gamma(\theta,\widetilde\theta).
\end{aligned}
\]
This proves \eqref{eq:energy-second-order} and \eqref{eq:energy-second-order2}.

It remains to prove the Lipschitz continuity. Let
$\gamma\in\Gamma_o(\mu,\nu)$. The first-order estimates for the feature
map give
\[
\|u_\nu-u_\mu\|_{H^1(\Omega)}
\leq
\int_{\Theta^2}
\|\Phi_\tau(\widetilde\theta;\cdot)
       -\Phi_\tau(\theta;\cdot)\|_{H^1(\Omega)}
\,d\gamma(\theta,\widetilde\theta)
\leq C(\tau,r)W_2(\mu,\nu).
\]
Moreover, $u_\mu$ and $u_\nu$ are uniformly bounded in $H^1(\Omega)$
when $\mu,\nu\in\mathcal P_2(K_r)$. Therefore,
\[
\begin{aligned}
|{\cal E}_\tau(\nu)-{\cal E}_\tau(\mu)|
&=
\frac12
\left|
\mathfrak a(u_\nu+u_\mu,u_\nu-u_\mu)
\right| \\
&\leq
C(\|W\|_{L^\infty(\Omega)})
\bigl(
 \|u_\nu\|_{H^1(\Omega)}
 +\|u_\mu\|_{H^1(\Omega)}
\bigr)
\|u_\nu-u_\mu\|_{H^1(\Omega)} \\
&\leq
C(\tau,r,\|W\|_{L^\infty(\Omega)})
W_2(\mu,\nu),
\end{aligned}
\]
which proves the claimed Lipschitz continuity.
\end{proof}

\begin{proposition}[Lower semicontinuity and compactness]\label{prop:lsc_compact_corrected}
The functional \({\cal E}^C_{\tau,r}:{\cal P}_2(\Theta)\to\mathbb R\cup\{+\infty\}\) is lower
semicontinuous. Moreover, its sublevels are compact in \({\cal P}_2(\Theta)\).
\end{proposition}

\begin{proof}
Let $(\mu_n)_{n\in\mathbb N}\subset\mathcal P_2(\Theta)$ be such that $\mu_n\to\mu$ in $W_2$. If
$\liminf_{n\to\infty}{\cal E}^C_{\tau,r}(\mu_n)=+\infty$, there is nothing to prove. Otherwise, set
$
\ell:=\liminf_{n\to\infty}{\cal E}^C_{\tau,r}(\mu_n)<+\infty$ 
and extract a subsequence, still denoted by $(\mu_n)_{n\in\mathbb N}$ for the sake of simplicity, such that ${\cal E}^C_{\tau,r}(\mu_n)\longrightarrow\ell$ as $n$ goes to infinity. After discarding finitely many terms, one has
${\cal E}^C_{\tau,r}(\mu_n)<+\infty$, and therefore
$\mu_n\in X_{\tau,r}$ for every $n \in \mathbb{N}$.
Since $K_r$ is closed, $\mu_n(K_r)=1$ for every $n$, and
$\mu_n\to\mu$ narrowly, the Portmanteau theorem gives
$
1=\limsup_{n\to\infty}\mu_n(K_r)\leq \mu(K_r)$.
Thus $\mu(K_r)=1$, and consequently
$\operatorname{Supp}(\mu)\subset K_r$.

By the continuity of ${\cal C}_\tau$ on $\mathcal P_2(K_r)$,
${\cal C}_\tau(\mu)=0$, hence $\mu\in X_{\tau,r}$. Moreover, the continuity
of ${\cal E}_\tau$ on $\mathcal P_2(K_r)$ gives
\[
{\cal E}^C_{\tau,r}(\mu)
=
{\cal E}_\tau(\mu)
=
\lim_{n\to\infty}{\cal E}_\tau(\mu_n)
=
\ell.
\]
This proves the lower semicontinuity.

Finally, every sublevel of \({\cal E}^C_{\tau,r}\) is contained in \({\cal P}_2(K_r)\), which is compact;
lower semicontinuity makes the sublevel closed, hence compact.
\end{proof}

\subsection{Local slope and Lagrange multiplier}\label{sec:localSlope}

The aim of this section is to prove Proposition~\ref{prop:existenceCurveSupportCond}, which states the existence of a constrained curve on a truncated domain. 

For the particular functional $\mathcal E_{\tau,r}^C$, the
definition~\eqref{eq:relaxedSlope} agrees on its effective domain with the usual relaxed slope
defined by means of sequences with uniformly bounded energy. Indeed,
if
$
\mu_n\to\mu$ in $W_2$ and
$\sup_n\mathcal E_{\tau,r}^C(\mu_n)<+\infty$,
then $\mu_n\in X_{\tau,r}$ for every $n$. Since $X_{\tau,r}$ is closed
and ${\cal E}_\tau$ is continuous on $\mathcal P_2(K_r)$, one has
$\mu\in X_{\tau,r}$ and
$
\mathcal E_{\tau,r}^C(\mu_n)
={\cal E}_\tau(\mu_n)
\longrightarrow
{\cal E}_\tau(\mu)
=\mathcal E_{\tau,r}^C(\mu)$.

\begin{lemma}[Nondegeneracy of the constraint]
\label{lem:constraint_nondegenerate_corrected}
Let $\mu\in\mathcal P_2(\Theta)$ be such that $u_\mu\neq0$, and set
$
 A(\mu):=\int_\Theta a^2\,d\mu(a,w,b)$.
Then $A(\mu)>0$ and
\begin{equation}
\label{eq:constraint-nondegeneracy}
 \|\nabla_\Theta C_{\tau,\mu}\|_{L^2(\Theta;\mu)}^2
 \geq \frac{\|u_\mu\|_{L^2(\Omega)}^2}{A(\mu)}.
\end{equation}
In particular, if $\mu\in X_\tau$, then
\[
 \|\nabla_\Theta C_{\tau,\mu}\|_{L^2(\Theta;\mu)}^2
 \geq \frac1{A(\mu)},
\]
and the denominator in the definition of $\sigma_\mu$ is strictly
positive.
\end{lemma}

\begin{proof}
Since $\Phi_\tau(a,w,b;x)=a\sigma_{H,\tau}(w\cdot x+b)$, there exists a
function $c_\mu$ such that
\[
 C_{\tau,\mu}(a,w,b)=a c_\mu(w,b),
 \qquad
 \partial_a C_{\tau,\mu}(a,w,b)=c_\mu(w,b).
\]
Moreover,
\[
 \int_\Theta a\,\partial_a C_{\tau,\mu}\,d\mu
 =\int_\Theta C_{\tau,\mu}\,d\mu
 =\frac{\langle u_\mu,u_\mu\rangle_{L^2(\Omega)}}
        {\|u_\mu\|_{L^2(\Omega)}}
 =\|u_\mu\|_{L^2(\Omega)}.
\]
If $A(\mu)=0$, then $a=0$ $\mu$-a.e., hence $u_\mu=0$, a contradiction.
The Cauchy--Schwarz inequality therefore gives
\[
 \|u_\mu\|_{L^2(\Omega)}
 \leq A(\mu)^{1/2}
       \|\partial_a C_{\tau,\mu}\|_{L^2(\mu)}
 \leq A(\mu)^{1/2}
       \|\nabla_\Theta C_{\tau,\mu}\|_{L^2(\mu)},
\]
which proves \eqref{eq:constraint-nondegeneracy}.
\end{proof}

\begin{proposition}[Interior constrained local slope]\label{prop:interior_slope_corrected}
Let \(r>0\), \(0< \delta < r\), and \(\mu\in X_{\tau,r}\) be such that
$\operatorname{Supp}(\mu)\subset K_{r-\delta}$.
Then
\begin{equation}\label{eq:interior_slope_corrected}
  |\partial {\cal E}^C_{\tau,r}|(\mu)
  =\|\nabla_\Theta V_{\tau,\mu}-\sigma_\mu\nabla_\Theta C_{\tau,\mu}\|_{L^2(\Theta;\mu)}
  =\|F_\mu\|_{L^2(\Theta;\mu)}.
\end{equation}
\end{proposition}

\begin{proof}
Let \(\nu\in X_{\tau,r}\), let \(\gamma\in\Gamma_o(\mu,\nu)\), and assume that
\(W_2(\mu,\nu)\) is small. By Lemmas~\ref{lem:C_first_variation_corrected} and
\ref{lem:E_first_variation_corrected},
\begin{align*}
  0={\cal C}_\tau(\nu)-{\cal C}_\tau(\mu)
  &=\int_{\Theta^2}\nabla_\Theta C_{\tau,\mu}(\theta)\cdot
    \xi_{\theta,\widetilde\theta}\,d\gamma+o(W_2(\mu,\nu)),\\
  {\cal E}_\tau(\nu)-{\cal E}_\tau(\mu)
  &=\int_{\Theta^2}\nabla_\Theta V_{\tau,\mu}(\theta)\cdot
    \xi_{\theta,\widetilde\theta}\,d\gamma+o(W_2(\mu,\nu)).
\end{align*}
Therefore, for any \(\lambda\in\mathbb R\),
\begin{equation*}
  {\cal E}_\tau(\nu)-{\cal E}_\tau(\mu)
  =\int_{\Theta^2}
     \bigl(\nabla_\Theta V_{\tau,\mu}-\lambda\nabla_\Theta C_{\tau,\mu}\bigr)(\theta)
     \cdot\xi_{\theta,\widetilde\theta}\,d\gamma
   +o(W_2(\mu,\nu)).
\end{equation*}
Choosing \(\lambda=\sigma_\mu\) and applying Cauchy--Schwarz gives
\begin{equation*}
  \bigl({\cal E}^C_{\tau,r}(\mu)-{\cal E}^C_{\tau,r}(\nu)\bigr)^+
  \le \|F_\mu\|_{L^2(\Theta; \mu)} W_2(\mu,\nu)+o(W_2(\mu,\nu)),
\end{equation*}
and hence
\begin{equation*}
  |\partial {\cal E}^C_{\tau,r}|(\mu)\le \|F_\mu\|_{L^2(\Theta; \mu)}.
\end{equation*}

We now prove the reverse inequality. If $F_\mu=0$, there is nothing
to prove. Assume therefore that $F_\mu\neq0$. Since $F_\mu$ is continuous, it is bounded on the compact set
$\operatorname{Supp}(\mu)$. Choose $s_0>0$ such that
\[
s_0\|F_\mu\|_{L^\infty(\operatorname{Supp}(\mu))}
 <\min\{\delta,\pi\}.
\]
Then, for every $s\in(0,s_0)$ and
$\theta\in\operatorname{Supp}(\mu)$, the displacement
$-sF_\mu(\theta)$ lies within the injectivity radius of $\Theta$ and
cannot reach the boundary of $K_r$. Hence the map
\[
T_s(\theta):=\exp_\theta\bigl(-sF_\mu(\theta)\bigr)
\]
maps $\operatorname{Supp}(\mu)$ into $K_r$, and the minimizing
logarithm appearing in Lemmas~\ref{lem:C_first_variation_corrected} and~\ref{lem:E_first_variation_corrected} is
\[
\xi_{\theta,T_s(\theta)}=-sF_\mu(\theta).
\]

Set $\widetilde\mu_s:=(T_s)_\#\mu$.Using Lemma~\ref{lem:C_first_variation_corrected} with the coupling $\pi_s:=(\operatorname{id},T_s)_\#\mu$,
we obtain
\[
\begin{aligned}
{\cal C}_\tau(\widetilde\mu_s)-{\cal C}_\tau(\mu)
&=
-s
\langle
\nabla_\Theta C_{\tau,\mu},F_\mu
\rangle_{L^2(\Theta;\mu)}
+
O(s^2).
\end{aligned}
\]
By the definition of $\sigma_\mu$,
\[
\langle
F_\mu,\nabla_\Theta C_{\tau,\mu}
\rangle_{L^2(\Theta;\mu)}=
\langle
\nabla_\Theta V_{\tau,\mu}
-\sigma_\mu\nabla_\Theta C_{\tau,\mu},
\nabla_\Theta C_{\tau,\mu}
\rangle_{L^2(\Theta;\mu)}
=0.
\]
It follows that ${\cal C}_\tau(\widetilde\mu_s)=O(s^2)$.

Since $\operatorname{Supp}(\mu)\subset K_{r-\delta}$ and
$\displaystyle  \sup_{\theta\in\operatorname{Supp}(\mu)}
  d\bigl(T_s(\theta),\theta\bigr)=O(|s|)$,
there exist $s_1>0$ and $\delta_1>0$ such that
$
  \operatorname{Supp}(\widetilde\mu_s)\subset K_{r-\delta_1}$
  for every $|s|<s_1$. Moreover, $\|u_{\widetilde\mu_s}\|_{L^2(\Omega)}$ remains bounded away
from zero for $|s|<s_1$. By continuity, after decreasing $s_1$ and the
range of $t$ if necessary, the same is true for
$\|u_{\mu_{s,t}}\|_{L^2(\Omega)}$. The maps $s\longmapsto C_{\tau,\widetilde\mu_s}$ and $s\longmapsto \nabla_\Theta C_{\tau,\widetilde\mu_s}$ are therefore of class $C^1$, with values respectively in
$C^2(K_r)$ and $C^1(K_r;T\Theta)$. After reducing $s_1$ if necessary,
there exists $t_1>0$ such that, for $|s|<s_1$ and $|t|<t_1$, the map
\[
  S_{s,t}(\vartheta)
  :=\exp_\vartheta\!\left(
       t\nabla_\Theta C_{\tau,\widetilde\mu_s}(\vartheta)
     \right)
\]
maps $\operatorname{Supp}(\widetilde\mu_s)$ into $K_r$. Set
\[
  \mu_{s,t}:=(S_{s,t})_\#\widetilde\mu_s,
  \qquad
  G(s,t):={\cal C}_\tau(\mu_{s,t}).
\]
Writing both push-forwards with respect to the fixed measure $\mu$, namely
\[
  G(s,t)
  =\left\|
      \int_\Theta
      \Phi_\tau\!\left(
        S_{s,t}\bigl(T_s(\theta)\bigr);\,\cdot
      \right)d\mu(\theta)
    \right\|_{L^2(\Omega)}-1,
\]
shows, by differentiation under the integral sign on the compact support
of $\mu$, that $G$ is of class $C^1$ near $(0,0)$. Furthermore,
\[
  G(0,0)=0,
  \qquad
  \partial_tG(0,0)
  =\int_\Theta
     |\nabla_\Theta C_{\tau,\mu}|^2\,d\mu
  >0
\]
by Lemma~\ref{lem:constraint_nondegenerate_corrected}. The implicit function theorem thus provides a $C^1$ function
$t=t(s)$, defined near $s=0$, such that $t(0)=0$ and
$G(s,t(s))=0$. Since $G(s,0)=O(s^2)$ and $\partial_tG$ is bounded away
from zero near $(0,0)$, the identity
\[
  0=G(s,0)+t(s)\int_0^1
       \partial_tG\bigl(s,\lambda t(s)\bigr)\,d\lambda
\]
implies that $t(s)=O(s^2)$.

Set $\mu_s:=\mu_{s,t(s)}$. For $s>0$ sufficiently small, $\mu_s\in X_{\tau,r}$. The explicit transport plans associated with $T_s$ and
$S_{s,t(s)}$ then give
\[
\begin{aligned}
W_2(\mu_s,\mu)
&\leq
W_2(\mu_s,\widetilde\mu_s)
+
W_2(\widetilde\mu_s,\mu) \\
&\leq
|t(s)|
\|\nabla_\Theta C_{\tau,\widetilde\mu_s}\|_
 {L^2(\Theta;\widetilde\mu_s)}
+
s\|F_\mu\|_{L^2(\Theta;\mu)} \\
&\leq
s\|F_\mu\|_{L^2(\Theta;\mu)}
+
O(s^2).
\end{aligned}
\]
%Notice that only this upper bound is required.

Using Lemma~\ref{lem:E_first_variation_corrected} with the transport plan induced by $T_s$, and using
the Lipschitz continuity of ${\cal E}_\tau$,
we obtain
\[
\begin{aligned}
{\cal E}_\tau(\mu_s)-{\cal E}_\tau(\mu)
&=
-s
\langle
\nabla_\Theta V_{\tau,\mu},F_\mu
\rangle_{L^2(\Theta;\mu)}
+
O(s^2) \\
&=
-s\|F_\mu\|_{L^2(\Theta;\mu)}^2
+
O(s^2).
\end{aligned}
\]
Indeed,
\[
\begin{aligned}
\langle\nabla_\Theta V_{\tau,\mu},F_\mu\rangle
&=
\langle
F_\mu+\sigma_\mu\nabla_\Theta C_{\tau,\mu},
F_\mu
\rangle \\
&=
\|F_\mu\|_{L^2(\Theta;\mu)}^2.
\end{aligned}
\]
%Consequently,
%\[
%\begin{aligned}
%|\partial {\cal E}^C_{\tau,r}|(\mu)
%&\geq
%\limsup_{s\downarrow0}
%\frac{{\cal E}_\tau(\mu)-{\cal E}_\tau(\mu_s)}
%     {W_2(\mu_s,\mu)} \\
%&\geq
%\lim_{s\downarrow0}
%frac{
%s\|F_\mu\|_{L^2(\Theta;\mu)}^2+O(s^2)
%}{
%s\|F_\mu\|_{L^2(\Theta;\mu)}+O(s^2)
%} \\
%&=
%\|F_\mu\|_{L^2(\Theta;\mu)}.
%\end{aligned}
%\]
%Together with the previously established upper bound, this proves
%the result.

Set $A:=\|F_\mu\|_{L^2(\Theta;\mu)}>0$. The preceding estimates imply
that, for some constant $C>0$ and all sufficiently small $s>0$,
\[
{\cal E}_\tau(\mu)-{\cal E}_\tau(\mu_s)
 \geq sA^2-Cs^2>0,
\qquad
W_2(\mu_s,\mu)\leq sA+Cs^2.
\]
Therefore the positive part in the definition of the local slope can be
removed along this perturbation, and
\begin{align*}
|\partial {\cal E}^C_{\tau,r}|(\mu)
&\geq \limsup_{s\downarrow0}
 \frac{\bigl({\cal E}_\tau(\mu)-{\cal E}_\tau(\mu_s)\bigr)^+}
      {W_2(\mu_s,\mu)}\\
&\geq \liminf_{s\downarrow0}
 \frac{sA^2-Cs^2}{sA+Cs^2}
 =A.
\end{align*}
Together with the previously established upper bound, this proves the
claimed identity.
\end{proof}

\begin{proposition}[No relaxation in the interior]\label{prop:no_relaxation_corrected}
Let \(r>0\), \(0< \delta< r\), and \(\mu\in X_{\tau,r}\) be such that
\(\operatorname{Supp}(\mu)\subset K_{r-\delta}\). Then
\begin{equation}\label{eq:no_relaxation_corrected}
  |\partial {\cal E}^C_{\tau,r}|_{\rm rel}(\mu)
  =|\partial {\cal E}^C_{\tau,r}|(\mu)
  =\|F_\mu\|_{L^2(\Theta;\mu)}.
\end{equation}
\end{proposition}

\begin{proof}
The inequality \(|\partial {\cal E}^C_{\tau,r}|_{\rm rel}(\mu)\le |\partial {\cal E}^C_{\tau,r}|(\mu)\)
follows from the constant recovery sequence. We prove the converse inequality.

Let \((\mu_n)_n\subset X_{\tau,r}\) be such that \(\mu_n\to\mu\) in \(W_2\) and
\({\cal E}_\tau(\mu_n)\to {\cal E}_\tau(\mu)\). It is enough to show that
\begin{equation*}
  \liminf_{n\to\infty}|\partial {\cal E}^C_{\tau,r}|(\mu_n)
  \ge \|F_\mu\|_{L^2(\Theta;\mu)}.
\end{equation*}
Choose \(\zeta\in C^\infty(\Theta;[0,1])\) such that \(\zeta=1\) on a neighborhood of
\(\operatorname{Supp}(\mu)\) and \(\operatorname{Supp}(\zeta)\subset K_{r-\delta/2}\).

The first-order estimates for the feature map yield
\[
\|u_{\mu_n}-u_\mu\|_{H^1(\Omega)}
 \leq C(\tau,r)W_2(\mu_n,\mu)\longrightarrow0.
\]
Since $\mu_n,\mu\in X_{\tau,r}$, their represented functions have unit
$L^2$ norm. The uniform bounds on the differentiated features therefore
give
\[
\|\nabla_\Theta V_{\tau,\mu_n}
     -\nabla_\Theta V_{\tau,\mu}\|_{C^0(K_r)}
+
\|\nabla_\Theta C_{\tau,\mu_n}
     -\nabla_\Theta C_{\tau,\mu}\|_{C^0(K_r)}
\longrightarrow0.
\]
Introduce
\[
\mathcal N_n
 :=\int_\Theta
   \nabla_\Theta V_{\tau,\mu_n}\cdot
   \nabla_\Theta C_{\tau,\mu_n}\,d\mu_n,
\qquad
\mathcal D_n
 :=\int_\Theta
   |\nabla_\Theta C_{\tau,\mu_n}|^2\,d\mu_n.
\]
Uniform convergence of the integrands and narrow convergence of
$\mu_n$ imply
\[
\mathcal N_n\longrightarrow
 \int_\Theta
 \nabla_\Theta V_{\tau,\mu}\cdot
 \nabla_\Theta C_{\tau,\mu}\,d\mu,
\qquad
\mathcal D_n\longrightarrow
 \|\nabla_\Theta C_{\tau,\mu}\|_{L^2(\Theta;\mu)}^2.
\]
Furthermore, Lemma~\ref{lem:constraint_nondegenerate_corrected} and
$\int_\Theta a^2\,d\mu_n\leq r^2$ give
$\mathcal D_n\geq r^{-2}$. Consequently,
\[
\sigma_{\mu_n}\longrightarrow\sigma_\mu,
\qquad
\|F_{\mu_n}-F_\mu\|_{C^0(K_r)}\longrightarrow0.
\]
In particular, all these quantities are uniformly bounded on $K_r$.

%On \(K_r\), the quantities
%\(\nabla_\Theta V_{\tau,\mu_n}\), \(\nabla_\Theta C_{\tau,\mu_n}\), \(\sigma_{\mu_n}\), and
%\(F_n\) are uniformly bounded; moreover they converge uniformly on compact subsets of
%\(K_r\) to the corresponding quantities associated with \(\mu\). This follows from the
%first-order estimates above and from Lemma~\ref{lem:constraint_nondegenerate_corrected},
%since \(\int_\Theta a^2\,d\mu_n\le r^2\).

Set
\(F_n:=F_{\mu_n}\) and \(C_n:=C_{\tau,\mu_n}\). Define
\[
D_n
:=
\int_\Theta
\zeta|\nabla_\Theta C_n|^2\,\mathrm d\mu_n,
\qquad
B_n
:=
\int_\Theta
\zeta F_n\cdot\nabla_\Theta C_n\,\mathrm d\mu_n,
\]
and
\[
H_n:=\zeta\nabla_\Theta C_n,
\qquad
\alpha_n:=\frac{B_n}{D_n},
\qquad
G_n:=\zeta F_n-\alpha_nH_n.
\]
For all sufficiently large $n$, the quantity $D_n$ is bounded from
below by a positive constant. Indeed,
$
D_n
\longrightarrow
\|\nabla_\Theta C_{\tau,\mu}\|_{L^2(\Theta;\mu)}^2>0$,
because $\zeta=1$ on a neighborhood of
$\operatorname{Supp}(\mu)$. Moreover,
$
B_n
\longrightarrow
\langle
F_\mu,\nabla_\Theta C_{\tau,\mu}
\rangle_{L^2(\Theta;\mu)}
=0$. Therefore, $\alpha_n\longrightarrow0$, and
\begin{equation}
\|G_n\|_{L^2(\Theta;\mu_n)}
\longrightarrow
\|F_\mu\|_{L^2(\Theta;\mu)},
\qquad
\langle
G_n,\nabla_\Theta C_n
\rangle_{L^2(\Theta;\mu_n)}
=0.
\label{eq:Gn-properties}
\end{equation}

We next verify that the first-order decrease of the energy is
nonnegative. Set $A_n
:=
\int_\Theta\zeta|F_n|^2\,\mathrm d\mu_n$. By construction,
 $\langle F_n,G_n\rangle_{L^2(\Theta;\mu_n)}=
A_n-\frac{B_n^2}{D_n}$. The weighted Cauchy-Schwarz inequality gives $B_n^2\leq A_nD_n$,and hence
\begin{equation}
\langle F_n,G_n\rangle_{L^2(\Theta;\mu_n)}
\geq0.
\label{eq:positive-descent}
\end{equation}

The vector fields $G_n$ and $H_n$ are supported at positive distance
from $\partial K_r$. 

%For each fixed $n$, define
%$
%T_{n,s}(\theta):=\exp_\theta(-sG_n(\theta))$ and
%$\widetilde\mu_{n,s}:=(T_{n,s})_\#\mu_n$.
%Using Lemma~\ref{lem:C_first_variation_corrected} with the transport plan induced by $T_{n,s}$ and
%\eqref{eq:Gn-properties}, we obtain
%$C_\tau(\widetilde\mu_{n,s})=O_n(s^2)$.

For each fixed $n$, define
\[
 T_{n,s}(\theta):=\exp_\theta(-sG_n(\theta)),
 \qquad
 \widetilde\mu_{n,s}:=(T_{n,s})_\#\mu_n.
\]
Since $G_n$ is supported at positive distance from $\partial K_r$, this
map preserves $K_r$ for $s$ small.  Lemma~\ref{lem:C_first_variation_corrected}  and
$\langle G_n,\nabla_\Theta C_n\rangle_{L^2(\mu_n)}=0$ give
\[
 C_\tau(\widetilde\mu_{n,s})=O_n(s^2).
\]
For $s,t$ small, set
\[
 S_{n,s,t}(\theta):=\exp_\theta(tH_n(\theta)),
 \qquad
 \mu_{n,s,t}:=(S_{n,s,t})_\#\widetilde\mu_{n,s},
 \qquad
 G_n^{\mathrm{corr}}(s,t):=C_\tau(\mu_{n,s,t}).
\]
Because $G_n$ and $H_n$ are supported at positive distance from
$\partial K_r$, for each fixed $n$ there exists $\varepsilon_n>0$ such
that the maps $T_{n,s}$ and $S_{n,s,t}$ preserve $K_r$ whenever
$|s|+|t|<\varepsilon_n$. Writing the two push-forwards with respect to
the fixed measure $\mu_n$, exactly as in the proof of Proposition~\ref{prop:interior_slope_corrected},
shows that $G_n^{\mathrm{corr}}$ is of class $C^1$ near $(0,0)$. Moreover,
\[
  \partial_tG_n^{\mathrm{corr}}(0,0)
  =\int_\Theta \nabla_\Theta C_n\cdot H_n\,d\mu_n
  =D_n>0.
\]
The implicit function theorem therefore gives a function $t_n(s)$ such
that $t_n(0)=0$, $G_n^{\mathrm{corr}}(s,t_n(s))=0$, and
$t_n(s)=O_n(s^2)$.
Setting
$\mu_{n,s}:=\mu_{n,s,t_n(s)}$, we obtain $\mu_{n,s}\in X_{\tau,r}$ and
\[
 W_2(\mu_{n,s},\mu_n)
 \leq s\|G_n\|_{L^2(\mu_n)}+O_n(s^2),
\]
whereas Lemma~\ref{lem:E_first_variation_corrected} yields
\[
 {\cal E}_\tau(\mu_{n,s})-{\cal E}_\tau(\mu_n)
 =-s\langle\nabla_\Theta V_{\tau,\mu_n},G_n\rangle_{L^2(\mu_n)}
  +O_n(s^2).
\] 

%Applying the same implicit-function argument as in the proof of
%Proposition~\ref{prop:interior_slope_corrected}, with $H_n$ as correction field, gives a function
%$t_n(s)=O_n(s^2)$ and admissible measures
%$\mu_{n,s}\in X_{\tau,r}$ such that
%\[
%W_2(\mu_{n,s},\mu_n)
%\leq
%s\|G_n\|_{L^2(\Theta;\mu_n)}
%+
%O_n(s^2)\]
%and
%\[
%{\cal E}_\tau(\mu_{n,s})-{\cal E}_\tau(\mu_n)
%=
%-s
%\langle
%\nabla_\Theta V_{\tau,\mu_n},G_n
%\rangle_{L^2(\Theta;\mu_n)}
%+
%O_n(s^2).
%\]

Since $\langle
G_n,\nabla_\Theta C_n
\rangle_{L^2(\Theta;\mu_n)}
=0$, we have
\[
\langle
\nabla_\Theta V_{\tau,\mu_n},G_n
\rangle_{L^2(\Theta;\mu_n)}
=
\langle F_n,G_n\rangle_{L^2(\Theta;\mu_n)}.
\]
If $F_\mu=0$, the desired lower bound is immediate. Assume therefore
that $F_\mu\neq0$. Then $\|G_n\|_{L^2(\Theta;\mu_n)}>0$ for all
sufficiently large $n$. Using \eqref{eq:positive-descent} and letting
$s\downarrow0$, we obtain
\[
|\partial {\cal E}^C_{\tau,r}|(\mu_n)
\geq
\frac{
\langle F_n,G_n\rangle_{L^2(\Theta;\mu_n)}
}{
\|G_n\|_{L^2(\Theta;\mu_n)}
}.
\]
Finally,
\[
A_n\longrightarrow\|F_\mu\|_{L^2(\Theta;\mu)}^2,
\qquad
B_n\longrightarrow0,
\qquad
D_n\longrightarrow
\|\nabla_\Theta C_{\tau,\mu}\|_{L^2(\Theta;\mu)}^2,
\]
and consequently
\[
\langle F_n,G_n\rangle_{L^2(\Theta;\mu_n)}
\longrightarrow
\|F_\mu\|_{L^2(\Theta;\mu)}^2.
\]
Together with \eqref{eq:Gn-properties}, this yields
\[
\liminf_{n\to\infty}
|\partial {\cal E}^C_{\tau,r}|(\mu_n)
\geq
\|F_\mu\|_{L^2(\Theta;\mu)}.
\]
This proves the converse inequality and therefore~\eqref{eq:no_relaxation_corrected}.

\end{proof}

\begin{proof}[Proof of Proposition~\ref{prop:existenceCurveSupportCond}]
For $h>0$, set $\mu_h^0:=\mu_0$ and choose recursively
\begin{equation}
\label{eq:implicit-euler-truncated}
 \mu_h^n\in\operatorname*{argmin}_{\nu\in\mathcal P_2(K_r)}
 \left\{
 {\cal E}^C_{\tau,r}(\nu)+\frac1{2h}W_2^2(\nu,\mu_h^{n-1})
 \right\}.
\end{equation}
The metric space $(\mathcal P_2(K_r),W_2)$ is compact and complete.
Moreover, ${\cal E}^C_{\tau,r}$ is proper, lower semicontinuous, and bounded from
below, since $\mu_0\in X_{\tau,r}$ and
\[
{\cal E}_\tau(\mu)\geq-\frac12\|W\|_{L^\infty(\Omega)}
 \qquad\text{for every }\mu\in X_{\tau,r}.
\]
Thus every minimization problem \eqref{eq:implicit-euler-truncated} has a
solution.  The standard compactness theorem for minimizing movements gives
a sequence $h_k\downarrow0$ such that the piecewise constant interpolants
converge to a generalized minimizing movement
$\mu^r\in AC^2_{\mathrm{loc}}([0,+\infty);\mathcal P_2(K_r))$.
Since $X_{\tau,r}$ is closed, $\mu^r_t\in X_{\tau,r}$ for every $t\geq0$.

The discrete energy estimate, written with the De Giorgi variational
interpolation, passes to the limit and first gives the dissipation inequality
on intervals whose endpoints are limits of grid points.  Approximating
arbitrary $0\leq s\leq t$ by such endpoints, using the absolute continuity
of $\mu^r$ and the continuity of ${\cal E}_\tau$ on $\mathcal P_2(K_r)$, yields
\[
 {\cal E}_\tau(\mu^r_t)
 +\frac12\int_s^t|\dot\mu^r_q|^2\,dq
 +\frac12\int_s^t
   |\partial {\cal E}^C_{\tau,r}|_{\mathrm{rel}}^2(\mu^r_q)\,dq
 \leq {\cal E}_\tau(\mu^r_s).
\]
Here no upper-gradient property is used.  Indeed, every discrete state with
finite energy belongs to $X_{\tau,r}$, and on this closed set the finite part
of ${\cal E}^C_{\tau,r}$ is the continuous functional ${\cal E}_\tau$.  Hence the limiting
discrete energy coincides with ${\cal E}_\tau(\mu^r_t)$ at every time.  This is
precisely \eqref{eq:truncated-edi}; see, for instance,
\cite[Corollary~3.3.4 and the proof of Theorem~2.3.3, in particular
inequality~(3.4.1)]{AGS08}.

Finally, view $\mu^r$ as an absolutely continuous curve in
$\mathcal P_2(\Theta)$.  The velocity representation for Wasserstein
absolutely continuous curves yields, for a.e. $t>0$, a Borel tangent vector
field $v^r_t\in L^2(\mu^r_t;T\Theta)$ satisfying the continuity equation and
$
 \|v^r_t\|_{L^2(\mu^r_t)}=|\dot\mu^r_t|$.
The representation is applied on the ambient manifold $\Theta$; no smooth
manifold structure is asserted for the set $K_r$, which has boundary and
corners.
\end{proof}

\begin{comment}
\cve{Encore quelques petites choses à changer pour terminer la preuve de la Proposition 1}

\begin{proof}[Proof of Proposition~\ref{prop:existenceCurveSupportCond}]
The metric space \({\cal P}_2(K_r)\), endowed with \(W_2\), is compact and complete. By
Proposition~\ref{prop:lsc_compact_corrected}, \({\cal E}^C_{\tau,r}\) is proper and lower
semicontinuous on this compact metric space. The compactness theorem for generalized
minimizing movements therefore yields a generalized minimizing movement starting from
\(\mu_0\), and this curve is a curve of maximal slope for \({\cal E}^C_{\tau,r}\) with respect
to the relaxed local slope~\cite{}.

The last assertion follows from the standard velocity representation for absolutely
continuous curves in \({\cal P}_2(K_r)\) over a smooth finite-dimensional Riemannian manifold: for
a.e. \(t\), there is a tangent velocity field \(v_t^r\in L^2(\mu_t^r;T\Theta)\) satisfying
the continuity equation and realizing the metric derivative.
\end{proof}
\end{comment}

\subsection{Identification of the velocity in the interior}

\begin{lemma}[Chain rules along absolutely continuous curves]\label{lem:chain_rules_corrected}
Let $(\mu_t)_{t\in[0,T]}\subset \mathcal P_2(K_r)$ be absolutely
continuous, and let
$v_t\in L^2(\Theta,\mu_t;T\Theta)$ be a velocity field associated with
$(\mu_t)_{t\in[0,T]}$. Assume in addition that there exists $c>0$ such
that
\[
    \inf_{t\in[0,T]}\|u_{\mu_t}\|_{L^2(\Omega)}\geq c.
\]
Then $t\mapsto u_{\mu_t}$ is absolutely continuous with values in
$H^1(\Omega)$ and, for a.e. $t\in(0,T)$,
\begin{equation} \label{eq:chain_C_corrected}
\frac{d}{dt}{\cal C}_\tau(\mu_t)
 =\int_\Theta \nabla_\Theta C_{\tau,\mu_t}(\theta)
       \cdot v_t(\theta)\,d\mu_t(\theta),
\end{equation}
\begin{equation}  \label{eq:chain_E_corrected}
\frac{d}{dt}{\cal E}_\tau(\mu_t)
 =\int_\Theta \nabla_\Theta V_{\tau,\mu_t}(\theta)
       \cdot v_t(\theta)\,d\mu_t(\theta).
\end{equation}
In particular, if $\mu_t\in X_\tau$ for every $t\in[0,T]$, then
\begin{equation}\label{eq:constraint_compatibility_corrected}
\int_\Theta \nabla_\Theta C_{\tau,\mu_t}(\theta)
       \cdot v_t(\theta)\,d\mu_t(\theta)=0
\end{equation}
for a.e. $t\in(0,T)$.
\end{lemma}

\begin{proof}
Since the map
$K_r\ni\theta\longmapsto \Phi_\tau(\theta;\cdot)\in H^1(\Omega)$
is Lipschitz continuous, there exists $L_r>0$ such that
\[
\|u_{\mu_t}-u_{\mu_s}\|_{H^1(\Omega)}
 \leq L_r W_2(\mu_t,\mu_s)
 \leq L_r\int_s^t |\dot\mu_q|\,dq.
\]
Thus $t\mapsto u_{\mu_t}$ is absolutely continuous in $H^1(\Omega)$.
For a.e. $t\in(0,T)$, the Bochner integral
\[
  g_t:=\int_\Theta
       D_\Theta\Phi_\tau(\theta;\cdot)[v_t(\theta)]\,d\mu_t(\theta)
  \in H^1(\Omega)
\]
is well defined and satisfies
$
  \|g_t\|_{H^1(\Omega)}
  \leq L_r\|v_t\|_{L^1(\mu_t)}
  \leq L_r\|v_t\|_{L^2(\mu_t)}$.

Since $H^1(\Omega)$ is separable, choose a countable norm-dense family
$(\ell_k)_{k\geq1}$ in the unit ball of $(H^1(\Omega))^*$. Applying the
scalar differentiation formula for the continuity equation to the
observable $
  \Theta \ni \theta\longmapsto \ell_k\bigl(\Phi_\tau(\theta;\cdot)\bigr)$ gives, for each $k$ and for a.e. $t$,
\[
  \frac{d}{dt}\ell_k(u_{\mu_t})=\ell_k(g_t).
\]
Taking the intersection of the corresponding full-measure sets over
$k\geq1$, and using the absolute continuity of
$t\mapsto u_{\mu_t}$ in $H^1(\Omega)$, we obtain on a single
full-measure set
\[
  \ell_k(\dot u_{\mu_t}-g_t)=0
  \qquad\text{for every }k\geq1.
\]
Density then yields the identity in $H^1(\Omega)$
\[
  \dot u_{\mu_t}
  =\int_\Theta
       D_\Theta\Phi_\tau(\theta;\cdot)[v_t(\theta)]\,d\mu_t(\theta)
\]
for a.e. $t\in(0,T)$.
Since ${\cal E}_\tau(\mu_t)=\frac12 \mathfrak a(u_{\mu_t},u_{\mu_t})$, we obtain
\[
\begin{aligned}
\frac{d}{dt}{\cal E}_\tau(\mu_t)
 &=\mathfrak a(u_{\mu_t},\dot u_{\mu_t})\\
 &=\int_\Theta
   \mathfrak a\!\left(u_{\mu_t},
   D_\Theta\Phi_\tau(\theta;\cdot)[v_t(\theta)]\right)
   d\mu_t(\theta)\\
 &=\int_\Theta \nabla_\Theta V_{\tau,\mu_t}(\theta)
       \cdot v_t(\theta)\,d\mu_t(\theta),
\end{aligned}
\]
which proves~\eqref{eq:chain_E_corrected}. The lower bound on
$\|u_{\mu_t}\|_{L^2(\Omega)}$ allows us to differentiate the
$L^2$-norm, and therefore
\[
\begin{aligned}
\frac{d}{dt}{\cal C}_\tau(\mu_t)
 &=\frac{\langle u_{\mu_t},\dot u_{\mu_t}\rangle_{L^2(\Omega)}}
         {\|u_{\mu_t}\|_{L^2(\Omega)}}\\
 &=\int_\Theta \nabla_\Theta C_{\tau,\mu_t}(\theta)
       \cdot v_t(\theta)\,d\mu_t(\theta),
\end{aligned}
\]
which proves~\eqref{eq:chain_C_corrected}. Finally, if $\mu_t\in X_\tau$ for every
$t$, then ${\cal C}_\tau(\mu_t)=0$ for every $t$, and~\eqref{eq:constraint_compatibility_corrected} follows.
\end{proof}

\begin{comment}
\begin{lemma}[Chain rules along absolutely continuous curves]\label{lem:chain_rules_corrected}
Let \((\mu_t)_{t\in[0,T]}\subset {\cal P}_2(K_r)\) be absolutely continuous and let
\(v_t\in L^2(\Theta,\mu_t;T\Theta)\) be its velocity field. Then, for a.e. \(t\),
\begin{align}
  \frac{d}{dt}{\cal C}_\tau(\mu_t)
  &=\int_\Theta \nabla_\Theta C_{\tau,\mu_t}(\theta)
     \cdot v_t(\theta)\,d\mu_t(\theta),
  \label{eq:chain_C_corrected}\\
  \frac{d}{dt}{\cal E}_\tau(\mu_t)
  &=\int_\Theta \nabla_\Theta V_{\tau,\mu_t}(\theta)
     \cdot v_t(\theta)\,d\mu_t(\theta).
  \label{eq:chain_E_corrected}
\end{align}
In particular, if \(\mu_t\in X_\tau\) for every \(t\geq 0\), then
\begin{equation}\label{eq:constraint_compatibility_corrected}
  \int_\Theta \nabla_\Theta C_{\tau,\mu_t}(\theta)
  \cdot v_t(\theta)\,d\mu_t(\theta)=0
\end{equation}
for a.e. \(t\geq 0\).
\end{lemma}

\begin{proof}
The maps \({\cal E}_\tau\) and \({\cal C}_\tau\) are differentiable on \({\cal P}_2(K_r)\) in the sense of
Lemmas~\ref{lem:C_first_variation_corrected} and \ref{lem:E_first_variation_corrected},
with gradients bounded on \(K_r\). The Wasserstein chain rule for such regular
functionals gives \eqref{eq:chain_C_corrected} and \eqref{eq:chain_E_corrected}. If
\(\mu_t\in X_\tau\), then \({\cal C}_\tau(\mu_t)=0\) for all \(t\geq 0\), and
\eqref{eq:constraint_compatibility_corrected} follows.
\end{proof}
\end{comment}

\begin{proposition}[Interior characterization of the velocity]\label{prop:velocity_identification_corrected}
Let \((\mu_t^r)_{t\ge0}\) be the truncated curve given by
Proposition~\ref{prop:existenceCurveSupportCond}. Assume that for some
\(T>0\) and \(\delta>0\),
\begin{equation*}
  \operatorname{Supp}(\mu_t^r)\subset K_{r-\delta}
  \qquad\text{for every }t\in[0,T].
\end{equation*}
Then, for a.e. \(t\in[0,T]\),
\begin{equation}\label{eq:velocity_identification_corrected}
  v_t^r=-\nabla_\Theta V^C_{\tau,\mu_t^r}
  \qquad \mu_t^r\text{-a.e.}
\end{equation}
Consequently,
\begin{equation}\label{eq:energy_dissipation_corrected}
  \frac{d}{dt}{\cal E}_\tau(\mu_t^r)
  =-\int_\Theta |\nabla_\Theta V^C_{\tau,\mu_t^r}(\theta)|^2
  \,d\mu_t^r(\theta)
\end{equation}
for a.e. \(t\in[0,T]\).
\end{proposition}

\begin{proof}
%Fix a time \(t\) at which the metric derivative, the velocity field and the metric energy
%inequality are all defined. Set
%$F_t^r:=\nabla_\Theta V^C_{\tau,\mu_t^r}$.
%By Lemma~\ref{lem:chain_rules_corrected} and the constraint compatibility relation,
%\begin{equation}\label{eq:chain_with_F_corrected}
%  \frac{d}{dt}{\cal E}_\tau(\mu_t^r)
%  =\int_\Theta F_t^r(\theta)\cdot v_t^r(\theta)\,d\mu_t^r(\theta).
%\end{equation}
%Since the support of \(\mu_t^r\) is contained in \(K_{r-\delta}\),
%Propositions~\ref{prop:interior_slope_corrected} and \ref{prop:no_relaxation_corrected}
%give
%\begin{equation*}
%  |\partial {\cal E}^C_{\tau,r}|_{\rm rel}(\mu_t^r)
%  =\|F_t^r\|_{L^2(\mu_t^r)}.
%\end{equation*}
%The curve of maximal slope satisfies the pointwise energy dissipation inequality
%\begin{equation}\label{eq:metric_EDI_corrected}
%  \frac{d}{dt}{\cal E}_\tau(\mu_t^r)
%  \le
%  -\frac12|\dot\mu_t^r|^2
%  -\frac12|\partial {\cal E}^C_{\tau,r}|_{\rm rel}(\mu_t^r)^2.
%\end{equation}

Since $E_\tau$ is Lipschitz continuous on $\mathcal P_2(K_r)$ and
$\mu^r\in AC^2_{\mathrm{loc}}([0,+\infty);\mathcal P_2(K_r))$, the map
$t\mapsto E_\tau(\mu_t^r)$ is locally absolutely continuous. Fix a time
$t\in(0,T)$ at which the metric derivative, the velocity
representation, the chain rules of Lemma~\ref{lem:chain_rules_corrected} , and the Lebesgue
pointwise values used below are all defined. Set
\[
F_t^r:=\nabla_\Theta V^C_{\tau,\mu_t^r}.
\]
By Lemma~\ref{lem:chain_rules_corrected}  and the constraint compatibility relation,
\begin{equation}\label{eq:chain_with_F_corrected}
\frac{d}{dt}{\cal E}_\tau(\mu_t^r)
 =\int_\Theta F_t^r(\theta)\cdot v_t^r(\theta)\,d\mu_t^r(\theta).
\end{equation}
Since $\operatorname{Supp}(\mu_t^r)\subset K_{r-\delta}$,
Propositions~\ref{prop:interior_slope_corrected} and~\ref{prop:no_relaxation_corrected} give
\[
|\partial {\cal E}^C_{\tau,r}|_{\mathrm{rel}}(\mu_t^r)
 =\|F_t^r\|_{L^2(\Theta;\mu_t^r)}.
\]
On the other hand, applying the dissipation inequality~\eqref{eq:truncated-edi} on
$[t,t+h]$, dividing by $h>0$, and letting $h\downarrow0$, the Lebesgue
differentiation theorem yields
\begin{equation}\label{eq:metric_EDI_corrected}
\frac{d}{dt}{\cal E}_\tau(\mu_t^r)
 \leq -\frac12|\dot\mu_t^r|^2
       -\frac12|\partial {\cal E}^C_{\tau,r}|_{\mathrm{rel}}^2(\mu_t^r).
\end{equation}
Using also
$|\dot\mu_t^r|=\|v_t^r\|_{L^2(\Theta;\mu_t^r)}$, inequality~\eqref{eq:metric_EDI_corrected}
becomes
\[
\frac{d}{dt}{\cal E}_\tau(\mu_t^r)
 \leq -\frac12\|v_t^r\|_{L^2(\Theta;\mu_t^r)}^2
       -\frac12\|F_t^r\|_{L^2(\Theta;\mu_t^r)}^2.
\]
Combining~\eqref{eq:chain_with_F_corrected} with~\eqref{eq:metric_EDI_corrected}, we obtain
\[
\begin{aligned}
0
&\leq \frac12\int_\Theta
       |v_t^r(\theta)+F_t^r(\theta)|^2\,d\mu_t^r(\theta)\\
&=\frac12\|v_t^r\|_{L^2(\mu_t^r)}^2
  +\frac12\|F_t^r\|_{L^2(\mu_t^r)}^2
  +\int_\Theta F_t^r\cdot v_t^r\,d\mu_t^r\\
&=\frac12\|v_t^r\|_{L^2(\mu_t^r)}^2
  +\frac12\|F_t^r\|_{L^2(\mu_t^r)}^2
  +\frac{d}{dt}E_\tau(\mu_t^r)
\leq 0.
\end{aligned}
\]
Consequently,
\[
    v_t^r=-F_t^r
    \qquad \mu_t^r\text{-a.e.}
\]
for a.e. $t\in(0,T)$. Substituting this identity into~\eqref{eq:chain_with_F_corrected}
gives~\eqref{eq:energy_dissipation_corrected}.

%On the other hand, %\eqref{eq:chain_with_F_corrected} and Cauchy-Schwarz imply
%\begin{equation*}
%  \frac{d}{dt}{\cal E}_\tau(\mu_t^r)
%  \ge
%  -\frac12\|v_t^r\|_{L^2(\mu_t^r)}^2
%  -\frac12\|F_t^r\|_{L^2(\mu_t^r)}^2.
%\end{equation*}
%Therefore equality holds, and equality in Cauchy--Schwarz yields
%\begin{equation*}
%  v_t^r=-F_t^r
%  \qquad \mu_t^r\text{-a.e.}
%\end{equation*}
%The energy identity \eqref{eq:energy_dissipation_corrected} follows immediately from
%\eqref{eq:chain_with_F_corrected}.
\end{proof}

\subsection{Characteristic flow and stability estimates}

The aim of this section is to give the proof of Proposition~\ref{prop:propositionCharacteristics}. We begin with some preliminary lemma.

\begin{lemma}[Regularity of the constrained vector field]\label{lem:vector_field_regular_corrected}
Let \(R>0\). There exists a constant \(L_R>0\), depending only on
\(R,\tau,\|W\|_{L^\infty}\), such that for every \(\mu,\nu\in X_\tau\) with
\(\operatorname{Supp}(\mu)\cup\operatorname{Supp}(\nu)\subset K_R\), the vector fields
\begin{equation*}
  b_\mu:=-\nabla_\Theta V^C_{\tau,\mu},
  \qquad
  b_\nu:=-\nabla_\Theta V^C_{\tau,\nu}
\end{equation*}
satisfy, for every \(\theta,\theta_1,\theta_2\in K_R\),
\begin{align}
  |b_\mu(\theta)|&\le L_R,
  \label{eq:b_bounded_on_KR_corrected}\\
  |b_\mu(\theta_1)-b_\mu(\theta_2)|&\le L_R d(\theta_1,\theta_2),
  \label{eq:b_lip_theta_corrected}\\
  |b_\mu(\theta)-b_\nu(\theta)|&\le L_R W_2(\mu,\nu).
  \label{eq:b_lip_mu_corrected}
\end{align}
Moreover, if \(\operatorname{Supp}(\mu)\subset K_R\), then for every \(\theta\in\Theta\),
\begin{equation}\label{eq:b_linear_growth_corrected}
  |b_\mu(\theta)|\le L_R(1+|a|),
  \qquad \theta=(a,w,b).
\end{equation}
\end{lemma}

\begin{proof}
Write $\omega=(w,b)$ and set $\varphi_\omega(x):=\sigma_{H,\tau}(w\cdot x+b)$ for all $x\in \Omega$.
For $\mu\in X_\tau$, define
\[
c_\mu(\omega):=\langle u_\mu,\varphi_\omega\rangle_{L^2(\Omega)},
\qquad
q_\mu(\omega):=a(u_\mu,\varphi_\omega).
\]
Since $\|u_\mu\|_{L^2(\Omega)}=1$, one has
\[
C_{\tau,\mu}(a,w,b)=a\,c_\mu(w,b),
\qquad
V_{\tau,\mu}(a,w,b)=a\,q_\mu(w,b).
\]
Let $\psi_\mu:=q_\mu-\sigma_\mu c_\mu$.
Then
\begin{equation}\label{eq:49a}
 b_\mu(a,w,b)
 =\Bigl(-\psi_\mu(w,b),
        -aP_{w^\perp}\nabla_w\psi_\mu(w,b),
        -a\partial_b\psi_\mu(w,b)\Bigr).
\end{equation}

Assume that $\operatorname{Supp}(\mu)\subset K_R$. The boundedness of
the derivatives of $\sigma_{H,\tau}$ and the estimate
\[
\|u_\mu\|_{H^1(\Omega)}
 \leq \int_\Theta
       \|\Phi_\tau(\theta;\cdot)\|_{H^1(\Omega)}\,d\mu(\theta)
 \leq C(\tau,R)
\]
show that $c_\mu$ and $q_\mu$, together with their derivatives up to
order two in $(w,b)$, are bounded uniformly in $\mu$. Furthermore,
Lemma~\ref{lem:constraint_nondegenerate_corrected} gives
\[
D_\mu:=\|\nabla_\Theta C_{\tau,\mu}\|_{L^2(\mu)}^2
 \geq R^{-2}.
\]
The numerator
\[
N_\mu:=\int_\Theta
 \nabla_\Theta V_{\tau,\mu}\cdot
 \nabla_\Theta C_{\tau,\mu}\,d\mu
\]
is uniformly bounded. Hence $|\sigma_\mu|\leq C(\tau,R,\|W\|_{L^\infty})$.
Formula~\eqref{eq:49a} now yields~\eqref{eq:b_bounded_on_KR_corrected},~\eqref{eq:b_lip_theta_corrected}, and the
growth estimate~\eqref{eq:b_linear_growth_corrected}. In the proof of~\eqref{eq:b_lip_theta_corrected}, one also
uses the fact that the map $w\mapsto P_{w^\perp}=I-w\otimes w$ is
smooth on $\mathbb S^{d-1}$.

It remains to prove the dependence on the measure. Let
$\mu,\nu\in X_\tau$ have supports contained in $K_R$. If
$\gamma\in\Gamma_o(\mu,\nu)$, the differentiated feature estimates
give
\begin{equation}\label{eq:49b}
\|u_\mu-u_\nu\|_{H^1(\Omega)}
 \leq C(\tau,R)W_2(\mu,\nu).
\end{equation}
Consequently,
\begin{equation}\label{eq:49c}
\|c_\mu-c_\nu\|_{C^2}
 +\|q_\mu-q_\nu\|_{C^2}
 \leq C(\tau,R,\|W\|_{L^\infty})W_2(\mu,\nu),
\end{equation}
where the $C^2$-norm may be taken on
$\mathbb S^{d-1}\times[-R,R]$.

Set
\[
G_\rho:=\nabla_\Theta C_{\tau,\rho},
\qquad H_\rho:=\nabla_\Theta V_{\tau,\rho}
\qquad (\rho=\mu,\nu).
\]
The functions $H_\rho\cdot G_\rho$ and $|G_\rho|^2$ are uniformly
Lipschitz on $K_R$. Thus, using~\eqref{eq:49c} and
$W_1\leq W_2$, we obtain
\[
|N_\mu-N_\nu|+|D_\mu-D_\nu|
 \leq C(\tau,R,\|W\|_{L^\infty})W_2(\mu,\nu).
\]
Since $D_\mu,D_\nu\geq R^{-2}$, the quotient formula
$\sigma_\rho=N_\rho/D_\rho$ gives
\begin{equation}\label{eq:49d}
|\sigma_\mu-\sigma_\nu|
 \leq C(\tau,R,\|W\|_{L^\infty})W_2(\mu,\nu).
\end{equation}
Combining~\eqref{eq:49a},~\eqref{eq:49c}, and~\eqref{eq:49d} proves~\eqref{eq:b_lip_mu_corrected}.
\end{proof}\begin{remark}[Extension to a neighbourhood of the constraint set]
Let $R>0$ and $c>0$, and define
\[
\mathcal U_{R,c}:=
\left\{\mu\in\mathcal P_2(K_R):
       \|u_\mu\|_{L^2(\Omega)}\geq c\right\}.
\]
For every $\mu\in\mathcal U_{R,c}$, Lemma~\ref{lem:constraint_nondegenerate_corrected} gives
\[
\|\nabla_\Theta C_{\tau,\mu}\|_{L^2(\mu)}^2
 \geq \frac{c^2}{R^2}.
\]
Therefore $\sigma_\mu$ and $b_\mu$ are well defined on
$\mathcal U_{R,c}$. The proof of Lemma~\ref{lem:vector_field_regular_corrected}, with the factors
$\|u_\mu\|_{L^2}^{-1}$ kept in the definition of $c_\mu$, shows that
there exists $L_{R,c}>0$ such that, for all
$\mu,\nu\in\mathcal U_{R,c}$,
\[
\|b_\mu\|_{L^\infty(K_R)}
 +\operatorname{Lip}_{K_R}(b_\mu)\leq L_{R,c},
\qquad
\|b_\mu-b_\nu\|_{L^\infty(K_R)}
 \leq L_{R,c}W_2(\mu,\nu).
\]
Moreover, $|b_\mu(a,w,b)|\leq L_{R,c}(1+|a|)$ on $\Theta$.
\end{remark}

\medskip
\begin{lemma}[Geometric comparison estimate]\label{lem:geom}
Let $K\subset\Theta$ be compact, and let $B_t,\widetilde B_t$ be tangent
vector fields on $K$ such that, for a.e. $t\in[0,T]$,
\[
  \max\bigl\{\|B_t\|_{L^\infty(K)},\|\widetilde B_t\|_{L^\infty(K)}\bigr\}\leq M,
  \qquad
  |B_t(\theta)-B_t(\eta)|\leq Ld(\theta,\eta),
\]
and
\[
  \|B_t-\widetilde B_t\|_{L^\infty(K)}\leq\varepsilon(t),
  \qquad \varepsilon\in L^1(0,T).
\]
If $x_t,y_t\in K$ solve
\[
  \dot x_t=B_t(x_t),
  \qquad
  \dot y_t=\widetilde B_t(y_t),
\]
then there exists a constant $C=C(L,M)>0$ such that
\begin{equation}
\label{eq:geometric-comparison}
  d(x_t,y_t)
  \leq d(x_0,y_0)
  +C\int_0^t d(x_s,y_s)\,ds
  +\int_0^t\varepsilon(s)\,ds.
\end{equation}
\end{lemma}

\begin{proof}
Indeed, fix $r_*:=\pi/2$. Whenever $d(x_t,y_t)<r_*$, the minimizing
geodesic joining $y_t$ to $x_t$ is unique. If $P_{y_t\to x_t}$ denotes
parallel transport along this geodesic, the first-variation formula gives
\[
  D^+d(x_t,y_t)
  \leq
  \bigl|B_t(x_t)-P_{y_t\to x_t}\widetilde B_t(y_t)\bigr|.
\]
For the product manifold
$\Theta=\mathbb{R}\times\mathbb{S}^{d-1}\times\mathbb{R}$, parallel
transport acts trivially on the Euclidean factors. On the spherical
factor, if $w,z\in\mathbb{S}^{d-1}$ are non-antipodal and
$v\in T_z\mathbb{S}^{d-1}$, then
\[
  P_{z\to w}v
  =v-\frac{\langle v,w\rangle}{1+\langle w,z\rangle}(w+z),
  \qquad
  |P_{z\to w}v-v|
  \leq d_{\mathbb{S}^{d-1}}(w,z)|v|.
\]
Consequently,
\[
\begin{aligned}
  D^+d(x_t,y_t)
  &\leq
  |B_t(x_t)-P_{y_t\to x_t}B_t(y_t)|
  +|B_t(y_t)-\widetilde B_t(y_t)|\\
  &\leq (L+M)d(x_t,y_t)+\varepsilon(t).
\end{aligned}
\]
If $d(x_t,y_t)\geq r_*$, the metric derivative estimate gives
\[
  D^+d(x_t,y_t)
  \leq |B_t(x_t)|+|\widetilde B_t(y_t)|
  \leq 2M
  \leq \frac{2M}{r_*}d(x_t,y_t).
\]
Combining the two cases and integrating the resulting upper-Dini
inequality proves \eqref{eq:geometric-comparison}. In particular, if
$B_t=\widetilde B_t$, Gronwall's lemma yields
\[
  d(x_t,y_t)\leq e^{Ct}d(x_0,y_0).
\]
\end{proof}

We are now in position to prove Proposition~\ref{prop:propositionCharacteristics}. 

\begin{proof}[Proof of Proposition~\ref{prop:propositionCharacteristics}]
By Lemma~\ref{lem:vector_field_regular_corrected}, \(b_t^r\) is locally Lipschitz in
\(\theta\), uniformly in \(t\), and has at most linear growth. Since $t\mapsto\mu_t^r$ is absolutely continuous in $W_2$, it is
continuous. Lemma~\ref{lem:vector_field_regular_corrected} therefore implies that, for every compact
$K_R\subset\Theta$,
\[
\|b_t^r-b_s^r\|_{C^0(K_R)}
 \leq L_{\max\{R,r\}}W_2(\mu_t^r,\mu_s^r).
\]
In particular, $t\mapsto b_t^r$ is continuous with values in
$C^0(K_R)$, and hence is Borel measurable. Together with the local
Lipschitz estimate and the linear growth estimate of Lemma~\ref{lem:vector_field_regular_corrected}, this
verifies all the hypotheses of the Caratheodory
Cauchy-Lipschitz theorem. The latter theorem gives existence and
uniqueness of \eqref{eq:char_truncated_corrected}. 

Fix $R>0$. The linear-growth estimate provides a compact set
$K_{R'}$, depending only on $T,R,r$, which contains every trajectory
issued from $K_R$ during $[0,T]$. Applying
\eqref{eq:geometric-comparison} on $K_{R'}$ with
$B_q=\widetilde B_q=b_q^r$ and $\varepsilon=0$, and then applying
Gronwall's lemma, gives
\[
  d\bigl(\chi^r(t,s,\theta_1),\chi^r(t,s,\theta_2)\bigr)
  \leq e^{L_{T,R,r}|t-s|}d(\theta_1,\theta_2)
\]
for every $\theta_1,\theta_2\in K_R$ and every $s,t\in[0,T]$.

Set \(\nu_t:=\chi^r(t,0,\cdot)_\#\mu_0\). Then \((\nu_t)_{t\in[0,T]}\) solves
\begin{equation*}
  \partial_t\nu_t+\operatorname{div}(b_t^r\nu_t)=0,
  \qquad \nu_0=\mu_0.
\end{equation*}
By Proposition~\ref{prop:velocity_identification_corrected}, \((\mu_t^r)_{t\in[0,T]}\)
solves the same linear continuity equation. Uniqueness for the continuity equation with
a locally Lipschitz vector field gives \eqref{eq:pushforward_truncated_corrected}.
\end{proof}

\subsection{Construction and global well-posedness of the untruncated flow}

The aim of this section is to end the proof of Theorem~\ref{th:existenceNoSupport}. We begin with the following stability lemma. 

\begin{lemma}[Stability of bounded nonlinear solutions]\label{lem:nonlinear_stability_corrected}
Let $R>0$ and let
$(\mu_t)_{t\in[0,T]},(\nu_t)_{t\in[0,T]}\subset X_\tau$ be two
locally absolutely continuous solutions of
\[
\partial_t\rho_t+\operatorname{div}(b_{\rho_t}\rho_t)=0,
\qquad
b_\rho=-\nabla_\Theta V^C_{\tau,\rho},
\]
such that
\[
\operatorname{Supp}(\mu_t)\cup\operatorname{Supp}(\nu_t)
 \subset K_R
\qquad\text{for every }t\in[0,T].
\]
Then there exists $L_R>0$ such that
\begin{equation}
W_2(\mu_t,\nu_t)
 \leq e^{L_Rt}W_2(\mu_0,\nu_0)
\qquad\text{for every }t\in[0,T].
\tag{53}
\end{equation}
In particular, bounded-support solutions with the same initial datum
are unique.
\end{lemma}

\begin{proof}
By Lemma~\ref{lem:vector_field_regular_corrected}, the vector fields $b_{\mu_t}$ and $b_{\nu_t}$ are locally
Lipschitz in space, uniformly in time, and continuous in time on compact
sets. Uniqueness for the corresponding linear continuity equations gives
\[
  \mu_t=X_t^\mu{}_{\#}\mu_0,
  \qquad
  \nu_t=X_t^\nu{}_{\#}\nu_0,
\]
where $X^\mu$ and $X^\nu$ are the associated characteristic flows.

Let $\gamma_0\in\Gamma^o(\mu_0,\nu_0)$ and set
$\gamma_t:=(X_t^\mu,X_t^\nu)_\#\gamma_0$. For
$\gamma_0$-a.e. $(\theta,\eta)$, both trajectories remain in $K_R$.
Lemma~\ref{lem:vector_field_regular_corrected} and \eqref{eq:geometric-comparison}, applied with
\[
  B_s=b_{\mu_s},
  \qquad
  \widetilde B_s=b_{\nu_s},
  \qquad
  \varepsilon(s)=L_RW_2(\mu_s,\nu_s),
\]
therefore give, after enlarging $L_R$ if necessary,
\[
\begin{aligned}
  d\bigl(X_t^\mu(\theta),X_t^\nu(\eta)\bigr)
  &\leq d(\theta,\eta)\\
  &\quad+L_R\int_0^t
  \left[
    d\bigl(X_s^\mu(\theta),X_s^\nu(\eta)\bigr)
    +W_2(\mu_s,\nu_s)
  \right]ds.
\end{aligned}
\]
Define
\[
  D(t):=\left(
  \int_{\Theta^2}
  d\bigl(X_t^\mu(\theta),X_t^\nu(\eta)\bigr)^2
  \,d\gamma_0(\theta,\eta)
  \right)^{1/2}.
\]
Minkowski's inequality and the fact that $\gamma_t$ is a coupling of
$\mu_t$ and $\nu_t$ yield
\[
  D(t)
  \leq D(0)+L_R\int_0^t
  \bigl(D(s)+W_2(\mu_s,\nu_s)\bigr)\,ds
  \leq D(0)+2L_R\int_0^tD(s)\,ds.
\]
Hence
\[
  D(t)\leq e^{2L_Rt}D(0).
\]
Since $W_2(\mu_t,\nu_t)\leq D(t)$ and
$D(0)=W_2(\mu_0,\nu_0)$, the conclusion follows after renaming the
constant.
\end{proof}

We are now in a position to prove Theorem~\ref{th:existenceNoSupport}.

\begin{proof}[Proof of Theorem~\ref{th:existenceNoSupport}]
We first construct the solution on a short time interval by a fixed
point argument.

\medskip
\noindent\emph{Step 1: local construction.}
Choose $R_0>0$ such that
$\operatorname{Supp}(\mu_0)\subset K_{R_0}$, and fix $R>R_0+1$.
The feature estimates imply that
\[
\|u_\rho-u_{\mu_0}\|_{H^1(\Omega)}
 \leq C_R W_2(\rho,\mu_0)
 \qquad\text{for every }\rho\in\mathcal P_2(K_R).
\]
Choose $\varepsilon>0$ so small that
$C_R\varepsilon\leq \frac12$. For $T>0$, let
\[
\mathcal X_T:=
\left\{
\rho\in C([0,T];\mathcal P_2(K_R)):
\rho_0=\mu_0,
\ \sup_{t\in[0,T]}W_2(\rho_t,\mu_0)\leq\varepsilon
\right\},
\]
endowed with the complete distance
\[
 d_T(\rho,\widetilde\rho)
 :=\sup_{t\in[0,T]}W_2(\rho_t,\widetilde\rho_t).
\]
For every $\rho\in\mathcal X_T$ and every $t\in[0,T]$,
\[
\|u_{\rho_t}\|_{L^2(\Omega)}
 \geq \|u_{\mu_0}\|_{L^2(\Omega)}
      -\|u_{\rho_t}-u_{\mu_0}\|_{L^2(\Omega)}
 \geq \frac12.
\]
Hence the extension of Lemma~\ref{lem:vector_field_regular_corrected} applies uniformly to the curves in
$\mathcal X_T$. In particular, there exist constants $M,L>0$,
independent of $\rho\in\mathcal X_T$, such that on $K_R$
\[
|b_{\rho_t}|\leq M,
\qquad
|b_{\rho_t}(\theta)-b_{\rho_t}(\eta)|
 \leq Ld(\theta,\eta),
\]
and
\[
\|b_{\rho_t}-b_{\widetilde\rho_t}\|_{L^\infty(K_R)}
 \leq L W_2(\rho_t,\widetilde\rho_t).
\]
Moreover, $t\mapsto b_{\rho_t}$ is continuous on compact sets.
Let $X^\rho$ be the flow solving
\[
\partial_tX_t^\rho(\theta)=b_{\rho_t}(X_t^\rho(\theta)),
\qquad X_0^\rho(\theta)=\theta,
\]
and define
\[
\mathcal T(\rho)_t:=X_t^\rho{}_{\#}\mu_0.
\]
If $T$ is small enough that
\[
MT<\min\{R-R_0,\varepsilon\},
\]
then $X_t^\rho(\operatorname{Supp}\mu_0)\subset K_R$ and
\[
W_2(\mathcal T(\rho)_t,\mu_0)
 \leq
\left(\int_\Theta d(X_t^\rho(\theta),\theta)^2
                 \,d\mu_0(\theta)\right)^{1/2}
 \leq Mt.
\]
Thus $\mathcal T$ maps $\mathcal X_T$ into itself.

For $\rho,\widetilde\rho\in X_T$, apply
\eqref{eq:geometric-comparison} to the two characteristic systems, with
the same initial point $\theta\in\operatorname{Supp}(\mu_0)$. The uniform
bounds above and
\[
  \|b_{\rho_t}-b_{\widetilde\rho_t}\|_{L^\infty(K_R)}
  \leq L W_2(\rho_t,\widetilde\rho_t)
\]
give a constant $C_R>0$, independent of
$\rho,\widetilde\rho\in X_T$, such that
\[
\begin{aligned}
  d\bigl(X_t^\rho(\theta),X_t^{\widetilde\rho}(\theta)\bigr)
  &\leq C_R\int_0^t
  d\bigl(X_s^\rho(\theta),X_s^{\widetilde\rho}(\theta)\bigr)\,ds\\
  &\quad+L\int_0^tW_2(\rho_s,\widetilde\rho_s)\,ds.
\end{aligned}
\]
Gronwall's lemma therefore yields
\[
  d\bigl(X_t^\rho(\theta),X_t^{\widetilde\rho}(\theta)\bigr)
  \leq Lt e^{C_Rt}d_T(\rho,\widetilde\rho).
\]
Using the coupling
$(X_t^\rho,X_t^{\widetilde\rho})_\#\mu_0$, we obtain
\[
  d_T\bigl(\mathcal T(\rho),\mathcal T(\widetilde\rho)\bigr)
  \leq LT e^{C_RT}d_T(\rho,\widetilde\rho).
\]
Reducing $T$ if necessary, we may assume that
$LT e^{C_RT}<1$. The Banach fixed-point theorem then yields a unique
$\mu\in X_T$ such that
$\mu_t=X_t^\mu{}_{\#}\mu_0$.
In particular, $\mu$ is locally absolutely continuous and solves
\[
\partial_t\mu_t+\operatorname{div}(b_{\mu_t}\mu_t)=0.
\]

We now verify that the constraint is preserved. Lemma~\ref{lem:chain_rules_corrected} applies with
$c=\frac12$, and the definition of $\sigma_{\mu_t}$ gives
\[
\begin{aligned}
\frac{d}{dt}{\cal C}_\tau(\mu_t)
&=\int_\Theta \nabla_\Theta C_{\tau,\mu_t}
              \cdot b_{\mu_t}\,d\mu_t\\
&=-\int_\Theta \nabla_\Theta C_{\tau,\mu_t}
              \cdot\nabla_\Theta V_{\tau,\mu_t}\,d\mu_t
  +\sigma_{\mu_t}
   \int_\Theta|\nabla_\Theta C_{\tau,\mu_t}|^2\,d\mu_t
=0.
\end{aligned}
\]
Since ${\cal C}_\tau(\mu_0)=0$, we conclude that $\mu_t\in X_\tau$ on
$[0,T]$. The same computation and the orthogonality relation
\[
\int_\Theta \nabla_\Theta V^C_{\tau,\mu_t}
              \cdot\nabla_\Theta C_{\tau,\mu_t}\,d\mu_t=0
\]
yield the energy identity
\begin{equation}\label{eq:54a}
\frac{d}{dt}{\cal E}_\tau(\mu_t)
 =-\|\nabla_\Theta V^C_{\tau,\mu_t}\|_{L^2(\mu_t)}^2
\qquad\text{for a.e. }t.
\end{equation}

\medskip
\noindent\emph{Step 2: a priori estimates.}
Let $[0,T_*)$ be any interval on which the solution exists. From the
energy identity and the constraint, for every $t<T_*$,
\begin{equation}\label{eq:54}
-\frac12\|W\|_{L^\infty(\Omega)}
 \leq {\cal E}_\tau(\mu_t)\leq {\cal E}_\tau(\mu_0).
\end{equation}
Consequently, $u_{\mu_t}$ is bounded in $H^1(\Omega)$ uniformly on
bounded time intervals.

Set $A(t):=\int_\Theta a^2\,d\mu_t(a,w,b)$.
Writing $V^C_{\tau,\mu}(a,w,b)=a\psi_\mu(w,b)$ and using the
continuity equation, we obtain
\begin{equation}\label{eq:55}
\begin{aligned}
A'(t)
&=-2\int_\Theta a\,\partial_aV^C_{\tau,\mu_t}\,d\mu_t
 =-2\int_\Theta V^C_{\tau,\mu_t}\,d\mu_t\\
&=-4{\cal E}_\tau(\mu_t)+2\sigma_{\mu_t}
\end{aligned}
\qquad\text{for a.e. }t<T_*.
\end{equation}
Here we used
\[
\int_\Theta V_{\tau,\mu}\,d\mu=2{\cal E}_\tau(\mu),
\qquad
\int_\Theta C_{\tau,\mu}\,d\mu=1
\qquad (\mu\in X_\tau).
\]
Furthermore,
\[
|\sigma_\mu|
 \leq
\frac{\|\nabla_\Theta V_{\tau,\mu}\|_{L^2(\mu)}}
     {\|\nabla_\Theta C_{\tau,\mu}\|_{L^2(\mu)}}.
\]
The numerator is bounded by $C(1+A(\mu))^{1/2}$, whereas Lemma~\ref{lem:constraint_nondegenerate_corrected}
gives
$\|\nabla_\Theta C_{\tau,\mu}\|_{L^2(\mu)}^{-1}
 \leq A(\mu)^{1/2}$. Hence
\begin{equation}\label{eq:56}
|\sigma_\mu|\leq C\bigl(1+A(\mu)\bigr),
\end{equation}
where $C$ depends only on $\tau$, ${\cal E}_\tau(\mu_0)$, and
$\|W\|_{L^\infty}$. Combining \eqref{eq:54}--\eqref{eq:56} and applying
Gronwall's lemma, we obtain, for every finite $T<T_*$,
\begin{equation}\label{eq:57}
A(t)\leq A_T
\qquad\text{for every }t\in[0,T],
\end{equation}
where $A_T$ depends on the initial datum and on $T$, but not on any
truncation parameter.

The preceding estimates imply
\[
\sup_{t\in[0,T]}\|\psi_{\mu_t}\|_{C^1(\mathbb S^{d-1}\times\mathbb R)}
 \leq C_T.
\]
Along a characteristic $\theta(t)=(a(t),w(t),b(t))$,
\[
\dot a(t)=-\psi_{\mu_t}(w(t),b(t)),
\qquad
\dot b(t)=-a(t)\partial_b\psi_{\mu_t}(w(t),b(t)),
\]
while the equation for $w(t)$ is tangent to $\mathbb S^{d-1}$.
Therefore
\[
|\dot a(t)|\leq C_T,
\qquad
|\dot b(t)|\leq C_T|a(t)|.
\]
Since $\operatorname{Supp}(\mu_0)\subset K_{R_0}$, integration gives
a constant $R_T>0$ such that
\begin{equation}\label{eq:58}
\operatorname{Supp}(\mu_t)\subset K_{R_T}
\qquad\text{for every }t\in[0,T].
\end{equation}

\medskip
\noindent\emph{Step 3: global continuation.}
Let $T_{\max}$ be the maximal existence time of the solution obtained
above. Suppose by contradiction that $T_{\max}<+\infty$. Applying the
previous estimates with $T=T_{\max}$ shows that the supports remain in
a fixed compact set $K_R$ for all $t<T_{\max}$. On this compact set,
the vector field is uniformly bounded, and therefore
\[
|\dot\mu_t|
 \leq \|b_{\mu_t}\|_{L^2(\mu_t)}\leq C_R
\qquad\text{for a.e. }t<T_{\max}.
\]
Thus $\mu_t$ converges in $W_2$ as $t\uparrow T_{\max}$ to a measure
$\overline\mu\in\mathcal P_2(K_R)$. By continuity of $C_\tau$,
$\overline\mu\in X_\tau$. Applying the local construction with
initial datum $\overline\mu$ extends the solution beyond
$T_{\max}$, a contradiction. Hence $T_{\max}=+\infty$.

\medskip
\noindent\emph{Step 4: uniqueness, dissipation, and characteristics.}
On every bounded time interval, the solution has bounded support by~\eqref{eq:58}. Lemma~\ref{lem:nonlinear_stability_corrected} therefore gives uniqueness in the class of
locally absolutely continuous constrained solutions with locally
bounded support. Integrating~\eqref{eq:54a} and using the lower bound in~\eqref{eq:54} yields
\[
\int_0^T
\|\nabla_\Theta V^C_{\tau,\mu_t}\|_{L^2(\mu_t)}^2\,dt
 \leq {\cal E}_\tau(\mu_0)+\frac12\|W\|_{L^\infty(\Omega)}<+\infty.
\]
Finally, the construction gives
\[
\mu_t=\chi(t,0,\cdot)_{\#}\mu_0
\]
on every bounded time interval, where $\chi$ is the characteristic
flow of $-\nabla_\Theta V^C_{\tau,\mu_t}$. This proves Theorem~\ref{th:existenceNoSupport}.

\end{proof}

\section{Proof of convergence}\label{sec:proofConvergence}

In this section, we provide the proof of Theorem~\ref{theoremConvergence}. In the following, a LaSalle's principle argument is invoked in order to prove Theorem~\ref{theoremConvergence}. From now on, and in all this section, we assume that $d\geq 2$.

\subsection{Escape from non suitable critical points}

In this section, we use the notation $\omega := (w,b)$  so that:

$$
\theta = (a,w,b) = (a, \omega)
$$
to make the difference between "linear" variables and "nonlinear" ones. Moreover in order to ease the reading, we introduce the notation:

$$
v_\mu := \nabla_\Theta V^C_{\tau,\mu}. 
$$

%In this section, we use the homogeneity of $\Phi_\tau$ to prove that
%the gradient flow escape from bad critical points.
%
%\begin{definition}
%For $\mu \in \Prob_2(\Theta)$, one defines $h_1(\mu) \in \Prob_2(\Theta_1)$ with $\Theta_1 := \mathbb{R}^d \times \mathbb{R}$ as
%
%$$
%\forall \phi \in C_b(\Theta_1), \ \int_{\Theta_1} \phi(w,t) dh_1(\mu)(w,t) =  \int_{\Theta} c + a\phi(w,t) d\mu(c,a,w,t).
%$$
%\end{definition}

Let us first prove the following preliminary lemma.

\begin{lemma}
\label{lem:boundGradV}
Let \(\mu\in\mathcal P_2(\Theta)\), set \(u_\mu:=P_\tau\mu\), and let
\(\theta=(a,\omega) \in \Theta\). Then
\[
 |V_{\tau,\mu}(\theta)|+|\nabla_\Theta V_{\tau,\mu}(\theta)|
 \leq C\bigl(\tau,\|W\|_{L^\infty(\Omega)}\bigr)
 \left(\int_\Theta a'^2\,\mathrm d\mu(a',\omega')\right)^{1/2}
 (1+|a|).
\]
If \(u_\mu\neq0\), then \(C_{\tau,\mu}\) is well defined and
\[
 |C_{\tau,\mu}(\theta)|+|\nabla_\Theta C_{\tau,\mu}(\theta)|
 \leq C(\tau)(1+|a|).
\]
In particular, both estimates hold for every \(\mu\in X_\tau\).
\end{lemma}

\begin{proof}
For \(\omega=(w,b)\), set
\[
 \varphi_\omega(x):=\sigma_{H,\tau}(w\cdot x+b),
 \qquad
 \Phi_\tau(a,\omega;x)=a\varphi_\omega(x).
\]
The bounds on \(\sigma_{H,\tau}\) and its derivatives imply
\[
 \sup_{\omega}
 \Bigl(
   \|\varphi_\omega\|_{H^1(\Omega)}
   +\|D_\omega\varphi_\omega\|_{H^1(\Omega)}
 \Bigr)
 \leq C(\tau).
\]
Consequently,
\[
 \|u_\mu\|_{H^1(\Omega)}
 \leq C(\tau)\int_\Theta |a'|\,\mathrm d\mu(a',\omega)
 \leq C(\tau)
 \left(\int_\Theta a'^2\,\mathrm d\mu(a',\omega)\right)^{1/2}.
\]
The first estimate now follows directly from the definition of
\(V_{\tau,\mu}\), since differentiation with respect to \(a\) removes the
factor \(a\), whereas differentiation with respect to \(\omega\) leaves a
factor bounded by \(|a|\). The same argument applied to
\[
 C_{\tau,\mu}(a,\omega)
 =a\,\frac{\langle u_\mu,\varphi_\omega\rangle_{L^2(\Omega)}}
          {\|u_\mu\|_{L^2(\Omega)}}
\]
gives the second estimate.
\end{proof}

The velocity potential and the derived vector field are smooth in the sense give by Lemma \ref{lemmaContinuityPhi}.

\begin{lemma}[\(W_2\)-continuity of the reduced potential]\label{lemmaContinuityPhi}
Let \((\mu_n)_{n\geq1}\subset X_\tau\) and assume that
$
 \mu_n\longrightarrow\mu$ in $W_2(\Theta)$ for some \(\mu\in\mathcal P_2(\Theta)\). Then \(\mu\in X_\tau\) and
\[
 \sigma_{\mu_n}\longrightarrow\sigma_\mu,
 \qquad
 \|\psi_{\mu_n}-\psi_\mu\|_{C^1(B)}\longrightarrow0,
\]
where for all $\rho\in X_\tau$,
$ V^C_{\tau,\rho}(a,\omega)=a\psi_\rho(\omega)$.

Consequently, for every compact set \(K\subset\Theta\),
\[
 \sup_{\theta\in K}
 \left(
   |V^C_{\tau,\mu_n}(\theta)-V^C_{\tau,\mu}(\theta)|
   +|\nabla_\Theta V^C_{\tau,\mu_n}(\theta)
     -\nabla_\Theta V^C_{\tau,\mu}(\theta)|
 \right)
 \longrightarrow0.
\]
In particular, \(X_\tau\) is closed in \((\mathcal P_2(\Theta),W_2)\), and
\(\rho\mapsto\psi_\rho\) is continuous from \(X_\tau\) to \(C^1(B)\).
\end{lemma}

\begin{proof}
For \(\omega=(w,b)\), set
\[
 \varphi_\omega(x):=\sigma_{H,\tau}(w\cdot x+b),
 \qquad
 \Phi_\tau(a,\omega;x)=a\varphi_\omega(x),
 \qquad
 u_\rho:=P_\tau\rho.
\]
The boundedness of the derivatives of \(\sigma_{H,\tau}\) gives, for all
\(\theta=(a,\omega)\) and \(\widetilde\theta=(\widetilde a,\widetilde\omega)\),
\[
 \|\Phi_\tau(\theta;\cdot)-\Phi_\tau(\widetilde\theta;\cdot)\|_{H^1(\Omega)}
 \leq C(\tau)(1+|a|+|\widetilde a|)
 d(\theta,\widetilde\theta).
\]
Let \(\gamma_n\in\Gamma_o(\mu_n,\mu)\). By Minkowski's and
Cauchy--Schwarz inequalities,
\begin{align*}
 \|u_{\mu_n}-u_\mu\|_{H^1(\Omega)}
 &\leq C(\tau)
 \int_{\Theta^2}(1+|a|+|\widetilde a|)
 d(\theta,\widetilde\theta)\,\mathrm d\gamma_n \\
 &\leq C(\tau)
 \left(
   \int_{\Theta^2}(1+|a|+|\widetilde a|)^2\,\mathrm d\gamma_n
 \right)^{1/2}
 W_2(\mu_n,\mu).
\end{align*}
The first factor is bounded because \(W_2(\mu_n,\mu)\to0\). Hence $u_{\mu_n}\longrightarrow u_\mu$ in $H^1(\Omega)$. Since \(\|u_{\mu_n}\|_{L^2(\Omega)}=1\), it follows that
\(\|u_\mu\|_{L^2(\Omega)}=1\), and therefore \(\mu\in X_\tau\).

For \(\rho\in X_\tau\), define
\[
 c_\rho(\omega):=\langle u_\rho,\varphi_\omega\rangle_{L^2(\Omega)} \quad \mbox{ and } \quad
 q_\rho(\omega)
 :=\langle\nabla u_\rho,\nabla\varphi_\omega\rangle_{L^2(\Omega)}
   +\langle Wu_\rho,\varphi_\omega\rangle_{L^2(\Omega)}.
\]
Then
\[
 C_{\tau,\rho}(a,\omega)=a c_\rho(\omega),
 \qquad
 V_{\tau,\rho}(a,\omega)=a q_\rho(\omega),
 \qquad
 \psi_\rho=q_\rho-\sigma_\rho c_\rho.
\]
The functions \(\varphi_\omega\) and their first derivatives with respect to
\(\omega\), viewed as \(H^1(\Omega)\)-valued functions, are uniformly bounded
for \(\omega\in B\). Therefore
\[
 \|c_{\mu_n}-c_\mu\|_{C^1(B)}
 +\|q_{\mu_n}-q_\mu\|_{C^1(B)}
 \leq C\bigl(\tau,\|W\|_{L^\infty}\bigr)
 \|u_{\mu_n}-u_\mu\|_{H^1(\Omega)}
 \longrightarrow0.
\]
Moreover, by the support properties of \(\sigma_{H,\tau}\), these functions
and their first derivatives vanish whenever \(\omega\notin B\).

Set
\[
 D_\rho:=\|\nabla_\Theta C_{\tau,\rho}\|_{L^2(\Theta;\rho)}^2,
 \qquad
 N_\rho:=
 \langle\nabla_\Theta V_{\tau,\rho},
        \nabla_\Theta C_{\tau,\rho}\rangle_{L^2(\Theta;\rho)}.
\]
Since, for a function of the form \(f(a,\omega)=a r(\omega)\),
$
 \nabla_\Theta f(a,\omega)
 =\bigl(r(\omega),a\nabla_\omega r(\omega)\bigr)$, where $\nabla_\omega$ denotes the Riemannian gradient on $S^{d-1}\times \mathbb{R}$, one has
\[
 D_\rho
 =\int_\Theta
   \left(c_\rho(\omega)^2
   +a^2|\nabla_\omega c_\rho(\omega)|^2\right)\,\mathrm d\rho(a,\omega),
\]
and
\[
 N_\rho
 =\int_\Theta
   \left(q_\rho(\omega)c_\rho(\omega)
   +a^2\nabla_\omega q_\rho(\omega)\cdot
          \nabla_\omega c_\rho(\omega)\right)\,\mathrm d\rho(a,\omega).
\]
The integrands are continuous and have at most quadratic growth in \(a\),
uniformly for \(\rho=\mu_n\) and \(\rho=\mu\). The preceding
\(C^1(B)\)-convergence, together with \(W_2(\mu_n,\mu)\to0\), therefore yields
\[
 D_{\mu_n}\longrightarrow D_\mu,
 \qquad
 N_{\mu_n}\longrightarrow N_\mu.
\]
By Lemma~\ref{lem:constraint_nondegenerate_corrected}, \(D_\mu>0\). Consequently,
\[
 \sigma_{\mu_n}=\frac{N_{\mu_n}}{D_{\mu_n}}
 \longrightarrow
 \frac{N_\mu}{D_\mu}=\sigma_\mu.
\]
Combining this convergence with the convergence of \(c_{\mu_n}\) and
\(q_{\mu_n}\) proves that $\|\psi_{\mu_n}-\psi_\mu\|_{C^1(B)}\longrightarrow0$.
Finally,
\[
 V^C_{\tau,\rho}(a,\omega)=a\psi_\rho(\omega),
 \qquad
 \nabla_\Theta V^C_{\tau,\rho}(a,\omega)
 =\bigl(\psi_\rho(\omega),a\nabla_\omega\psi_\rho(\omega)\bigr),
\]
which gives the asserted locally uniform convergence on \(\Theta\).
\end{proof}

The main ingredient in our analysis is the following asymptotic sign separation property on transported slices. 

\begin{theorem}[Asymptotic sign separation on the transported slice]
\label{thm:asymptotic_sign_separation}
Assume Hypothesis~\ref{hyp:geometric_support}. Let $(\mu_t)_{t\ge0}$ be the solution
given by Theorem~\ref{th:existenceNoSupport}, and assume that
\[
\mu_t\longrightarrow \mu^\star
\qquad\text{in }W_2(\Theta)\quad\text{as }t\to+\infty.
\]
Write
\[
V_{\tau,\mu^\star}^C(a,\omega)=a\,\psi_{\mu^\star}(\omega),
\qquad \omega\in B.
\]

Let $\eta_->0$ and $\eta_+>0$ be such that $-\eta_-$ and $\eta_+$ are regular values of
$\psi_{\mu^\star}|_{\operatorname{int}(B)}$, and define
\[
O_-:=\{\omega\in B\; ;\; \psi_{\mu^\star}(\omega)<-\eta_-\},
\qquad
O_+:=\{\omega\in B\; ;\; \psi_{\mu^\star}(\omega)>\eta_+\}.
\]
Then there exists $T\ge0$ such that, for all $t\ge T$,
$$
\Sigma_t\cap(\mathbb R\times O_-)\subset [0,+\infty)\times O_-,
\qquad 
\mbox{ and } \qquad
\Sigma_t\cap(\mathbb R\times O_+)\subset (-\infty,0]\times O_+,
$$
where $\Sigma_t:=\chi(t,\{0\}\times B)$.
\end{theorem}

\begin{proof}
Set \(\psi_t:=\psi_{\mu_t}\) and \(\psi_\star:=\psi_{\mu^\star}\).
By Lemma~\ref{lemmaContinuityPhi},
$
 \|\psi_t-\psi_\star\|_{C^1(B)}\longrightarrow0$
as $t\to+\infty$. Since \(\psi_\star=0\) on \(\partial B\), the nonzero level sets considered
below are contained in \(\operatorname{int}(B)\).

We first treat the negative region; if \(O_-\) is empty, there is nothing to
prove. Since \(-\eta_-\) is a regular value of
\(\psi_\star|_{\operatorname{int}(B)}\),
$
 \partial O_-=\{\omega\in B:\psi_\star(\omega)=-\eta_-\}
$ is a compact \(C^1\) hypersurface contained in \(\operatorname{int}(B)\). Its
outward unit normal is
\[
 n_-(\omega)
 =\frac{\nabla_\omega\psi_\star(\omega)}
        {|\nabla_\omega\psi_\star(\omega)|}.
\]
Set
\[
 \beta_-:=\min_{\omega\in\partial O_-}
 |\nabla_\omega\psi_\star(\omega)|>0.
\]
The \(C^1(B)\)-convergence gives a time \(T_-\geq0\) such that, for all
\(t\geq T_-\),
\[
 \psi_t\leq-\frac{\eta_-}{2}
 \quad\text{on }\overline{O_-},
 \qquad
 \nabla_\omega\psi_t\cdot n_-
 \geq\frac{\beta_-}{2}
 \quad\text{on }\partial O_-.
\]

The positive region is treated in the same way. If \(O_+\neq\varnothing\), set
\[
 n_+(\omega)
 =-\frac{\nabla_\omega\psi_\star(\omega)}
         {|\nabla_\omega\psi_\star(\omega)|},
 \qquad
 \beta_+:=\min_{\omega\in\partial O_+}
 |\nabla_\omega\psi_\star(\omega)|>0.
\]
There exists \(T_+\geq0\) such that, for all \(t\geq T_+\),
\[
 \psi_t\geq\frac{\eta_+}{2}
 \quad\text{on }\overline{O_+},
 \qquad
 \nabla_\omega\psi_t\cdot n_+
 \leq-\frac{\beta_+}{2}
 \quad\text{on }\partial O_+.
\]
Here the last inequality is equivalent to
\(\nabla_\omega\psi_t\cdot\nabla_\omega\psi_\star>0\) on
\(\partial O_+\). If one of the two regions is empty, we simply take the
corresponding time equal to zero.

Set
\[
 T_0:=\max(T_-,T_+),
 \qquad
 M:=\max_{(a,\omega)\in\Sigma_{T_0}}|a|<+\infty,
\]
and
\[
 T:=T_0+\max\left(\frac{2M}{\eta_-},\frac{2M}{\eta_+}\right).
\]
Because \(t\mapsto\mu_t\) is continuous in \(W_2\), Lemma~\ref{lemmaContinuityPhi} implies that
\(t\mapsto\psi_t\) is continuous in \(C^1(B)\). Hence the vector field is
continuous along characteristics and the latter are \(C^1\).

We prove the negative inclusion. Suppose, by contradiction, that there exist
\(t\geq T\) and a characteristic
\[
 z(s)=(a(s),\omega(s))=\chi(s,(0,\omega_0))
\]
such that \(\omega(t)\in O_-\) and \(a(t)<0\). Let
\[
 I_-:=\{s\in[T_0,t]:\omega(s)\in O_-,\ a(s)<0\},
\]
and let \((\sigma,t]\) be the connected component of \(I_-\) containing
\(t\), with the convention \(\sigma=T_0\) if this component reaches the left
endpoint.

Assume that \(\sigma>T_0\). By continuity, either
$ a(\sigma)=0$ and $\omega(\sigma)\in\overline{O_-}$, or
$
 a(\sigma)<0$ and $\omega(\sigma)\in\partial O_-$.
The first alternative is impossible, since
\[
 \dot a(\sigma)=-\psi_\sigma(\omega(\sigma))
 \geq\frac{\eta_-}{2}>0,
\]
whereas \(a(s)<0\) for \(s>\sigma\) close to \(\sigma\). We are therefore in the second alternative. Set
$
 h_-(\omega):=\psi_\star(\omega)+\eta_-$.
Since \(\sigma\) is an entrance time into \(O_-=\{h_-<0\}\),
$
 \left.\frac{\mathrm d}{\mathrm ds}h_-(\omega(s))
 \right|_{s=\sigma+}\leq0$. On the other hand,
$$
 \left.\frac{\mathrm d}{\mathrm ds}h_-(\omega(s))
 \right|_{s=\sigma+}
 =-a(\sigma)\,
   \nabla_\omega\psi_\star(\omega(\sigma))\cdot
   \nabla_\omega\psi_\sigma(\omega(\sigma)) 
>0,
$$
where we used \(a(\sigma)<0\) and
\[
 \nabla_\omega\psi_\star\cdot\nabla_\omega\psi_\sigma
 =|\nabla_\omega\psi_\star|
  \bigl(n_-\cdot\nabla_\omega\psi_\sigma\bigr)>0
 \quad\text{on }\partial O_-.
\]
This contradiction proves that \(\sigma=T_0\).

Thus \(\omega(s)\in O_-\) and \(a(s)<0\) for every
\(s\in(T_0,t]\). Consequently,
\[
 \dot a(s)=-\psi_s(\omega(s))\geq\frac{\eta_-}{2}
 \quad\text{for a.e. }s\in[T_0,t],
\]
and hence
\[
 a(t)\geq a(T_0)+\frac{\eta_-}{2}(t-T_0)
 \geq-M+\frac{\eta_-}{2}(t-T_0)\geq0,
\]
which yields a contradiction with \(a(t)<0\). Therefore
\[
 \Sigma_t\cap(\mathbb R\times O_-)
 \subset[0,+\infty)\times O_-
 \qquad\text{for every }t\geq T.
\]
The proof of the positive inclusion is identical after reversing the signs.
\end{proof}

\medskip

We will also need the following proposition. 

\begin{proposition}[Finite-time escape from a Wasserstein neighbourhood]\label{prop:escape_signed_strip}
Let $\mu\in X_\tau$ and write
\[
  V^C_{\tau,\mu}(a,\omega)=a\,\psi_\mu(\omega),\qquad \omega\in B.
\]
Assume that $\psi_\mu$ takes a negative value, and let $-\eta<0$ be a
regular value in the range of $\psi_\mu|_{\operatorname{int}B}$. For
$\alpha>0$, define
\[
  O_-:=\{\omega\in B:\psi_\mu(\omega)<-\eta\},\qquad
  A^-_{\alpha,\eta}:=(\alpha,+\infty)\times O_-.
\]
Then there exists $\varepsilon_{\mu,\alpha,\eta}>0$ with the following
property. For every $t_0\geq0$ and every locally $W_2$-absolutely
continuous curve $(\nu_t)_{t\geq t_0}\subset X_\tau$ solving
\[
  \partial_t\nu_t+\operatorname{div}(b_t\nu_t)=0,
  \qquad b_t(\theta):=-\nabla_\Theta V^C_{\tau,\nu_t}(\theta),
\]
if
\[
  W_2(\nu_{t_0},\mu)\leq\varepsilon_{\mu,\alpha,\eta}
  \quad\text{and}\quad
  \nu_{t_0}(A^-_{\alpha,\eta})>0,
\]
then there exists $t_1>t_0$ such that
\[
  W_2(\nu_{t_1},\mu)>\varepsilon_{\mu,\alpha,\eta}.
\]

The symmetric statement holds when $\psi_\mu$ takes a positive value:
if $+\eta>0$ is a regular value in the range of
$\psi_\mu|_{\operatorname{int}B}$, set
\[
  O_+:=\{\omega\in B:\psi_\mu(\omega)>\eta\},\qquad
  A^+_{\alpha,\eta}:=(-\infty,-\alpha)\times O_+.
\]
\end{proposition}

\begin{proof}
We prove the negative case. By the support properties of
\(\sigma_{H,\tau}\), one has \(\psi_\mu=0\) on \(\partial B\). Hence the level set 
\(\{\psi_\mu=-\eta\}\subset\operatorname{int}(B)\). Set
\[
 h_-(\omega):=\psi_\mu(\omega)+\eta.
\]
Since \(0\) is a regular value of \(h_-|_{\operatorname{int}(B)}\) and
\(h_-^{-1}(0)\) is compact, there exist
\(\delta\in(0,\eta/2)\) and \(\beta_->0\) such that
\[
 |\nabla_\omega\psi_\mu(\omega)|\geq\beta_-
 \qquad\text{whenever }|h_-(\omega)|\leq\delta.
\]
Set $ U_-:=\{\omega\in B:|h_-(\omega)|\leq\delta\}$. This is a compact subset of \(\operatorname{int}(B)\).

By Lemma~\ref{lemmaContinuityPhi}, the map \(\rho\mapsto\psi_\rho\) is continuous from
\((X_\tau,W_2)\) to \(C^1(B)\). We may therefore choose
\(\varepsilon_{\mu,\alpha,\eta}>0\) such that, whenever
\(\rho\in X_\tau\) and
\(W_2(\rho,\mu)\leq\varepsilon_{\mu,\alpha,\eta}\),
\[
 \psi_\rho\leq-\frac{\eta}{2}
 \quad\text{on }\overline{O_-}
 \quad
\mbox{ and }
\quad \nabla_\omega\psi_\mu\cdot\nabla_\omega\psi_\rho
 \geq\frac{\beta_-^2}{2}
 \quad\text{on }U_-.
\]

Let a solution as in the statement be given. By inner regularity, there is a
compact set \(K\subset A^-_{\alpha,\eta}\) such that
$ m_0:=\nu_{t_0}(K)>0$. Set
$
 M_a:=\left(\int_\Theta a^2\,\mathrm d\mu\right)^{1/2}$, and choose \(T>t_0\) such that
\[
 \sqrt{m_0}\left(\alpha+\frac{\eta}{2}(T-t_0)\right)
 >M_a+\varepsilon_{\mu,\alpha,\eta}.
\]
Suppose, by contradiction, that
\[
 W_2(\nu_t,\mu)\leq\varepsilon_{\mu,\alpha,\eta}
 \qquad\text{for every }t\in[t_0,T].
\]
Then
\[
 \sup_{t\in[t_0,T]}\int_\Theta a^2\,\mathrm d\nu_t
 \leq\bigl(M_a+\varepsilon_{\mu,\alpha,\eta}\bigr)^2.
\]
Moreover,  Lemma~\ref{lemmaContinuityPhi} gives a uniform bound on
\(\|\psi_{\nu_t}\|_{C^1(B)}\). Since
\[
 b_t(a,\omega)
 =\bigl(-\psi_{\nu_t}(\omega),
        -a\nabla_\omega\psi_{\nu_t}(\omega)\bigr),
\]
we obtain
\[
 \int_{t_0}^T\int_\Theta |b_t(\theta)|^2
 \,\mathrm d\nu_t(\theta)\,\mathrm dt<+\infty.
\]
The map \((t,\theta)\mapsto b_t(\theta)\) is Borel, because
\(t\mapsto\nu_t\) is \(W_2\)-continuous and  Lemma~\ref{lemmaContinuityPhi} applies. The
superposition principle therefore yields a probability measure \(\Pi\) on
\(AC([t_0,T];\Theta)\) such that \(\nu_t=(e_t)_\#\Pi\) and, for
\(\Pi\)-a.e. trajectory \(\gamma(t)=(a(t),\omega(t))\),
\[
 \dot a(t)=-\psi_{\nu_t}(\omega(t)),
 \qquad
 \dot\omega(t)=-a(t)\nabla_\omega\psi_{\nu_t}(\omega(t)) \qquad \mbox{for a.e. }t\in[t_0,T].
\]

The set of trajectories satisfying \(\gamma(t_0)\in K\) has
\(\Pi\)-mass \(\nu_{t_0}(K) = m_0\). Consider one such trajectory. As long as
\(\omega(t)\in O_-\),
$\dot a(t)\geq\frac{\eta}{2}$, so in particular \(a(t)\geq\alpha\). Suppose that the trajectory leaves
\(O_-\) for the first time at \(\tau\leq T\). Then
\(h_-(\omega(\tau))=0\), and for \(t<\tau\) sufficiently close to \(\tau\)
one has \(\omega(t)\in U_-\). Therefore, for a.e. such \(t\),
\[
 \frac{\mathrm d}{\mathrm dt}h_-(\omega(t))
 =-a(t)\nabla_\omega\psi_\mu(\omega(t))\cdot
       \nabla_\omega\psi_{\nu_t}(\omega(t))
 \leq-\alpha\frac{\beta_-^2}{2}.
\]
Integrating up to \(\tau\) contradicts
\(h_-(\omega(t))<0\) for \(t<\tau\) and
\(h_-(\omega(\tau))=0\). Hence the trajectory never leaves \(O_-\), and
\[
 a(t)\geq\alpha+\frac{\eta}{2}(t-t_0)
 \qquad\text{for every }t\in[t_0,T].
\]
It follows that
\[
 \nu_T\left(
   \left\{(a,\omega):
   a\geq\alpha+\frac{\eta}{2}(T-t_0)\right\}
 \right)\geq m_0,
\]
and consequently
\[
 \int_\Theta a^2\,\mathrm d\nu_T
 \geq m_0\left(\alpha+\frac{\eta}{2}(T-t_0)\right)^2.
\]
On the other hand, since the projection \(\pi_a\) is \(1\)-Lipschitz,
\begin{align*}
 \left(\int_\Theta a^2\,\mathrm d\nu_T\right)^{1/2}
 &=W_2((\pi_a)_\#\nu_T,\delta_0) \\
 &\leq W_2((\pi_a)_\#\nu_T,(\pi_a)_\#\mu)
       +W_2((\pi_a)_\#\mu,\delta_0) \\
 &\leq\varepsilon_{\mu,\alpha,\eta}+M_a,
\end{align*}
contradicting the choice of \(T\). Thus the curve cannot remain in the
closed \(\varepsilon_{\mu,\alpha,\eta}\)-ball on the whole interval
\([t_0,T]\), and there exists \(t_1\in(t_0,T]\) such that
\[
 W_2(\nu_{t_1},\mu)>\varepsilon_{\mu,\alpha,\eta}.
\]
The positive case follows by reversing the signs.
\end{proof}

\begin{lemma}
\label{lemmaRegularValuesPhi}
For all $\mu \in X_\tau$, if $\psi_\mu < 0$ somewhere, there exists a strictly negative regular value $-\eta$ ($\eta>0$) of $\psi_\mu|_{\operatorname{int}(B)}$. Moreover, for all $\mu \in X_\tau$, if $\psi_\mu > 0$ somewhere, there exists a strictly positive regular value $\eta$  of $\psi_\mu|_{\operatorname{int}(B)}$. 
\end{lemma}

\begin{proof}
We first verify the regularity required in order to apply the
Morse-Sard theorem, which will be recalled below. For $\omega=(w,b)\in S^{d-1}\times\mathbb{R}$,
set $\phi_\omega(x):=\sigma_{H,\tau}(w\cdot x+b)$. Since $\sigma_{H,\tau}\in C_c^\infty(\mathbb{R})$ and $\Omega$ is
bounded, the map $ S^{d-1}\times\mathbb{R}\ni\omega
    \longmapsto \phi_\omega\in H^1(\Omega)$
is of class $C^\infty$. Indeed, in any smooth local chart of
$S^{d-1}$, every derivative with respect to $(w,b)$ is a finite
linear combination of functions of the form
\[
    x\longmapsto p(x)\,
    \sigma_{H,\tau}^{(j)}(w\cdot x+b),
\]
where $p$ is a polynomial of bounded degree; the same statement holds
after taking one derivative with respect to $x$. These expressions
depend continuously on $(w,b)$ with values in $H^1(\Omega)$.

Recall that, for $\mu\in X_\tau$,
\[
    c_\mu(\omega)
      =\langle u_\mu,\phi_\omega\rangle_{L^2(\Omega)},
    \qquad
    q_\mu(\omega)
      =a(u_\mu,\phi_\omega).
\]
The maps
\[
    v\longmapsto \langle u_\mu,v\rangle_{L^2(\Omega)}
    \quad\text{and}\quad
    v\longmapsto a(u_\mu,v)
\]
are continuous linear functionals on $H^1(\Omega)$. Consequently,
$c_\mu,q_\mu\in C^\infty(S^{d-1}\times\mathbb{R})$, and hence $
    \psi_\mu=q_\mu-\sigma_\mu c_\mu
    \in C^\infty(S^{d-1}\times\mathbb{R})$.

Next, we use the Morse-Sard theorem recalled below:

\begin{theorem}[Morse-Sard]
Let \(M\) be a smooth \(m\)-dimensional manifold without boundary, and let
\(f\in C^k(M;\mathbb R)\) with \(k\ge m\). Then the image of the critical points of $f$ (where the gradient is zero) 
has Lebesgue measure zero in \(\mathbb R\).

The proof of the positive case follows similar lines.  
\end{theorem}

In particular, $\psi_\mu|_{\operatorname{int}(B)}$ has the regularity
required by the Morse--Sard theorem, since
$
    \dim\bigl(\operatorname{int}(B)\bigr)
    =\dim(S^{d-1})+1=d$.

We now prove the negative case. Choose
$\omega_0=(w_0,b_0)\in B$ such that $\psi_\mu(\omega_0)<0$.
Since $\psi_\mu=0$ on $\partial B$, necessarily
$\omega_0\in\operatorname{int}(B)$. Set $
    L:=\sqrt d+2$
and consider the path
\[
    \gamma(s)
      :=\bigl(w_0,(1-s)b_0+sL\bigr),
    \qquad s\in[0,1].
\]
Then
\[
    \gamma([0,1))\subset\operatorname{int}(B),
    \qquad
    \gamma(1)\in\partial B,
    \qquad
    \psi_\mu(\gamma(1))=0.
\]
By the intermediate value theorem, every value in $\bigl(\psi_\mu(\omega_0),0\bigr)$
belongs to the range of
$\psi_\mu|_{\operatorname{int}(B)}$. By the Morse-Sard theorem, the set of critical values of
$\psi_\mu|_{\operatorname{int}(B)}$ has Lebesgue measure zero.
Since $\bigl(\psi_\mu(\omega_0),0\bigr)$ is a nonempty open interval,
we may choose a regular value $
    r\in\bigl(\psi_\mu(\omega_0),0\bigr)$.
Setting $\eta:=-r>0$ shows that $-\eta$ is a strictly negative
regular value in the range of
$\psi_\mu|_{\operatorname{int}(B)}$.

The positive case is identical: if $\psi_\mu(\omega_0)>0$, one
chooses a regular value $
    r\in\bigl(0,\psi_\mu(\omega_0)\bigr)$ 
and sets $\eta:=r$.
\end{proof}

\subsection{Convergence to a critical point}
%The proof of convergence is based on the following hypothesis on the initial measure $\mu_0$ : 
%
%\begin{hypothesis}\label{hypothesisSupport}
%The support of the measure $\mu_0$ verifies :
%$$
%\{0\} \times \{0\} \times S_{\R^d}(1) \times [-\sqrt{d} -2,\sqrt{d} + 2] \subset \text{supp}(\mu_0)
%$$
%\end{hypothesis}

We begin with two preliminary lemmas.

\begin{lemma}\label{lem:vanishing_outside_active_strip}
For every $\mu\in\mathcal X_\tau$ and every $\theta=(a,w,b)\in\Theta$ such that
\[
|b|\geq\sqrt d+2,
\]
one has
\[
V_{\tau,\mu}^C(\theta)=0,
\qquad
\nabla_\Theta V_{\tau,\mu}^C(\theta)=0.
\]
\end{lemma}

\begin{proof}
Let $x\in\Omega=[0,1]^d$. Since $|x|_\infty\le 1$ and $|w|_1\le \sqrt d$ for
$w\in S_{\mathbb R^d}(1)$,
\[
|w\cdot x+b|
\ge
|b|-|w|_1|x|_\infty
\geq
(\sqrt d+2)-\sqrt d
=
2.
\]
Because $\sigma_{H,\tau}$ and all its derivatives vanish outside $(-2,2)$, this implies $\sigma_{H,\tau}(w\cdot x+b)=0$ for all $x\in\Omega$.

Hence
\[
\Phi_\tau(a,w,b;\cdot)=0,
\qquad
\nabla_x\Phi_\tau(a,w,b;\cdot)=0,
\qquad
\nabla_\Theta\Phi_\tau(a,w,b;\cdot)=0.
\]
The conclusion follows directly from the definitions of
$V_{\tau,\mu}$, $C_{\tau,\mu}$ and $V_{\tau,\mu}^C$.
\end{proof}

In fact we can prove an important result on the support of the gradient curve $(\mu(t))_{t\geq 0}$ for all time. This is the result stated in the next lemma. 

\begin{lemma}[Propagation of the transported slice]\label{lem:qualitative_support_propagation}
Assume Hypothesis~\ref{hyp:geometric_support} and $d\geq 2$. For every $t\geq 0$, set
\[
    \Sigma_t:=\chi\bigl(t,\{0\}\times B\bigr).
\]
Then
\[
    \Sigma_t\subset\operatorname{Supp}(\mu_t)
    \qquad\text{and}\qquad
    \pi_{w,b}(\Sigma_t)=B.
\]
Consequently, for every nonempty relatively open set $O\subset B$,
\[
    \Sigma_t\cap(\mathbb{R}\times O)\neq\varnothing
    \qquad\text{and}\qquad
    \mu_t(\mathbb{R}\times O)>0.
\]
\end{lemma}

\begin{proof}
Set
\[
    L:=\sqrt d+2,
    \qquad
    B:=S^{d-1}\times[-L,L].
\]
Let $\chi(t,s,\theta)$ be the characteristic flow given by
Theorem~\ref{th:existenceNoSupport}, and write
\[
    \chi(t,\theta):=\chi(t,0,\theta).
\]
By Theorem~\ref{th:existenceNoSupport},
$
    \mu_t=\chi(t,\cdot)_{\#}\mu_0$. We first observe that $B$ is invariant under the projected flow
issued from the slice $\{0\}\times B$. By Lemma~\ref{lem:vanishing_outside_active_strip},
\[
    \nabla_\Theta V^C_{\tau,\mu_s}(a,w,b)=0
    \qquad\text{whenever } |b|\geq L.
\]
Thus, if a characteristic reaches $|b|=L$, its velocity vanishes
there. By uniqueness of characteristics, a trajectory starting with
$|b|\leq L$ cannot cross either component of the boundary
$|b|=L$.

Define
\[
    \xi_t:B\longrightarrow B,
    \qquad
    \xi_t(\omega)
      :=\pi_{w,b}\bigl(\chi(t,(0,\omega))\bigr).
\]
For fixed $t\geq0$, define also
\[
    H:[0,1]\times B\longrightarrow B,
    \qquad
    H(r,\omega)
      :=\pi_{w,b}\bigl(\chi(rt,(0,\omega))\bigr).
\]
Then
\[
    H(0,\omega)=\omega,
    \qquad
    H(1,\omega)=\xi_t(\omega).
\]
Moreover, if $\omega\in\partial B$, then the characteristic starting
from $(0,\omega)$ is constant, and therefore
\[
    H(r,\omega)=\omega
    \qquad
    \text{for every }r\in[0,1].
\]
Hence $\xi_t$ is homotopic to the identity relative to $\partial B$.

Through the homeomorphism
\[
    B\ni(w,b)
    \longmapsto
    \left(w,\frac{b+L}{2L}\right)
    \in S^{d-1}\times[0,1],
\]
the map $\xi_t$ satisfies the assumptions of Lemma~\ref{lemmaSurjectivity} (see below).
Consequently, $\xi_t$ is surjective, and thus
\[
    \pi_{w,b}(\Sigma_t)=\xi_t(B)=B.
\]

We next prove that
\[
    \Sigma_t\subset\operatorname{Supp}(\mu_t).
\]
We use the elementary fact that, if $F:X\to Y$ is continuous and
$x\in\operatorname{Supp}(\nu)$, then
\[
    F(x)\in\operatorname{Supp}(F_{\#}\nu).
\]
Indeed, if $U$ is an open neighborhood of $F(x)$, then
$F^{-1}(U)$ is an open neighborhood of $x$, so that
\[
    (F_{\#}\nu)(U)
      =\nu(F^{-1}(U))>0.
\]
Applying this fact to $F=\chi(t,\cdot)$ and using Hypothesis~\ref{hyp:geometric_support} gives
\[
    \Sigma_t
      =\chi\bigl(t,\{0\}\times B\bigr)
      \subset
      \chi\bigl(t,\operatorname{Supp}(\mu_0)\bigr)
      \subset\operatorname{Supp}(\mu_t).
\]

Finally, let $O\subset B$ be nonempty and relatively open. Such a set
meets the relative interior of $B$, so choose
\[
    \omega_1\in
    O\cap\bigl(S^{d-1}\times(-L,L)\bigr).
\]
By the surjectivity of $\xi_t$, there exists $\omega_0\in B$ such that
\[
    \xi_t(\omega_0)=\omega_1.
\]
Hence
\[
    z_1:=\chi(t,(0,\omega_0))
      \in\Sigma_t\cap(\mathbb{R}\times O),
\]
which proves the first conclusion.

Since $\omega_1\in\operatorname{int}(B)$ and $O$ is relatively open
in $B$, there exists an open neighborhood
$U\subset S^{d-1}\times\mathbb{R}$ such that
\[
    \omega_1\in U\subset O\cap\operatorname{int}(B).
\]
The set $\mathbb{R}\times U$ is therefore an open neighborhood of
$z_1$ in $\Theta$ and is contained in $\mathbb{R}\times O$.
Since $z_1\in\operatorname{Supp}(\mu_t)$,
\[
    \mu_t(\mathbb{R}\times O)
      \geq \mu_t(\mathbb{R}\times U)>0.
\]
\end{proof}

\medskip

The proof of the preceding lemma uses Lemma~\ref{lemmaSurjectivity} below, which gives conditions for the surjectivity of a continuous map on a cylinder.

\begin{lemma}\label{lemmaSurjectivity}[Surjectivity on the cylinder]
Let \(d\ge 2\), and set $C:=S^{d-1}\times[0,1]$.
Let \(f:C\to C\) be continuous. Assume that
\[
f(w,0)=(w,0),
\qquad
f(w,1)=(w,1)
\qquad
\forall w\in S^{d-1},
\]
and that \(f\) is homotopic to the identity relative to the boundary $\partial C=S^{d-1}\times\{0,1\}$.
Then \(f\) is surjective.
\end{lemma}

\begin{proof}
Let \(n:=d-1\). Since \(d\ge 2\), we have \(n\ge 1\).
We argue by contradiction and assume that \(f\) is not onto. Then there exists
$
p\in C\setminus f(C)$.
Since \(f\) fixes the boundary pointwise, no boundary point can be missed. Hence
$
p\in \operatorname{int}(C)=S^n\times(0,1)$.

Let $j:C\setminus\{p\}\hookrightarrow C$ 
be the canonical inclusion. Since \(p\notin f(C)\), the map \(f\) factors as $
f=j\circ g$, where
$
g:C\to C\setminus\{p\}
$ is the same map as \(f\), but with codomain restricted to \(C\setminus\{p\}\).
Moreover, because \(f\) fixes \(\partial C\) pointwise, \(g\) is a map of pairs
$
g:(C,\partial C)\to(C\setminus\{p\},\partial C)$. Similarly, \(j\) is a map of pairs
$
j:(C\setminus\{p\},\partial C)\to(C,\partial C)$.

Therefore, on relative homology,
$
f_\ast
=
j_\ast\circ g_\ast:
H_{n+1}(C,\partial C;\mathbb Z)
\to
H_{n+1}(C,\partial C;\mathbb Z)$.

Since \(C=S^n\times[0,1]\) is an oriented compact \((n+1)\)-manifold with boundary,
$
H_{n+1}(C,\partial C;\mathbb Z)\simeq \mathbb Z$.

Let \([C,\partial C]\) denote its relative fundamental class. Because \(f\) is homotopic
to the identity relative to \(\partial C\), the induced map on relative homology is the
identity:
$$
f_\ast [C,\partial C]=[C,\partial C]\neq 0.
$$

On the other hand,
\[
H_{n+1}(C\setminus\{p\},\partial C;\mathbb Z)=0.
\]
Indeed, \(C\setminus\{p\}\) has the homotopy type of an \(n\)-dimensional CW-complex
containing \(\partial C\), equivalently of a space obtained from the two boundary
spheres by adding \(n\)-dimensional cells. Hence there is no nontrivial relative
homology in degree \(n+1\).

Thus $g_\ast [C,\partial C]=0$, and consequently
\[
f_\ast [C,\partial C]
=
j_\ast g_\ast [C,\partial C]
=
0.
\]
This contradicts $f_\ast [C,\partial C]=[C,\partial C]\neq 0$. Therefore \(f\) must be surjective.
\end{proof}

\begin{comment}
\begin{lemma}\label{lemmaSurjectivity}
Let $f$ be a continuous map $f: S_{\R^d}(1) \times [0,1] \rightarrow S_{\R^d}(1) \times [0,1] =:C$, homotopic to the identity such that:

$$
\forall w \in S_{\R^d}(1), \ 
\left\{
\begin{aligned}
f(w,0) =& (w,0), \\
f(w,1) =& (w,1).
\end{aligned}
\right.
$$
Then $f$ is surjective.

\end{lemma}
\begin{proof}
Suppose that $f$ misses a point $p$, then necessarily $p = (w,t)$ with $0 < t < 1$. We can write:

$$
g: C \rightarrow C \setminus \{p\}
$$
the restriction of $f$ on its image. The induced homomorphism on homology groups writes:

$$
g_\star: H_{d-1}(C) \rightarrow H_{d-1}(C \setminus \{p\}).
$$

Aside that, we have the classic information on homology groups of $C$ and $C \setminus \{p\}$:

$$
\left\{
\begin{aligned}
H_{d-1}(C)&= H_{d-1}(S_{\R^d}(1)) &\simeq \mathbb{Z}, \\
H_{d-1}(C \setminus \{p\}) &= H_{d-1}(S_{\R^d}(1) \vee S_{\R^d}(1)) &\simeq \mathbb{Z}^2
\end{aligned}
\right.
$$
where $\vee$ designates the wedge sum. Thus, the homomorphism $g_\star$ can be written as: 
$$
g_\star: \mathbb{Z} \rightarrow \mathbb{Z}^2.
$$

As $g$ lets the two spheres $w \rightarrow (w,0), w \rightarrow (w,1)$ invariant, we have:

$$
g_\star(1) = (1,1).
$$
Now we note $i: C \setminus \{p\} \rightarrow C$ the canonical inclusion map. For all $(a,b)\in \mathbb{Z}^2$, 

$$
i_\star(a,b) = a + b.
$$

By hypothesis, $f$ is homotopic to the identity so $f_\star = I_\star$ and 
$f_\star(1) = 1$ but at the same time:

$$
f_\star(1) = i_\star g_\star (1) = i_\star((1,1)) = 2
$$
which gives a contradiction.

\end{proof}
\end{comment}

We then have the following corollary.

\begin{corollary}\label{cor:presence_in_signed_unstable_strip}
Assume Hypothesis~\ref{hyp:geometric_support}. Let $(\mu_t)_{t\ge0}$ be the solution
given by Theorem~\ref{th:existenceNoSupport}, and assume that
\[
\mu_t\to\mu^\star
\qquad\text{in }W_2(\Theta)\quad\text{as }t\to+\infty.
\]
Write
\[
V_{\tau,\mu^\star}^C(a,\omega)=a\,\psi_{\mu^\star}(\omega).
\]

Assume that $\psi_{\mu^\star}<0$ somewhere, let $\eta>0$ be such that $-\eta$ is a
regular value of $\psi_{\mu^\star}$, and define
\[
O_-:=\{\omega\in B\; ;\; \psi_{\mu^\star}(\omega)<-\eta\}.
\]
Then, for every $\alpha>0$, there exists $T_\alpha\ge0$ such that for all $t\ge T_\alpha$,
\[
\Sigma_t\cap\bigl((\alpha,+\infty)\times O_-\bigr)\neq\varnothing,
\]
where
\[
\Sigma_t:=\chi(t,\{0\}\times B).
\]
In particular,
\[
\mu_t\bigl((\alpha,+\infty)\times O_-\bigr)>0
\qquad\forall t\ge T_\alpha.
\]

The symmetric statement holds if $\psi_{\mu^\star}>0$ somewhere, using a regular value
$+\eta$, the set
\[
O_+:=\{\omega\in B\; ;\; \psi_{\mu^\star}(\omega)>\eta\},
\]
and the strip
\[
(-\infty,-\alpha)\times O_+.
\]
\end{corollary}

\begin{proof}
We only prove the negative case.

By Lemma~\ref{lem:qualitative_support_propagation}, for every $t\ge0$,
\[
\Sigma_t\cap(\mathbb R\times O_-)\neq\varnothing.
\]
By Theorem~\ref{thm:asymptotic_sign_separation}, there exists $T_0\ge0$ such that for all
$t\ge T_0$,
\[
\Sigma_t\cap(\mathbb R\times O_-)\subset [0,+\infty)\times O_-.
\]

By Lemma~\ref{lemmaContinuityPhi} and the convergence $\mu_t\to\mu^*$, after increasing $T_0$, we may assume, for every $t\ge T_0$,
$$
\psi_{\mu_t}\le -\eta/2 \mbox{ on }\overline O_-,
$$
and
$$
\nabla_\omega\psi_{\mu^*}\cdot\nabla_\omega\psi_{\mu_t}\ge c_->0,
$$   
on a fixed collar of $\partial O_-$. The barrier argument used in Proposition~\ref{prop:escape_signed_strip} then shows that $[0,+\infty)\times O_-$ is forward invariant and one has
\[
\psi_{\mu_t}(\omega)\le -\frac{\eta}{2},
\qquad\forall \omega\in\overline{O_-},\ \forall t\ge T_0.
\]

Fix $\alpha>0$ and set
\[
    \tau_\alpha:=\frac{4\alpha}{\eta},
    \qquad
    T_\alpha:=T_0+\tau_\alpha.
\]
Let $t\geq T_\alpha$ and define
\[
    s:=t-\tau_\alpha\geq T_0.
\]
By Lemma~\ref{lem:qualitative_support_propagation} and Theorem~\ref{thm:asymptotic_sign_separation}, we may choose
\[
    z_s\in
    \Sigma_s\cap\bigl([0,+\infty)\times O_-\bigr).
\]
Set
\[
    z_q:=\chi(q,s,z_s)
       =\bigl(a(q),\omega(q)\bigr),
    \qquad q\in[s,t].
\]
By the forward invariance of
$[0,+\infty)\times O_-$ established above,
\[
    z_q\in[0,+\infty)\times O_-
    \qquad\text{for every }q\in[s,t].
\]
Consequently,
\[
    \dot a(q)
      =-\psi_{\mu_q}(\omega(q))
      \geq\frac{\eta}{2}
\]
for almost every $q\in[s,t]$. Hence
\[
    a(t)
      \geq a(s)+\frac{\eta}{2}(t-s)
      \geq \frac{\eta}{2}\tau_\alpha
      =2\alpha>\alpha.
\]
Since $z_t=\chi(t,s,z_s)\in\Sigma_t$, we have proved that
\[
    \Sigma_t\cap\bigl((\alpha,+\infty)\times O_-\bigr)
    \neq\varnothing
    \qquad\text{for every }t\geq T_\alpha.
\]

Moreover, by Lemma~\ref{lem:qualitative_support_propagation}, $\Sigma_t\subset\operatorname{Supp}(\mu_t)$.

Since $O_-\Subset\operatorname{int}(B)$ is open,
$(\alpha,+\infty)\times O_-$ is open in $\Theta$. Therefore
\[
    \mu_t\bigl((\alpha,+\infty)\times O_-\bigr)>0
    \qquad\text{for every }t\geq T_\alpha.
\]

\begin{comment}
Fix $\alpha>0$ and set $
\tau_\alpha:=\frac{2\alpha}{\eta}$. Let $t\ge T_0+\tau_\alpha$ and define $s:=t-\tau_\alpha\ge T_0$.
Choose $z_s\in \Sigma_s\cap([0,+\infty)\times O_-)$ and let $z_t:=\chi(t,s,z_s)\in \Sigma_t$. Since the whole trajectory between $s$ and $t$ remains in $[0,+\infty)\times O_-$, its
$a$-component satisfies
\[
a(t)\ge a(s)+\frac{\eta}{2}(t-s)\ge \alpha.
\]
Therefore $z_t\in \Sigma_t\cap\bigl((\alpha,+\infty)\times O_-\bigr)$,which proves the first assertion.

Since $(\alpha,+\infty)\times O_-$is open and intersects the support of $\mu_t$, one has
\[
\mu_t\bigl((\alpha,+\infty)\times O_-\bigr)>0,
\qquad\forall t\ge T_0+\tau_\alpha.
\]
\end{comment}
\end{proof}

\medskip

We are now in a position to conclude on the convergence and prove Theorem~\ref{theoremConvergence}.

\begin{proof}[Proof of Theorem~\ref{theoremConvergence}]
Assume by contradiction that $V_{\tau,\mu^\star}^C\not\equiv 0$. Write
\[
V_{\tau,\mu^\star}^C(a,\omega)=a\,\psi_{\mu^\star}(\omega),
\qquad \omega\in B.
\]
Then either $\psi_{\mu^\star}<0$ somewhere or $\psi_{\mu^\star}>0$ somewhere.

\medskip

\noindent

We detail here the negative case, when $\psi_{\mu^\star}<0$ somewhere, the positive case being identical. Choose $\eta>0$ such that $-\eta$ is a regular value of $\psi_{\mu^\star}$, and define
\[
O_-:=\{\omega\in B\; ;\; \psi_{\mu^\star}(\omega)<-\eta\}.
\]
Fix $\alpha>0$.
By Corollary~\ref{cor:presence_in_signed_unstable_strip}, there exists $T_\alpha\ge0$ such
that
\[
\mu_t\bigl((\alpha,+\infty)\times O_-\bigr)>0,
\qquad\forall t\ge T_\alpha.
\]

Apply Proposition~\ref{prop:escape_signed_strip} with $\mu=\mu^*$ (using the fact that $X_\tau$ is closed in $W_2$ yields that $\mu_*\in X_\tau$) and let $\varepsilon:=\varepsilon_{\mu^*,\alpha,\eta}$. Since $\mu_t\to\mu^*$ in $W_2$, choose $t_0\ge T_\alpha$ such that
$W_2(\mu_t,\mu^*)\le\varepsilon$ for all $t\ge t_0$. Corollary~\ref{cor:presence_in_signed_unstable_strip} gives
$\mu_{t_0}((\alpha,+\infty)\times O_-)>0$. Proposition~\ref{prop:escape_signed_strip} then yields
$t_1>t_0$ with $W_2(\mu_{t_1},\mu^*)>\varepsilon$, a contradiction.

\medskip

Thus, the proof in the positive case being identical, neither $\psi_{\mu^\star}<0$ somewhere nor $\psi_{\mu^\star}>0$ somewhere is
possible. Therefore,
$
\psi_{\mu^\star}\equiv 0$ on $B$,
and hence
\[
V_{\tau,\mu^\star}^C(a,\omega)=a\,\psi_{\mu^\star}(\omega)\equiv 0
\qquad\text{on }\Theta.
\]
This concludes the proof.
\end{proof}

\section{Numerical simulations}\label{sec:sim}

In this section, we report proof-of-concept numerical experiments for a
finite-width stochastic optimization scheme motivated by the constrained
variational formulation introduced above. The purpose of these experiments is
to document the empirical behavior of this practical scheme on a collection of
smooth benchmark potentials. They are not intended to establish convergence,
ground-state selection, or scalability properties beyond the configurations
considered below. For all tests, we use the following parameters:

\begin{itemize}
\item We use the tensorflow/keras framework.
\item Two-layer neural networks are utilized with a network width of either $m = 100$ or $m = 1000$.
\item The dataset is made of $10^5$ points sampled uniformly from the domain $\Omega$.
\item Batches are made of $n:=100$ points taken from the dataset.
\item The energy $\E$, the constraint $\C$ and its derivatives are computed by Monte-Carlo approximation and automatic differentiation. 
\item The optimizer is the classical stochastic gradient descent (SGD) with learning rate equal to $\frac{1}{\tau m}$.
\item After each SGD update, the coefficients of the last linear layer are
rescaled so that the Monte Carlo approximation of
$\|u_\theta\|_{L^2(\Omega)}$ is equal to one. This rescaling can be viewed as a
practical projection onto the empirically approximated normalization constraint.
\item To evaluate the performance of our method, we use a finite difference algorithm that is capable of computing eigenfunctions and eigenvalues in the case where the potential depends only on the first two variables.
\item Note that because of Monte-Carlo sampling, the neural network algorithm is stochastic. This is why, we ran our algorithms 8 times to evaluate the mean and the variance of our results. 
\end{itemize}

\paragraph{Relation with the theoretical dynamics.}
The numerical procedure is motivated by the constrained Wasserstein dynamics
analyzed in the previous sections, but the two evolutions should not be
identified. The theoretical dynamics is an infinite-width, measure-valued flow
whose velocity is projected onto the tangent space of the constraint through
the nonlocal Lagrange multiplier $\sigma_\mu$. In contrast, the numerical
procedure uses a finite-width network, stochastic Monte Carlo approximations of
the energy and its derivatives, an explicit SGD update, and a subsequent
rescaling of the output layer. In particular, the rescaling step is not shown
here to coincide with a time discretization of the projected Wasserstein vector
field. We therefore regard the numerical algorithm as a related practical
optimization scheme inspired by the theoretical flow, rather than as a
consistent discretization established by the present analysis.

Establishing such a connection would require, among other ingredients, a
finite-particle or mean-field consistency analysis, together with a study of
the effects of the time step, the Monte Carlo approximation, and the
normalization procedure. These questions are beyond the scope of the present
work.

\medskip

Regarding the computation of the $L^2$ error, the exact solution is calculated via a finite difference scheme in low dimension ($1$ or $2$). The corresponding grid points and exact values are extracted from this resolution. Then, at each step, a set of elements from this low-dimensional grid of size equal to the batch size $n$ is randomly sampled. To evaluate the error in higher dimensions, we construct $d$-dimensional points by appending random coordinates drawn uniformly in $[0, 1]$ for the remaining dimensions, and the $L^2$ error is estimated on this resulting sample.

\medskip

In the first test, we take a potential such that the solution behaves well with respect to tensor trains method \cite{CoursMiSong} \textit{i.e.}:

$$
W(x) = 100 \cos(2 \pi x_1).
$$

\begin{figure}[H]
    \centering
    \begin{subfigure}[b]{0.45\textwidth}
        \centering
        \includegraphics[width=\textwidth]{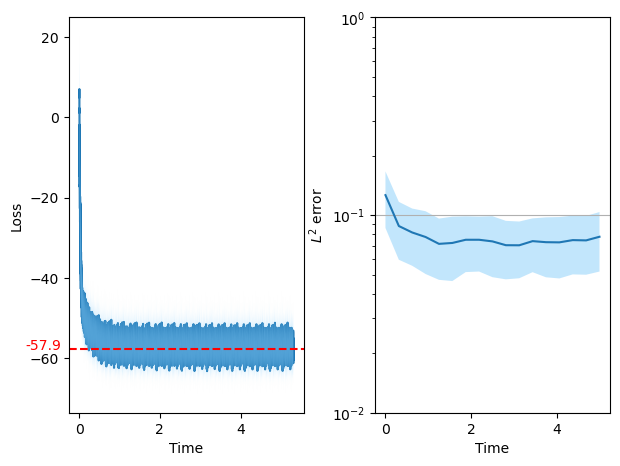}
        \caption{$m=100$}
    \end{subfigure}
    \hfill
    \begin{subfigure}[b]{0.45\textwidth}
        \centering
        \includegraphics[width=\textwidth]{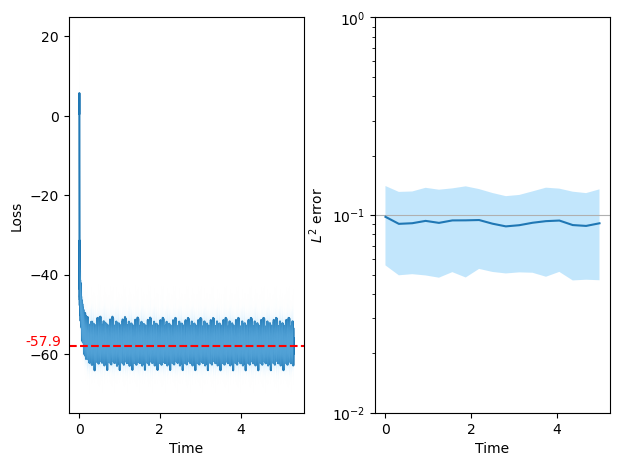}
        \caption{$m=1000$}
    \end{subfigure}
    \caption{Training loss and $L^2$ error for $W(x) = 100 \cos(2\pi x_1)$ and  $d=2$}    					\label{fig:cos_one_dir_d2}
\end{figure}

\begin{figure}[H]
    \centering
    \begin{subfigure}[b]{0.45\textwidth}
        \centering
        \includegraphics[width=\textwidth]{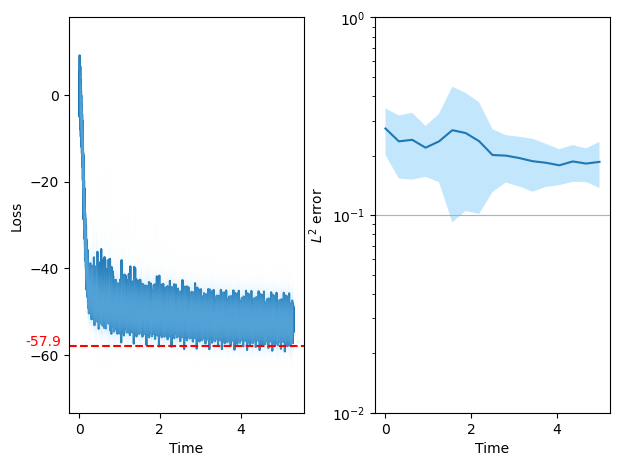}
        \caption{$m=100$}
    \end{subfigure}
    \hfill
    \begin{subfigure}[b]{0.45\textwidth}
        \centering
        \includegraphics[width=\textwidth]{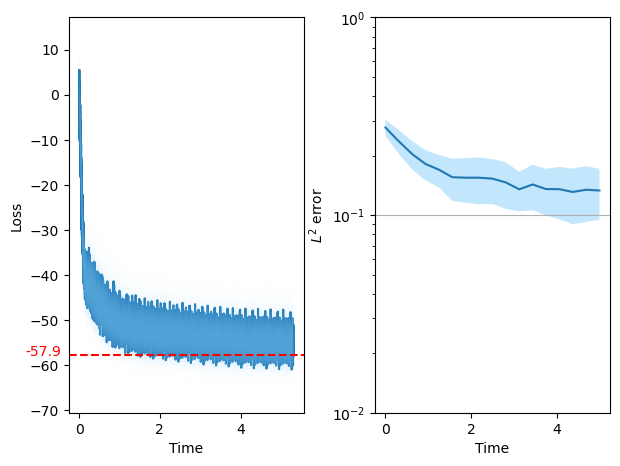}
        \caption{$m=1000$}
    \end{subfigure}
    \caption{Training loss and $L^2$ error for $W(x)= 100 \cos(2 \pi x_1)$ and $d=8$}    					\label{fig:cos_one_dir_d8}
\end{figure}

\begin{figure}[H]
    \centering
    \begin{subfigure}[b]{0.45\textwidth}
        \centering
        \includegraphics[width=\textwidth]{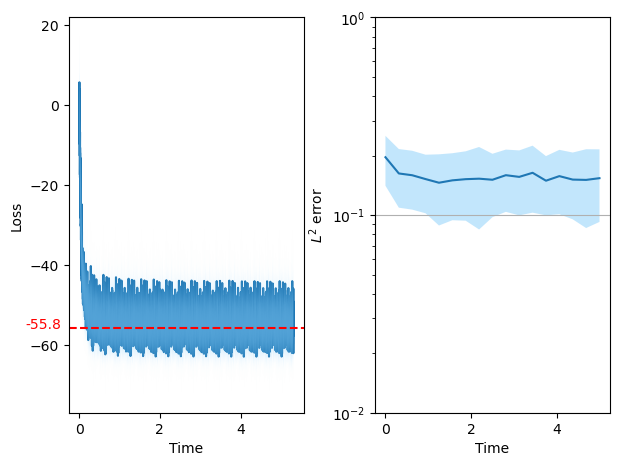}
        \caption{$d=2, m=100$}
    \end{subfigure}
    \hfill
    \begin{subfigure}[b]{0.45\textwidth}
        \centering
        \includegraphics[width=\textwidth]{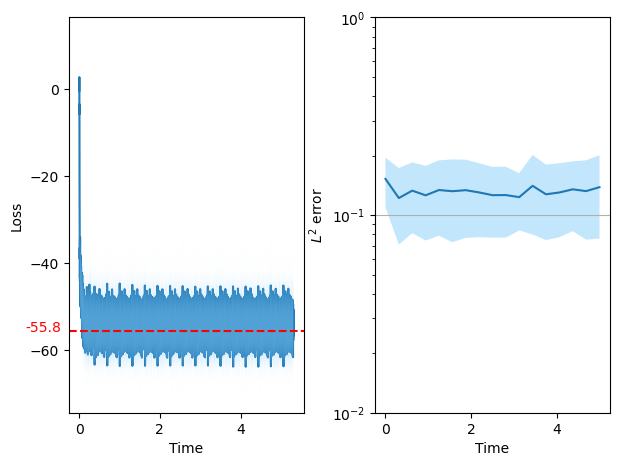}
        \caption{$d=2, m=1000$}
    \end{subfigure}
    \\
    \begin{subfigure}[b]{0.45\textwidth}
        \centering
        \includegraphics[width=\textwidth]{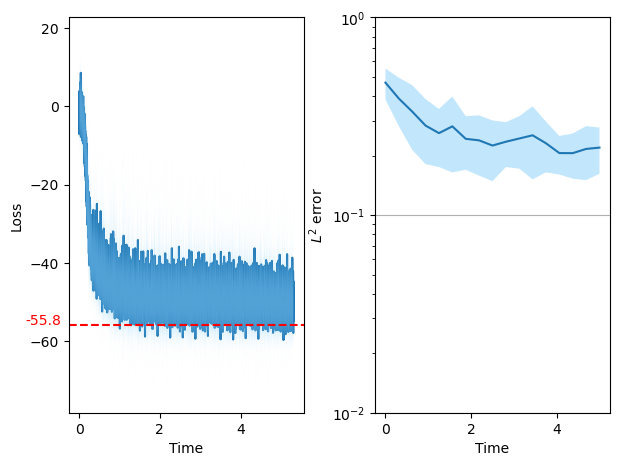}
        \caption{$d=8, m=100$}
    \end{subfigure}
    \hfill
    \begin{subfigure}[b]{0.45\textwidth}
        \centering
        \includegraphics[width=\textwidth]{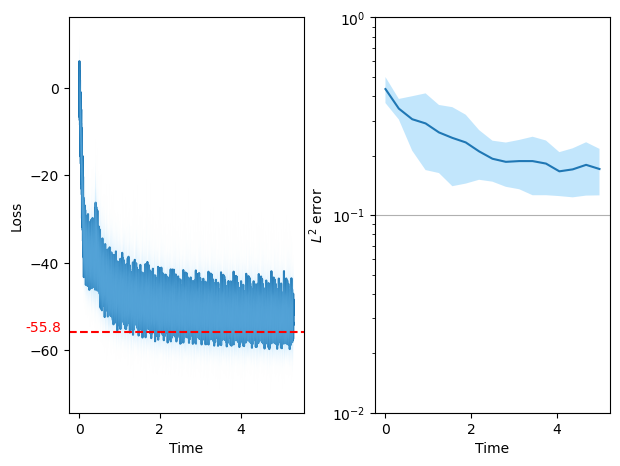}
        \caption{$d=8, m=1000$}
    \end{subfigure}
    \caption{Training loss and $L^2$ error for $W(x) = -100 \cos(2 \pi(x_1 - x_2))$}
    \label{fig:cos_diag}
\end{figure}

\begin{figure}[H]
    \centering
    \begin{subfigure}[b]{0.45\textwidth}
        \centering
        \includegraphics[width=\textwidth]{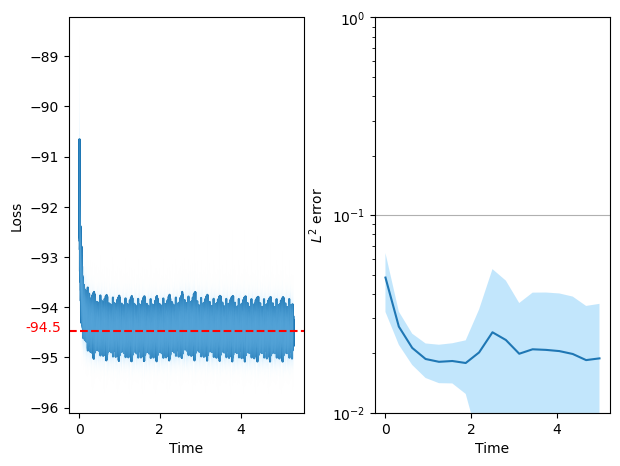}
        \caption{$d=2, m=100$}
    \end{subfigure}
    \hfill
    \begin{subfigure}[b]{0.45\textwidth}
        \centering
        \includegraphics[width=\textwidth]{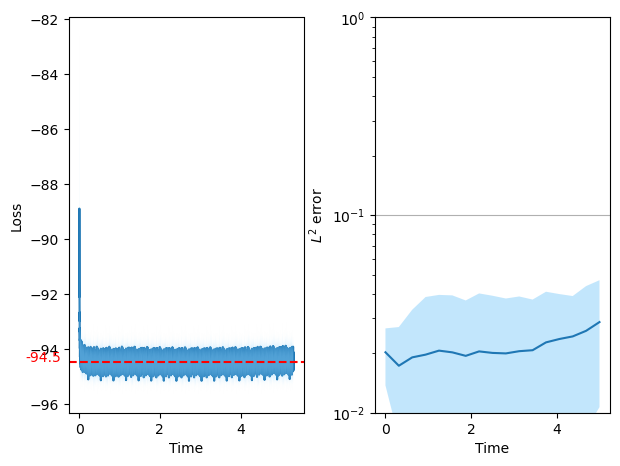}
        \caption{$d=2, m=1000$}
    \end{subfigure}
    \\
    \begin{subfigure}[b]{0.45\textwidth}
        \centering
        \includegraphics[width=\textwidth]{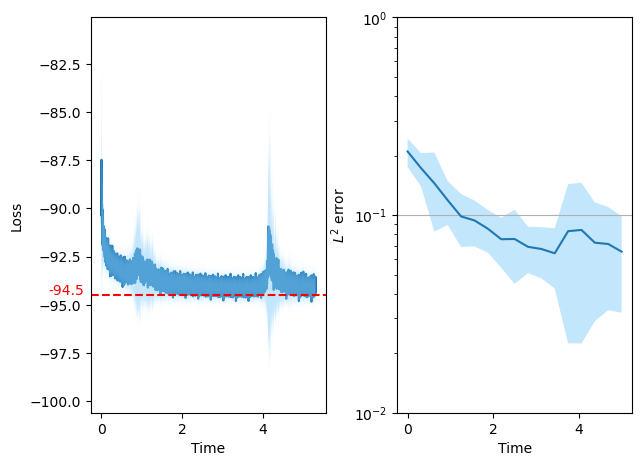}
        \caption{$d=8, m=100$}
    \end{subfigure}
    \hfill
    \begin{subfigure}[b]{0.45\textwidth}
        \centering
        \includegraphics[width=\textwidth]{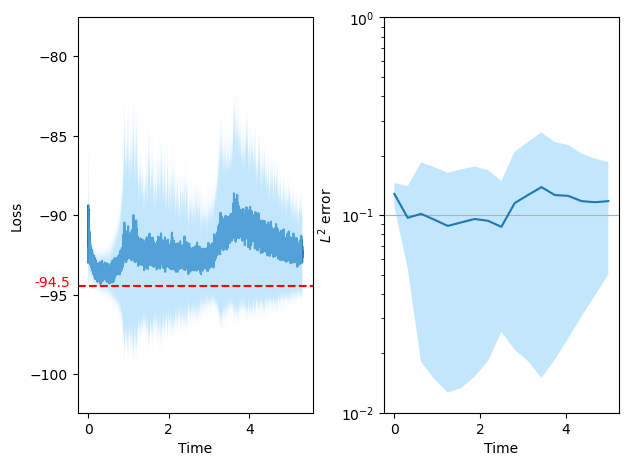}
        \caption{$d=8, m=1000$}
    \end{subfigure}
    \caption{Training loss and $L^2$ error for $W(x) = -100 e^{-\frac{1}{2} (x_1 - x_2)^2}$}
    \label{fig:exp_diag}
\end{figure}

Figures \ref{fig:cos_one_dir_d2}-\ref{fig:exp_diag} show the energy multiplied by a factor $2$ and $L^2$ error calculated thanks to the training dataset, during optimization. The red horizontal line corresponds to the first eigenvalue computed by the finite difference algorithm while the $L^2$ error is computed with respect to the first eigenfunction computed by this same algorithm. The thick blue line represents the mean value and the shaded area shows the variance of the quantity.

Next when one increases the dimension, curves are given in Figure \ref{fig:cos_one_dir_d8}. In all the experiments below, the potential depends either on $x_1$ alone or
on the pair $(x_1,x_2)$. Consequently, the associated ground-state problem has
an intrinsic spatial dimension of at most two: with Neumann boundary
conditions, the corresponding low-dimensional ground state can be extended
constantly in the remaining coordinates. The experiments with $d=8$ therefore
primarily assess the sensitivity of the training procedure to additional
inactive coordinates and to an enlarged ambient parametrization. They should
not be interpreted as tests on genuinely eight-dimensional coupled potentials,
nor as evidence that the method avoids the curse of dimensionality in general.

For both $d=2$ and $d=8$, the reported curves exhibit a similar qualitative
stabilization, and no systematic improvement is observed when the width is
increased from $m=100$ to $m=1000$. In view of the low intrinsic dimension of
this example, this observation is interpreted as robustness with respect to
additional inactive coordinates rather than as a general high-dimensional
scaling result.

\medskip
 
 Next we give a less obvious test where tensor techniques behave well. We call this test the "cos diagonal test".

$$
W(x) = -100 \cos(2 \pi (x_1 - x_2))
$$
Plots are given here for $d=2$ and $d=8$.

\medskip

Finally, we consider the following smooth nonseparable potential:
\[
    W(x)=-100\exp\left(-\frac{1}{2}(x_1-x_2)^2\right).
\]
We refer to this example as the ``exponential diagonal'' test. Unlike the
preceding one-dimensional example, this potential is not separable as a single
product of functions of $x_1$ and $x_2$. This benchmark is used to assess the
behavior of the neural-network optimization scheme in the presence of a
nonseparable dependence between two coordinates.

No direct numerical comparison with a tensor approximation method is performed
in the present work. Moreover, nonseparability by itself does not imply that a
function cannot be efficiently approximated by a low-rank tensor
representation. The results below should therefore not be interpreted as
evidence of a general advantage over tensor methods, but only as an evaluation
of the proposed neural-network scheme on this particular nonseparable
benchmark.

Plots are given in figure \ref{fig:exp_diag}. 

\medskip

\paragraph{Summary of the observations.}
Across the configurations considered in Figures~1--4, the estimated energy and
the reported $L^2$ error stabilize near the reference values obtained with the
finite-difference solver. On these particular examples, this behavior is
consistent with the numerical iterates approaching the reference ground state.
It does not, however, constitute a convergence result, nor does it imply that
the same optimization procedure must select the ground state for other
potentials, initializations, or hyperparameters. In particular, the theoretical
results of the previous sections identify convergent limits as eigenfunctions
under the stated assumptions, but do not establish convergence or ground-state
selection.

The comparison between $d=2$ and $d=8$ shows no pronounced degradation for the
examples considered here. As discussed above, these examples retain an
intrinsic dimension of at most two, so this observation is limited to
robustness with respect to additional inactive variables. The influence of the
network width is less systematic: increasing $m$ does not consistently reduce
the error or the variability over the tested range, and in some configurations
the wider network produces a larger reported error. Finally, since no tensor
solver is included in the comparison, the experiments do not support a
quantitative conclusion regarding the relative performance of neural-network
and tensor-based methods.

\medskip

\begin{figure}[H]
    \centering
    \includegraphics[width=0.6\textwidth]{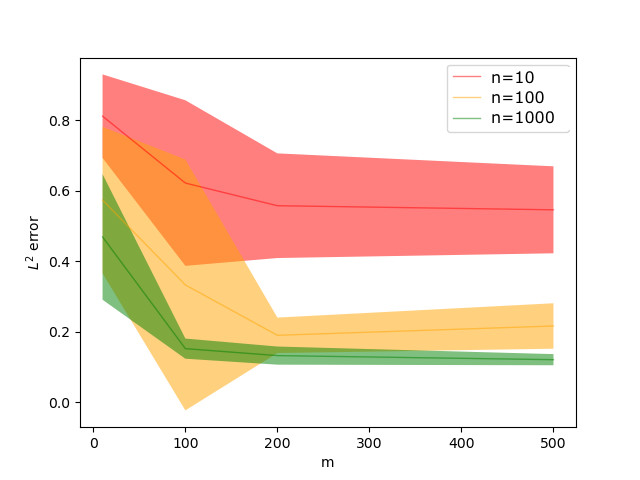}
    \caption{$L^2$ error for different batch sizes $n$ and widths $m$}
    \label{fig:effect_n_m}
\end{figure}

To compare more precisely the effect of increasing the width or the batch size, we compare the $L^2$ error after convergence for different parameters in Figure \ref{fig:effect_n_m} for the cos diagonal test. Within the range of parameters explored in Figure~5, increasing the network
width does not produce a monotone improvement of the final $L^2$ error. The
data instead suggest that, for this particular benchmark and optimization
setup, moderate widths are already sufficient to reach the observed accuracy
level. Increasing the batch size is associated with a lower empirical error and
a reduced variability across runs. These observations remain descriptive,
however: the limited number of configurations and independent runs does not
allow us to infer a general scaling law with respect to either the width or the
batch size.

\begin{figure}[H]
    \centering
    \includegraphics[width=0.6\textwidth]{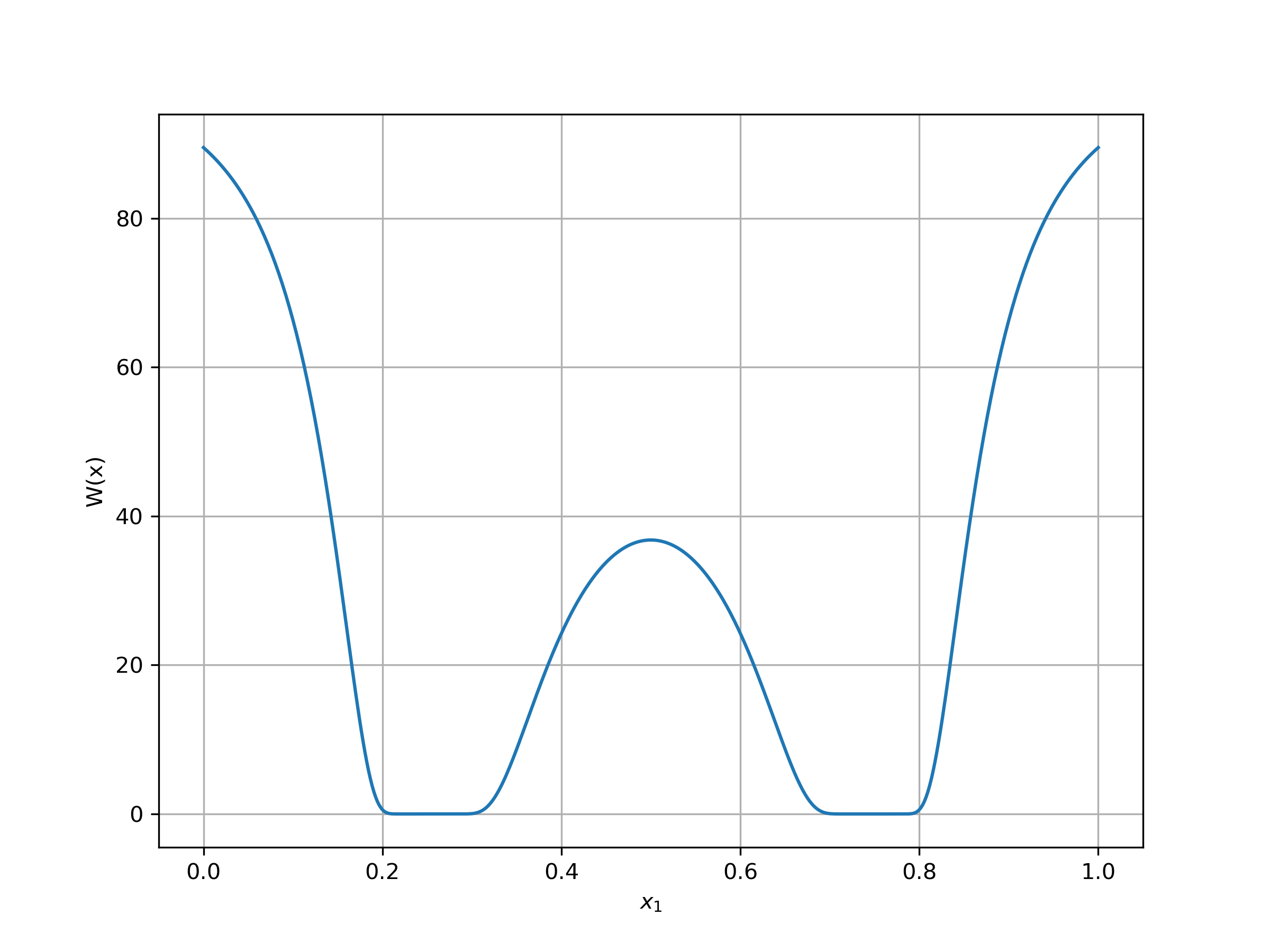}
    \caption{The double well potential}
    \label{fig:double_well_V}
\end{figure}

For the sake of completeness, we show a last case with a double well potential. First define $f(z) = (z^2 - 1)^{-2}$ and $g(z) = f(4(z-0.5))$, the potential $W$ is given by:

$$
W(x) = 100\exp(-g(x_1)).
$$
A visual representation of the potential is given in Figure \ref{fig:double_well_V}.

\begin{figure}[H]
    \centering
    \begin{subfigure}[b]{0.45\textwidth}
        \centering
        \includegraphics[width=\textwidth]{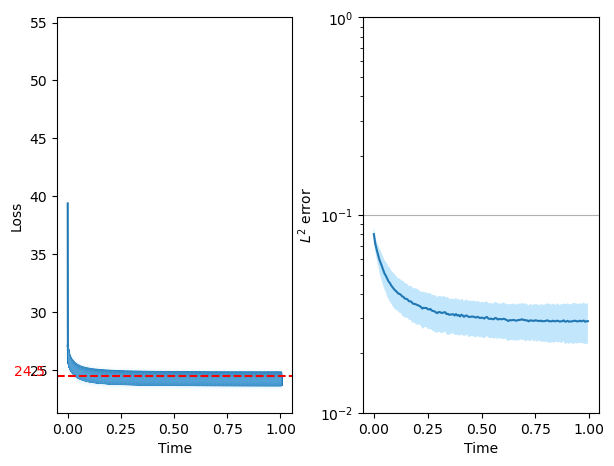}
        \caption{$d=2$}
    \end{subfigure}
    \hfill
    \begin{subfigure}[b]{0.45\textwidth}
        \centering
        \includegraphics[width=\textwidth]{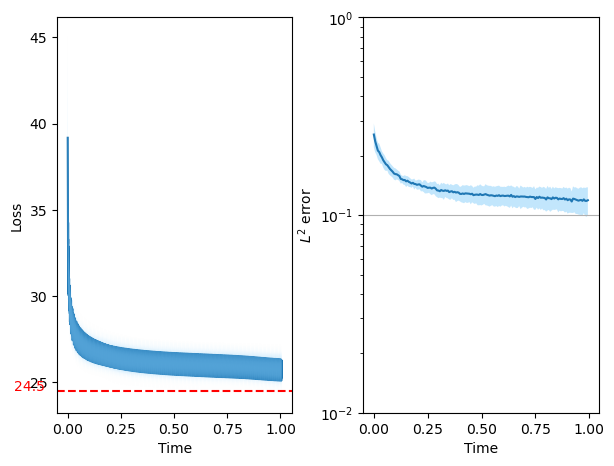}
        \caption{$d=8$}
    \end{subfigure}
    \caption{Training loss and $L^2$ error for $W(x) = 100 e^{-g(x_1)}$}    					
\label{fig:double_well}
\end{figure}

Owing to remarks related to Figure \ref{fig:effect_n_m}, we take a smaller network of size $m=200$ and a larger batch size with $n=1000$. For the selected values $m=200$ and $n=1000$, the energy and $L^2$-error curves
reported in Figure~\ref{fig:double_well} reach levels comparable to those observed in the preceding
benchmarks. This additional experiment suggests that the practical
optimization scheme can also handle this smooth double-well example with the
chosen hyperparameters. It should not be interpreted as a general robustness
result for multimodal, singular, or less regular potentials.

\begin{remark}[Scope of the numerical study]\label{rem:low_dimension}
The largest ambient dimension considered in these experiments is $d=8$, which
is modest compared with the dimensions encountered in many-body quantum
chemistry. In addition, all the benchmark potentials depend on at most two
coordinates, so their intrinsic spatial dimension remains low. The purpose of
these experiments is therefore to provide a proof of concept in a standard,
non-high-performance-computing environment and to illustrate several aspects
of the optimization procedure. They do not assess computational complexity for
genuinely $d$-dimensional coupled potentials, nor do they resolve the question
of scalability to the very-high-dimensional regime.
\end{remark}

%\begin{remark}\label{rem:low_dimension}
%It should be noted that in previous experiments, the dimension $d$ is relatively low, at most $8$, in comparison to applications in quantum chemistry where the dimension is three times the number of electrons of the atom or the molecule considered \cite{LeBris2006}. The examples displayed are intended to illustrate the capabilities of the method in a non-high-performance computing environment. This work does not address the fundamental question of scalability to the very high dimensional regime.
%\end{remark}

\section{Conclusion and perspectives}

In this work, an eigencouple of the Schrödinger operator is approximated by a two-layer neural network of infinite width, which can be represented by a probability measure. Subsequently, a Wasserstein gradient flow with respect to a constrained energy is introduced, and an existence result is provided. Assuming the convergence of the gradient flow, it is shown that the gradient flow converges to a suitable measure that represents an eigenfunction. The numerical experiments complement the theoretical analysis at a qualitative
level. For the selected smooth benchmark potentials and the reported
hyperparameters, the finite-width stochastic training procedure produces
states whose energies and represented functions are close to the
finite-difference ground-state references. In this restricted empirical sense,
the results support the practical relevance of the normalized variational
formulation. These observations should nevertheless be distinguished from the theoretical
results. The numerical algorithm, based on finite-width SGD and output-layer
rescaling, is inspired by but not proved to be a discretization of the
constrained Wasserstein flow analyzed in this work. Moreover, the theoretical
result assumes convergence of that flow and identifies a convergent limit as an
eigenfunction; it neither proves convergence nor guarantees selection of the
lowest eigenvalue. The apparent selection of the ground state in the reported
experiments is therefore an empirical observation rather than a consequence of
the present theory.
Establishing a rigorous connection between the finite-width stochastic
algorithm and the infinite-width constrained dynamics, as well as identifying
conditions that imply convergence and ground-state selection, constitute
natural directions for future work. Further numerical investigations should
also include genuinely high-dimensional coupled potentials, independent test
samples, and direct comparisons with alternative approximation methods at
matched accuracy and computational cost.

Let us point out that the convergence result we prove assumes the convergence itself, which represents a technical limitation. Finally, numerical tests were conducted on very smooth potentials for which the result from \cite{LuLu2022} holds, meaning that the solution can be approximated by two-layer neural networks. In order to address the Schrödinger multibody problem presented in the introduction, it is necessary to consider less regular potentials in a higher-dimensional context, as well as more complex neural networks. Furthermore, to be relevant for applications, the results obtained should be subsequently more accurate than those obtained in the present framework. One potential avenue for improvement is the use of the natural gradient technique, as introduced in \cite{MullerZeinhofer2023} and the references therein. We intend to investigate this path in a future research work.

\section*{Statements and Declarations}

\paragraph*{Funding}
The research leading to these results received funding from the European Research Council (ERC) under the European Union's Horizon 2020 Research and Innovation Program (Grant Agreement n$^\circ$101077204 High-LEAP).

\paragraph*{Competing Interests}
The authors have no relevant financial or non-financial interests to disclose.

\paragraph*{Data Availability}
The datasets and code generated during and/or analysed during the current study are available from the corresponding author on reasonable request.

\bibliographystyle{unsrt}
\bibliography{sample}

\newpage

\end{document}